\documentclass[invmat,final,numbook,francais]{svjour}
\usepackage{times}
\usepackage[T1]{fontenc}
\usepackage{amsmath}
\usepackage{amssymb}
\usepackage{amsfonts}
\usepackage{eucal}
\usepackage[all]{xy}
\usepackage[frenchb]{babel}

\newtheorem{prop}{Proposition}[subsection]
\newtheorem{theo}[prop]{Théor\`eme}
\newtheorem{coro}[prop]{Corollaire}
\newtheorem{lemm}[prop]{Lemme}
\newtheorem{lemm*}[equation]{Lemme}

\spnewtheorem{conj}[prop]{Conjecture}{\bf}{\rm}
\spnewtheorem{vide}[prop]{}{\bf}{\rm}
\spnewtheorem{defi}[prop]{Définition}{\bf}{\rm}

\newtheorem{rema}[prop]{Remarques}

\newtheorem{nota}[prop]{Notations}

\numberwithin{equation}{prop}

\newcommand{\riso}{ \tilde {\rightarrow}\, }

\newcommand{\Spec}{\mathrm{Spec}\,}
\newcommand{\Spf}{\mathrm{Spf}\,}
\newcommand{\Spff}{\mathrm{Spff}\,}
\renewcommand{\sp}{\mathrm{sp}}
\renewcommand{\det}{\mathrm{det}}

\newcommand{\FF}{{\mathcal{F}}}
\newcommand{\B}{{\mathcal{B}}}

\newcommand{\E}{{\mathcal{E}}}

\newcommand{\M}{{\mathcal{M}}}
\newcommand{\NN}{{\mathcal{N}}}
\newcommand{\D}{{\mathcal{D}}}
\newcommand{\I}{{\mathcal{I}}}
\newcommand{\PP}{{\mathcal{P}}}
\newcommand{\QQ}{{\mathcal{Q}}}
\renewcommand{\O}{{\mathcal{O}}}

\newcommand{\V}{\mathcal{V}}
\renewcommand{\S}{\mathcal{S}}
\newcommand{\Y}{\mathcal{Y}}
\newcommand{\ZZ}{\mathcal{Z}}
\newcommand{\X}{\mathfrak{X}}

\newcommand{\U}{\mathfrak{U}}

\newcommand{\A}{\mathbb{A}}
\renewcommand{\P}{\mathbb{P}}
\newcommand{\F}{\mathbb{F}}
\newcommand{\C}{\mathbb{C}}
\newcommand{\DD}{\mathbb{D}}
\renewcommand{\L}{\mathbb{L}}
\newcommand{\R}{\mathbb{R}}
\newcommand{\Q}{\mathbb{Q}}
\newcommand{\Z}{\mathbb{Z}}
\newcommand{\N}{\mathbb{N}}

\newcommand{\hdag}{  \phantom{}{^{\dag} }    }

\begin{document}
\title{Dévissages des $F$-complexes de $\mathcal{D}$-modules arithmétiques en $F$-isocristaux surconvergents}
\titlerunning{$F$-complexes de $\mathcal{D}$-modules arithmétiques dévissables}

\author{D. Caro 
}                     
%
%
\institute{NWF I - Mathematik
Universität Regensburg
Universitätsstra\ss e 31 \\
93053 Regensburg
DEUTSCHLAND}
\date{Received: / Revised version:}
%
\maketitle
\begin{abstract}
Nous d\'efinissons et \'etudions
les d\'evissages des $F$-complexes de $\mathcal{D}$-modules arithm\'etiques en $F$-isocristaux surconvergents.
Nous prouvons que les $F$-complexes surholonomes sont d\'evissables en $F$-isocristaux surconvergents.
On \'etablit ensuite une formule cohomologique,
\'etendant celle d'\'Etesse et Le Stum,
des fonctions $L$ associ\'ees aux duaux des $F$-complexes de $\mathcal{D}$-modules arithm\'etiques
d\'evissables en $F$-isocristaux surconvergents.
Puis, nous obtenons un analogue $p$-adique de Weil II g\'en\'eralisant celui de Kedlaya.
\end{abstract}
\keywords{Arithmetical $\D$-modules, overconvergent isocrystal, Frobenius, holonomicity, L-function, weight}


\tableofcontents


\section*{Introduction}

Soient $X$ un schéma lisse sur le corps des complexes et $\D _X$ l'anneau des opérateurs différentiels sur $X$.
Pour tout $\D _X$-module (toujours à gauche par défaut) holonome $\E$,
il existe un ouvert dense sur lequel $\E$ devient $\O _X$-cohérent.
Plus généralement, si $\E$ est un complexe de $\D _X$-modules à cohomologie bornée et holonome,
il existe une partition de $X = \cup _r X _r$ en sous-schémas lisses
telle que la restriction de $\E$ sur chaque strate $X _r$ soit à cohomologie $\O _{X _r}$-cohérente.
On remarque qu'un analogue $l$-adique est le dévissage
au dessus d'une stratification d'un faisceau constructible en faisceaux lisses.

Le but de ce travail est de fournir un analogue $p$-adique de ce {dévissage}.
Précisons d'abord celui-ci.

Soient $\V$ un anneau de valuation discrète complet d'inégales caractéristiques
($0$, $p$), de corps résiduels $k$, de corps des fractions $K$,
$\PP$ un $\V$-schéma formel lisse, $P$ sa fibre spéciale, $T$ un diviseur de $P$
et $U$ l'ouvert de $P$ complémentaire de $T$.
Berthelot a construit le faisceau des opérateurs différentiels d'ordre infini et de niveau fini sur $\PP$
à singularités surconvergentes le long de $T$,
$\D ^\dag _{\PP} (\hdag T) _{\Q}$ (voir \cite{Be1}).
Ce faisceau est l'analogue $p$-adique de $\D _U$.
Pour illustrer ce fait, rappelons qu'il existe un foncteur canonique pleinement fidèle de
la catégorie des ($F$-)isocristaux sur $U$ surconvergent le long de $T$ dans celle
des ($F$-)$\D ^\dag _{\PP} (\hdag T) _{\Q}$-modules cohérents
(\cite[4.4.5]{Be1} et \cite{Be2}).

Berthelot conjecture que l'image essentielle des $F$-isocristaux
sur $U$ surconvergent le long de $T$
est incluse dans celle des $F\text{-}\D ^\dag _{\PP, \Q}$-modules holonomes
(pour la structure induite
par l'extension $\D ^\dag _{\PP, \Q}\rightarrow \D ^\dag _{\PP} (\hdag T) _{\Q}$).
Réciproquement, soit $\E$ un $F\text{-}\D ^\dag _{\PP,\Q}$-module holonome.
Il existe alors un diviseur $T _{\E}$ de $P$ tel que la {\og restriction \fg}
(avec le sens de \cite[2.2.6]{caro_surcoherent})
de $\E$ en dehors de $T _{\E}$
soit un $F \text{-}$isocristal sur $P \setminus T _{\E}$ surconvergent le long
de $T _{\E}$.

Pour donner un sens aux dévissages de $F$-complexes de $\D$-modules arithmétiques,
il s'agit dans un premier temps d'exprimer, avec le langage des $\D$-modules arithmétiques,
un analogue $p$-adique des {\og $\D _Y$-modules $\O _Y$-cohérents \fg},
où $Y$ est un schéma lisse sur le corps des complexes.
Or, les {\og $F$-isocristaux surconvergents sur $Y$ \fg},
avec $Y$ un $k$-schéma séparé et lisse,
forment de manière conjecturale cet analogue $p$-adique.

Pour valider cette analogie, il suffirait de construire un foncteur pleinement fidèle
de la catégorie des $F$-isocristaux surconvergents sur $Y$ dans celle des
$F\text{-}\D$-modules arithmétiques {\og sur $Y$ \fg}. Comme on vient de le voir,
un tel foncteur existe déjà lorsque $Y$ est la fibre spéciale de l'ouvert (complémentaire
d'un diviseur) d'un $\V$-schéma formel propre et lisse.
L'étude du cas général fait l'objet de ce travail.
Voici en gros le résultat obtenu.

Nous prouvons que, pour tout $k$-schéma lisse $Y$,
il existe un ouvert dense $\widetilde{Y}$ de $Y$ et
un foncteur canonique (explicite) pleinement fidèle
de la catégorie des $F$-isocristaux surconvergents sur $\widetilde{Y}$
dans celle des $F\text{-}\D$-modules arithmétiques sur {$\widetilde{Y}$}.
Ainsi, cette construction est {\og génériquement \fg} résolue.

Cela suffit ensuite à définir convenablement les
{\og $F$-complexes de $\D$-modules arithmétiques dévissables en $F$-isocristaux surconvergents\fg}.

L'intérêt fondamental de ces dévissages est de ramener
l'étude des $F$-complexes de $\D$-modules arithmétiques dévissables à celle des $F$-isocristaux surconvergents.
Voici un exemple d'une telle application.

Deligne a défini la notion de {\og faisceaux constructibles
(pour la topologie étale $l$-adique)
de poids mixte \fg} (voir \cite{deligne-weil-II}) et a prouvé que le $i $-ième
espace de cohomologie de l'image directe extraordinaire d'un faisceau mixte de poids $\leq n$ est un faisceau mixte
de poids $\leq n+i$. Un tel résultat est nommé {\og Weil II \fg} en raison des conjectures de Weil qui en découlent.
Or, Kedlaya a prouvé un analogue $p$-adique de Weil II pour les $F$-isocristaux surconvergents (\cite{kedlaya-weilII}).
Nous en déduisons alors par dévissage un analogue $p$-adique de Weil II
pour les $F$-complexes de $\mathcal{D}$-modules arithmétiques dévissables.

Décrivons à présent le contenu des différentes parties de ce travail.
Dans une première partie, nous étudions les $\D$-modules arithmétiques sur les $\V$-schémas formels faibles lisses.
Nous comparons alors l'image directe, l'image inverse extraordinaire et le foncteur cohomologique local
avec ceux définis sur les $\V$-schémas formels déduits par complétion.

  Soient $P ^\dag$ un $\V$-schéma formel faible lisse, $P $ sa fibre spéciale,
  $\PP$ son complété $p$-adique,
  $T$ un diviseur de $P $, $U ^\dag$ l'ouvert de
$P ^\dag $ complémentaire de $T $, $Y ^\dag \hookrightarrow U ^\dag$ une immersion fermée
de $\V$-schémas formels faibles et $Y$ la fibre spéciale de $Y ^\dag$.
On suppose en outre $Y ^\dag$ affine et lisse.
Nous terminons ce chapitre par la construction et l'étude
d'un foncteur canonique de la catégorie des isocristaux surconvergents sur $Y$
dans celle des $\D ^\dag _{\PP} (\hdag T) _{\Q}$-modules cohérents à support dans l'adhérence de $Y$ dans $P$.
Avant de passer à la deuxième partie qui s'intéresse à une construction inverse,
signalons que ce premier travail s'inspire de la théorie des $\mathcal{D}$-modules arithmétiques
élaborée par Mebkhout et Narv\'aez-Macarro (\cite{Mebkhout-Narvaez-Macarro_D})
ainsi que du théorème de comparaison de Noot-Huyghe (\cite{huyghe-comparaison}).

Soient $X$ un sous-schéma fermé d'un $\V$-schéma formel propre et lisse $\PP$, $T$ un diviseur de $P$,
tel que $Y:= X \setminus T $ soit lisse.
Dans une deuxième partie, nous construisons, une sous-catégorie pleine
de celle des $F\text{-}\D ^\dag _{\PP} (\hdag T) _{\Q}$-modules cohérents à support dans $X$.
Nous nommerons (voir \ref{isosurcohpf}) ses objets {\og $F$-isocristaux surcohérents sur $Y$ \fg}
(cette catégorie ne dépend que de $Y$).
On construit ensuite un foncteur pleinement fidèle de la catégorie des
$F$-isocristaux surcohérents sur $Y$ dans celle des $F$-isocristaux surconvergents sur $Y$.
Lorsque $Y$ satisfait certaines propriétés géométriques, celles-ci étant
au moins validées sur un ouvert dense de $Y$, on obtient en fait,
via les deux foncteurs précédemment construits, des équivalences de catégorie quasi-inverses
entre les  $F$-isocristaux surcohérents sur $Y$ et les $F$-isocristaux surconvergents sur $Y$.

Dans une dernière partie, nous définissons
les {\og $F$-complexes de $\D$-modules arithmétiques dévissables en $F$-isocristaux surconvergents \fg}.
Puis, nous décrivons une condition suffisante pour obtenir un tel dévissage.
On en tire en particulier que
les $F$-complexes surholonomes (voir \cite{caro_surholonome}) se dévissent en $F$-isocristaux surconvergents.

Il faut souligner qu'une des conjectures de Berthelot sur la stabilité de l'holonomie
(\cite[5.3.6.D]{Beintro2})
implique l'équivalence entre {\og surholonomie \fg} et {\og holonomie \fg}, et entraîne également
celle entre {\og $F$-complexes dévissables \fg} et {\og $F$-complexes holonomes \fg}.

On termine ce travail en étendant quelques propriétés, notamment les formules cohomologiques des fonctions $L$
et {\og Weil II \fg}, des $F$-isocristaux surconvergents aux $F$-complexes de $\D$-modules arithmétiques dévissables.

\bigskip

\noindent \textbf{Notations}
La lettre $\V$ désignera un anneau de valuation discrète complet,
de corps résiduel parfait $k$ de caractéristique $p>0$, de corps de
fractions $K$ de caractéristique $0$, d'idéal maximal $\mathfrak{m}$ et $\pi$ une uniformisante.
Les $\V$-schémas formels seront notés par des lettres calligraphiques ou
gothiques et leur fibre spéciale par les lettres romanes
correspondantes.
De plus, $s\geq 1$ sera un entier naturel et $F$ la puissance
$s$-ième de l'endomorphisme de Frobenius. Les modules sont par défaut des modules à gauche.
Si $\E$ est un faisceau abélien, $\E _\Q$ désignera $\E \otimes _\Z \Q$.

Si $f$ : $\X ' \rightarrow \X$ est un morphisme de $\V$-schémas formels lisses,
$f _0$ (ou $f$) : $X' \rightarrow X$ sera le morphisme induit.
Lorsque $T$ et $T'$ sont respectivement des diviseurs de $X$ et $X'$ tels que
$f ( X '\setminus T') \subset X \setminus T$, nous désignerons par
$f ^! _{T',T}$ et $f _{T,T',+}$ les foncteurs
image inverse extraordinaire, image directe
par $f$ (voir \cite[3.4, 3.5, 4.3]{Beintro2} et \cite[1.1.5]{caro_courbe})
à singularités surconvergentes le long de $T$ et $T'$.
Pour ne pas les confondre avec d'autres opérations cohomologiques analogues, il nous
arrivera d'ajouter le symbole $\dag$ à ceux-ci.
Si $T' = f ^{-1} (T)$, on les
notera $f ^! _T$ et $f _{T,+}$, ou simplement $f ^!$ et $f _{+}$.
Les produits tensoriels externes et internes seront notés respectivement
$\overset{\L ^\dag}{\boxtimes}$ et
$\overset{\L ^\dag}{\otimes}  $ (\cite[4.3]{Beintro2}).
En outre, si $Z$ est un sous-schéma fermé de $X$,
$\R \underline{\Gamma} ^\dag _Z $ désignera le foncteur cohomologique local
à support strict dans $Z$ (au sens de \cite[2.2.6]{caro_surcoherent})
et $(\hdag Z)$ le foncteur restriction (\cite[2.2.6]{caro_surcoherent}).
Si $T' \subset T$ sont deux diviseurs de $X$, on notera abusivement
$(\hdag T)$ à la place de $(\hdag T ,T')$.

Pour tout diviseur $T$ de $X$, nous désignerons par $\DD _{\X,T}$ ou $\DD _T$,
le foncteur dual $\D ^\dag _{\X} (\hdag T) _{\Q}$-linéaire (voir \cite[I.3.2]{virrion} pour la définition
des foncteurs duaux).
Dans la catégorie des $F$-complexes de $\D ^\dag _{\PP} (\hdag T) _{\Q}$-modules à cohomologie cohérente
  et bornée,
  on note $\DD ^* _T$ le foncteur
  $\DD _T (-)\otimes _{\O _{\X } ( \hdag T ) _{\Q}} \DD _T ( \O _{\X } ( \hdag T ) _{\Q} ) ^{\vee}$,
  où $\vee $ désigne le dual $\O _{\X } ( \hdag T ) _{\Q}$-linéaire.
On dispose d'un isomorphisme $\DD ^* _T \riso \DD _T $ dans la catégorie des
complexes de $\D ^\dag _{\PP} (\hdag T) _{\Q}$-modules à cohomologie cohérente
  et bornée, ce dernier n'étant pas à priori compatible à Frobenius (voir \cite{caro_comparaison}).

Enfin, si aucune confusion n'est à craindre, il nous arriva d'omettre le foncteur
canonique
$\underset{\longrightarrow}{\lim}$ : $\smash{\underset{^{\longrightarrow}}{LD}} ^{\mathrm{b}} _{\Q
,\mathrm{qc}} ( \smash{\widehat{\D}} _{\X} ^{(\pmb{\cdot})}(T))
\rightarrow D (\D ^\dag _{\X,\Q} (\hdag T))$ (voir \cite[4.2.2]{Beintro2} lorsque
le diviseur est vide, mais la construction est analogue).

Si $T$ est l'ensemble vide, nous omettrons de l'indiquer dans toutes ses opérations.
Sauf mention contraire, on supposera (sans nuire à la généralité) les $k$-schémas réduits.

\section{Isocristaux surconvergents sur les schémas affines et lisses}
\subsection{$\D$-modules sur les $\V$-schémas formels faibles lisses}
On appellera {\it $\V$-algèbre f.c.t.f}, une $\V$-algèbre commutative faiblement complète de type fini, i.e.,
un quotient de $\V [ t _1, \dots, t _d] ^\dag$.
Une $\V$-algèbre f.c.t.f. est munie de la topologie $p$-adique.
Un morphisme d'algèbres f.c.t.f. est une application entre deux $\V$-algèbres f.c.t.f. qui est un morphisme d'anneaux
(unitaires). Un morphisme de $\V$-algèbres f.c.t.f. est donc toujours continu.
Si $A \rightarrow B$ et $A \rightarrow C$ sont deux morphismes de $\V$-algèbres f.c.t.f.,
le complété faible de $B \otimes _A C$ sera noté $B \otimes ^\dag _A C$.
Comme le morphisme canonique $B \otimes _\V A \rightarrow B \otimes _A C$ est surjectif,
on vérifie que la $\V$-algèbre $B \otimes ^\dag _A C$ est f.c.t.f.

Soit $A$ une $\V$-algèbre f.c.t.f. On note $A _0 = A /\pi A$ et, pour tout $f \in A$,
$\overline{f}$ désigne l'image de $f$ dans $A _0$ et $A [f]$, la complétion faible de $A _f$.
Suivant Meredith le {\it schéma formel faible affine} associé à $A$,
noté $\mathrm{Spff}  (A)$, est l'espace annelé dont l'espace topologique sous-jacent correspond à
$X _0 =\Spec A _0$ (cet ensemble est aussi égal à celui des idéaux premiers ouverts de $A$) et
dont les sections du faisceau structural $\O _X$ sur un ouvert principal $X _{\overline{f}}$, avec $f \in A$, coïncident
avec $A [f] $ (voir \cite{meredith-weakformalschemes}). On notera $D (f)$ l'ouvert de $\mathrm{Spff}  (A)$
dont l'espace topologique est $X _{0 \overline{f}}$. Les ouverts de la forme $D(f)$ seront dits {\it principaux}.
Nous disposons du théorème de type $A$ (\cite[3.3]{meredith-weakformalschemes}) sur $X=\Spff (A)$ :
les foncteurs $M \mapsto \widetilde{M}:= \O _X \otimes _A M$ et $\M \mapsto \Gamma (X,  \M)$ sont des
équivalences quasi-inverses entre la catégorie des $\O _X$-modules cohérents et celle des
$A$-modules de type fini. En outre, on bénéficie du théorème de type $B$ (\cite[2.14]{meredith-weakformalschemes}) :
pour tout $\O _X$-module cohérent $\M$, pour tout entier $i>0$, $H ^i(X, \M)=0$.

Soit $\mathfrak{P}$ un idéal ouvert de $A$. On pose $A _{[\mathfrak{P}]}=\underset{\longrightarrow}{\lim}
_{f \not \in \mathfrak{P}} A [f]$ (le système inductif est filtrant). Celui-ci est muni de l'idéal
$\mathfrak{P} A _{[\mathfrak{P}]}=\underset{\longrightarrow}{\lim}
_{f \not \in \mathfrak{P}} \mathfrak{P} A [f]$. Pour tout $f \not \in \mathfrak{P}$,
$\mathfrak{P}  A [f] / \pi  A [f] \riso
(\mathfrak{P}  /\pi A ) _{\overline{f}}\neq (A _0 ) _{\overline{f}}$. Cela entraîne que
$\mathfrak{P} A [f] \neq A [f]$ et a fortiori
$\mathfrak{P} A _{[\mathfrak{P}]} \neq A _{[\mathfrak{P}]}$
(on vérifie que $1\not \in \mathfrak{P} A _{[\mathfrak{P}]}$).
La proposition qui suit, analogue à \cite[0.6.7.17]{EGAI},  précise ce dernier fait.

\begin{prop}\label{a[p]local}
  L'anneau $A _{[\mathfrak{P}]}$ est local et d'idéal maximal $\mathfrak{P}A _{[\mathfrak{P}]}$.
  De plus, son corps résiduel est isomorphe au corps des fractions de $A / \mathfrak{P}$.
\end{prop}
\begin{proof}
Soit $x$ un élément de $A _{[\mathfrak{P}]}$ n'appartenant pas à $\mathfrak{P}A _{[\mathfrak{P}]}$.
Pour établir la première assertion, il s'agit de vérifier que $x $ est inversible.
Il existe $f  \not \in \mathfrak{P}$ tel que
$x$ provienne d'un élément $y$ de $A [f ]$. Par hypothèse, l'image canonique $\overline{y}$
de $y$ sur $A [f]/\pi A [f] \riso (A _0 ) _{\overline{f}}$ n'appartient pas à
$(\mathfrak{P}  /\pi A ) _{\overline{f}}$. Cela implique donc que
$\overline{y} = \overline{a} / \overline{f} ^r$, où $a \in A \setminus \mathfrak{P}$ et $r \in \N$
Ainsi, $g := a f \not \in \mathfrak{P}$ et l'image canonique de $\overline{y}$ sur
$(A _0 ) _{\overline{g}}$ est inversible.
    Comme $\pi A[g]$ est inclus dans l'idéal de Jacobson de $A[g]$, il en découle que
    l'image canonique de $y$ sur $A [g]$ est inversible. L'élément $x$ est par conséquent inversible.

    Soit $S = \widehat{A} \setminus \widehat{\mathfrak{P}}$. En vertu de \cite[0.6.7.17]{EGAI},
    $\widehat{A} \{ S \}$ est un anneau local d'idéal maximal $\widehat{\mathfrak{P}} \widehat{A} \{ S \}$.
    L'homomorphisme canonique $A _{[\mathfrak{P}]} \rightarrow\widehat{A} \{ S \}$ est local.
    Il induit un morphisme
    $A _{[\mathfrak{P}]} \rightarrow\widehat{A} \{ S \}/\widehat{\mathfrak{P}} \widehat{A} \{ S \}
    \tilde{\leftarrow} \mathrm{Frac} (A / \mathfrak{P})$ (l'isomorphisme se prouve grâce
    à \cite[0.6.7.17]{EGAI} et via $A / \mathfrak{P} \riso \widehat{A} /\widehat{\mathfrak{P}} \widehat{A}$).
    Un élément de $\mathrm{Frac} (A / \mathfrak{P})$ provient aisément d'un élément de
    $A _{[\mathfrak{P}]}$, i.e., ce dernier homomorphisme est surjectif. D'où l'isomorphisme
    de corps résiduels :
    $A _{[\mathfrak{P}]} / \mathfrak{P}A _{[\mathfrak{P}]}
    \riso
    \widehat{A} \{ S \}/\widehat{\mathfrak{P}} \widehat{A} \{ S \}$.
\hfill \hfill \qed \end{proof}

\begin{vide}
  \label{defaphiphitilde}
  Soient $\phi $ : $A \rightarrow B$ un morphisme de $\V$-algèbres f.c.t.f.,
  $(X, \O _X) := \Spff A$ et $(Y,\O _Y) := \Spff B$.
  L'image inverse par $\phi$ d'un idéal premier ouvert est un idéal premier ouvert.
  On obtient un morphisme d'espaces topologiques $\overset{a}{} \phi$ : $Y \rightarrow X$.
  Pour tout élément $f$ de $A$, on dispose d'un  morphisme canonique
  $A [f]\rightarrow B [\phi (f)]$ et,
  pour tout multiple $f'$ de $f$, du diagramme commutatif :
  $$\xymatrix @R=0,3cm {
  {A [f]} \ar[r] \ar[d]
  &
  {B [\phi (f)]}
  \ar[d]
  \\
  {A [f']} \ar[r]
  &
  {B [\phi (f')].}
  }$$
  Comme $\overset{a}{} \phi ^{-1} D(\phi (f))= D (f)$,
  ces homomorphismes définissent donc un homomorphisme de faisceaux d'anneaux
  $\widetilde{\phi}$ : $ \O _X \rightarrow \overset{a}{} \phi _* \O _Y$.
  On a donc construit un morphisme $(\overset{a}{} \phi, \widetilde{\phi})$ : $\Spff B \rightarrow \Spff A$
  d'espaces annelés. De plus, pour tout idéal premier $\mathfrak{Q}$ de $B$,
  on dispose d'un homomorphisme d'anneaux locaux
  $A _{[\phi ^{-1}(\mathfrak{Q})]} \rightarrow B _{[\mathfrak{Q}]}$.
  L'homomorphisme $(\overset{a}{} \phi, \widetilde{\phi})$ est ainsi un homomorphisme d'espaces localement annelés.
  Via la proposition ci-après, ils sont tous de cette forme.
\end{vide}

\begin{prop}
\label{egaI.10.2.1}
  Soient $A$ et $B$ deux $\V$-algèbres f.c.t.f. et soient $X = \Spff A$, $Y= \Spff B$.
  Pour qu'un morphisme $u=(\psi,\theta)$ : $Y \rightarrow X$ d'espaces annelés soit
  de la forme $ (\overset{a}{} \phi , \widetilde{\phi})$, où $\phi $ est un homomorphisme
  d'anneaux $A \rightarrow B$, il faut et il suffit que $u$ soit un morphisme d'espaces localement annelés.
\end{prop}
\begin{proof}
Par le truchement de \ref{a[p]local}, on reprend les arguments de \cite[10.2.2]{EGAI}.
\hfill \hfill \qed \end{proof}

\begin{defi}\label{defschff}
Un {\it $\V$-schéma formel faible} est un espace localement annelé en $\V$-algèbres $(X,\O _X)$ localement
isomorphe à un $\V$-schéma formel faible affine (voir \cite{meredith-weakformalschemes}).
On pourra le noter $X$. Pour tout entier $i$,
le $\V / \pi ^{i+1} \V$-schéma induit par réduction modulo $\pi ^{i+1}$ sera noté
$X _i$. De plus, $\widehat{X}$ ou $\X$ désignera le $\V$-schéma formel déduit par complétion $p$-adique de $X$.

  Un {\it morphisme} $f$ : $Y \rightarrow X$ de $\V$-schémas formels faibles est
un morphisme d'espaces localement annelés.
On pose $f _i$ : $Y _i \rightarrow X _i$ et $\hat{f}$ : $\Y \rightarrow \X$ les morphismes induits.
Par abus de notations, on pourra parfois les noter $f$.
Un tel morphisme est dit {\it lisse} (resp. {\it étale}) si pour tout
$i$, les $f _i$ sont des morphismes lisses (resp. étales).
Enfin, $f$ est {\it séparé} lorsque $f _0$ est séparé. Un $\V$-schéma formel faible $X$
est séparé si son morphisme structural $X \rightarrow \Spff \V$ l'est.
Les propriétés usuelles des morphismes séparés
restent valables.

Soit $Y$ un ouvert de l'espace topologique sous-jacent à
un $\V$-schéma formel faible $X$. On munit $Y$ de la structure canonique
de $\V$-schéma formel faible en posant par abus de notations $Y := (Y, \O _X |_Y)$.
De plus, on dispose d'un morphisme canonique $Y \rightarrow X$. Le composé d'un tel morphisme
avec un isomorphisme sera nommé {\it immersion ouverte}.
On appelle {\it ouvert affine} de $X$, un ouvert de $Y$ induisant un $\V$-schéma formel faible affine.
Nous nous écartons de la terminologie de Meredith qui les nomme {\og affine wf open \fg} afin de les
démarquer des ouverts affines de $X _0$ (\cite[4]{meredith-weakformalschemes}).

Un morphisme $f$ : $Y \rightarrow X$ de $\V$-schémas formels faibles est dit
{\it affine} si, pour tout ouvert affine $X'$ de $X$, $f ^{-1} (X')$ est un ouvert affine de $Y$.
\end{defi}

\begin{prop}\label{egaI.10.4.6}
  Soient $Y$ un $\V$-schéma formel faible, $X=\Spff A$ un $\V$-schéma formel faible affine
  d'anneau $A$. Il existe une correspondance biunivoque canonique entre les morphismes de
  $\V$-schémas formels faibles de la forme $Y \rightarrow X$
  et les homomorphismes d'anneaux de la forme de $A\rightarrow \Gamma (Y,\O _Y)$.
\end{prop}
\begin{proof}
  De manière analogue à \cite[2.2.4 ou 10.4.6]{EGAI}, cela résulte de \ref{egaI.10.2.1}.
\hfill \hfill \qed \end{proof}

\begin{prop}
  La catégorie des $\V$-schémas formels faibles possède des produits fibrés.
\end{prop}
\begin{proof}
  Le cas affine se déduit de \cite[1.5]{MonWas68} : si $A \rightarrow B$ et $A\rightarrow C$ sont
  deux morphismes de $\V$-algèbres f.c.t.f.,
  le schéma faiblement formel
  $\Spff ( B) \times _{\Spff (A)} \Spff (C): = \Spff ( B \otimes ^\dag _A C)$
  vérifie la propriété universelle des produits fibrés.
  Pour le cas général, on procède par recollement (de manière analogue au cas des schémas : voir \cite{EGAI}).
\hfill \hfill \qed \end{proof}

Passons maintenant aux immersions fermées.
Afin d'y voir plus clair, la proposition suivante est utile.
\begin{prop}\label{sous-sch-id}
  Soient $X$ un $\V$-schéma formel faible et $\I$ un idéal cohérent de $\O _X$.
  Si $Y$ est le support (fermé) de $\O _X/\I$,
  l'espace topologique annelé $(Y, (\O _X / \I) |_Y)$ est un $\V$-schéma formel faible.
\end{prop}
\begin{proof}
Comme $\O _X$ est un anneau cohérent, $\O _X / \I$ est cohérent et donc d'après \cite[0.5.2.2]{EGAI},
le support $Y$ est bien un fermé de $X$.
La proposition est locale. Supposons $X =\Spff A$.
Grâce au théorème de type $A$ pour les $\O _X$-modules cohérents,
en notant $I := \Gamma ( X, \I)$, le morphisme canonique
$I \otimes _A \O _X \rightarrow \I$ est un isomorphisme. On obtient
$A /I \otimes _A \O _X \riso \O _X / \I$. Prouvons maintenant que l'on dispose d'un isomorphisme
canonique $(Y, (A /I \otimes _A \O _X) |_Y) \riso \Spff A /I$.
Il s'agit alors de vérifier que l'on dispose, pour tout $f \in A$ et tout $g\in A$ multiple de $f$, d'un isomorphisme
canonique $A / I \otimes _A A [f] \riso A /I [\overline{f}]$ induisant le diagramme commutatif
$$\xymatrix @R=0,3cm {
{A / I \otimes _A A [f]}
\ar[r] \ar[d] ^-\sim
&
{A / I \otimes _A A [g]}
\ar[d] ^-\sim
\\
{A /I [\overline{f}]}
\ar[r]
&
{A /I [\overline{g}],}
}$$
où $\overline{f}, \overline{g} \in A /I$
sont les images canoniques respectives de $f$ et $g$.

L'épimorphisme $A _f \twoheadrightarrow (A /I) _{\overline{f}}$ induit le suivant
$(A _f ) ^\dag \twoheadrightarrow ((A /I) _{\overline{f}}) ^\dag $
puis $A/I \otimes _A (A _f ) ^\dag \twoheadrightarrow ((A /I) _{\overline{f}}) ^\dag $.
La fonctorialité en $f$ de ce dernier est immédiate.
Il reste à s'assurer de son injectivité ce qui résulte
du diagramme commutatif
$$\xymatrix @R=0,3cm {
{A/I \otimes _A (A _f ) ^\dag }
\ar[r] \ar @{^{(}->} [d]
&
{((A /I) _{\overline{f}}) ^\dag }
\ar @{^{(}->} [d]
\\
{A/I \widehat{\otimes} _A (A _f ) ^\dag }
\ar[r] ^-\sim
&
{((A /I) _{\overline{f}}) ^\wedge,}
}$$
dont les flèches verticales sont injectives et celle du bas est un isomorphisme.
On a ainsi prouvé
$\Spff A / I \riso (Y, (\O _X /\I) |_Y)$.
\hfill \hfill \qed \end{proof}

Soit $\I$ un idéal cohérent de $\O _X$. On dispose d'un morphisme canonique
$(Y, (\O _X / \I) |_Y) \rightarrow (X, \O _X)$. De plus,
avec la preuve de \ref{sous-sch-id}, on vérifie que lorsque $X=\Spff A$ et $I := \Gamma (X, \I)$,
ce morphisme correspond via \ref{egaI.10.2.1} au morphisme canonique $A \rightarrow A /I$.
En particulier, $(Y, (\O _X / \I) |_Y) \rightarrow (X, \O _X)$ est affine.
\begin{defi}
  On appelle sous-$\V$-schéma formel faible fermé d'un $\V$-schéma formel faible $X$ tout
  $\V$-schéma formel faible $(Y, (\O _X / \I) |_Y)$, où $\I$ est un idéal cohérent de $\O _X$ ; on dit
  que celui-ci est le sous-$\V$-schéma formel faible fermé défini par $\I$.

Un morphisme $f$ : $Y \rightarrow X$ de $\V$-schémas formels faibles
est une {\og immersion fermée \fg}
s'il se factorise en  $Y \overset{g}{\rightarrow} X ' \overset{u}{\rightarrow} X$,
où $X'$ est un sous-$\V$-schéma formel fermé de $X$ et $g$ est un isomorphisme.
Une immersion fermée est un morphisme affine.

\end{defi}
\begin{vide}\label{caraimmfer}
  De manière analogue à \cite[10.14.4]{EGAI}, on a la caractérisation suivante
  d'une immersion fermée. Soit $f$ : $Y\rightarrow X$ un morphisme de $\V$-schémas formels faibles,
  et soit $(X _\alpha)$ un recouvrement de $f (Y)$ par des ouverts affines de $X$, tels que les
  $f^{-1} (X _\alpha)$ soient des ouverts affines de $Y$. Pour que $f$ soit une immersion fermée,
  il faut et il suffit que
  $f (Y)$ soit une partie fermée de $X$ et que, pour tout $\alpha$,
  l'homomorphisme
  $\Gamma (X _\alpha , \O _X) \rightarrow \Gamma ( f ^{-1} (X _\alpha ), \O _Y)$,
  induit via \ref{egaI.10.2.1} par la restriction de $f$ à $f ^{-1}(X _\alpha )$,
  soit surjective.
\end{vide}

\begin{rema}\label{remaaffinrel}
  Soient $X$ et $Y$ deux schémas formels faibles affines. Si $X$ est lisse et si $f _0$ : $ Y _0 \rightarrow X _0$
  est un morphisme de $k$ schémas, alors
  il existe un morphisme $Y \rightarrow X$ de $\V$-schémas formels faibles relevant $f _0$.
  En effet, puisque $X$ est lisse, il existe un relèvement $\widehat{Y} \rightarrow \widehat{X}$ de $f _0$.
  On conclut ensuite en vertu de \cite[2.4.3]{vanderPutMonsky-Washnitzer}.
\end{rema}

\begin{prop}\label{stabimmfer}
  (i) Si $f$ : $Z \rightarrow Y$ et $g$ : $Y \rightarrow X$ sont des immersions fermées de $\V$-schémas formels
  faibles, $g \circ f$ est une immersion fermée.

  (ii) Soient $X$, $Y$, $Z$ trois $\V$-schémas formels faibles, $f$ : $Y \rightarrow X$ une immersion fermée
  et $Z \rightarrow X$ un morphisme. Le morphisme $ Y \times _X Z \rightarrow Z$ est une immersion fermée.

  (iii) Soient $X$ un $\V$-schéma formel faible, $f$ : $Y \rightarrow Y'$ et
  $g$ : $Z \rightarrow Z'$ des $X$-morphismes qui soient des immersions fermées, alors
  $f \times _X g$ est une immersion fermée.
\end{prop}
\begin{proof}
  Soient $A \rightarrow C$ un morphisme de $\V$-algèbres f.c.t.f., $I$ un idéal de $A$
  et $B := A /I$. L'algèbre $B \otimes _A C$ est une $\V$-algèbre f.c.t.f. Ainsi,
$C /I C \riso B \otimes _A C \riso B \otimes ^\dag _A C$.
  Avec cette remarque, la preuve est analogue à \cite[10.14.5]{EGAI}
\hfill \hfill \qed \end{proof}

La proposition qui suit donne des exemples d'immersions fermées.
\begin{prop}
  Soient $f$ : $Z\rightarrow Y$ et $g$ : $Y \rightarrow X$ des morphismes de $\V$-schémas formels faibles.
  Si $g$ est séparé alors le graphe de $f$, $\Gamma _f = (1, f) _X$ : $Z \rightarrow Z \times _X Y$,
  est une immersion fermée.
\end{prop}
\begin{proof}
  Analogue à \cite[10.15.4]{EGAI}.
\hfill \hfill \qed \end{proof}

\begin{defi}
  Une {\og immersion \fg} est le composé d'une immersion fermée suivi d'une immersion ouverte.

\end{defi}

\begin{prop}\label{stabimm}
  (i) Si $f$ : $Z \rightarrow Y$ et $g$ : $Y \rightarrow X$ sont des immersions de $\V$-schémas formels
  faibles, $g \circ f$ est une immersion.

  (ii) Soient $X$, $Y$, $Z$ trois $\V$-schémas formels faibles, $f$ : $Y \rightarrow X$ une immersion
  et $Z \rightarrow X$ un morphisme. Le morphisme $ Y \times _X Z \rightarrow Z$ est une immersion.

  (iii) Soient $X$ un $\V$-schéma formel faible, $f$ : $Y \rightarrow Y'$ et
  $g$ : $Z \rightarrow Z'$ des $X$-morphismes qui soient des immersions, alors
  $f \times _X g$ est une immersion.
\end{prop}
\begin{proof}
  Traitons d'abord $(i)$.
Soient $j$ : $Z \rightarrow Y$ une immersion ouverte et $u$ : $Y \rightarrow X$ une immersion fermée.
L'image de $Z$ par $j$ est un ouvert de $Y$. Il existe donc un ouvert $X'$ de $X$ tel que
$u ^{-1} (X') = j (Z)$. En notant $u'$ : $j(Z) \rightarrow X'$ l'immersion fermée induite par $u$,
$ u \circ j$ se décompose en $Z \overset{j}{\rightarrow} j (Z) \overset{u'}{\rightarrow} X' \subset X$.
Ainsi, $u \circ j$ est le composé d'une immersion fermée suivi d'une immersion ouverte.
Les assertions (ii) et (iii) découlent de \ref{stabimmfer}.
\hfill \hfill \qed \end{proof}

\begin{prop}\label{diagimm}
  Soit $f$ : $X \rightarrow Y$ un morphisme de $\V$-schémas formels faibles.
Le morphisme canonique $\delta=(1, 1)_Y$ : $X \rightarrow X \times _Y X$ est une immersion.
On l'appellera {\og l'immersion diagonale \fg}.
\end{prop}
\begin{proof}
  Soient $(X _\alpha)$ et $(Y _\alpha )$ des recouvrements respectifs de $X$ et $Y$ par des ouverts affines
tels que $f$ se factorise par $X _\alpha \rightarrow Y _\alpha$.
 Par construction du produit fibré $X \times _Y X$, $X _\alpha\times _{Y _\alpha } X _\alpha$
est un ouvert de $X \times _Y X$ et $\delta ^{-1} (X _\alpha\times _{Y _\alpha } X _\alpha) = X _\alpha$.
Notons alors $Y '$ l'ouvert de $Y$ réunion des $X _\alpha\times _{Y _\alpha } X _\alpha$.
Le morphisme $\delta $ se factorise alors en un morphisme
$\delta '$ : $X \rightarrow Y'$.
D'après la caractérisation \ref{caraimmfer} des immersions fermées, pour prouver que
$\delta '$ est une immersion fermée, il suffit de vérifier que
les homomorphismes canoniques
$\Gamma ( X _\alpha , \O _{X _\alpha})
\rightarrow
\Gamma ( X _\alpha \times _{Y _\alpha} X_\alpha , \O _{X _\alpha \times _{Y _\alpha} X_\alpha})$
sont surjectifs, ce qui est immédiat.

\hfill \hfill \qed \end{proof}

\begin{coro}
  Soient $f$ : $X \rightarrow S$, $g$ : $Y \rightarrow S$ et $\phi$ : $S \rightarrow T$ des morphismes
de $\V$-schémas formels faibles. Le morphisme canonique
$X \times _S Y \rightarrow X \times _T Y$ induit par $\phi$ est une immersion.
En particulier, lorsque $S =Y$, le graphe de $f$, $X \rightarrow X \times _T Y$, est une immersion.
\end{coro}
\begin{proof}
  Cela résulte de
$X \times _S Y \riso S \times _{(S \times _T S)} (X \times _T Y)$, de \ref{diagimm} et
du fait que les immersions fermées sont stables par changement de base (\ref{stabimm}).
\hfill \hfill \qed \end{proof}

\begin{vide}

Soit $X$ un $\V$-schéma de type fini. Meredith (\cite[4]{meredith-weakformalschemes})
construit le faisceau $\O ^\dag _X$ de la façon suivante :
si $U \subset V \subset X$ sont des ouverts affines de $X$, alors
$\Gamma (U, \O ^\dag _X) := \Gamma (U, \O _X) ^\dag$ et
$\Gamma (V, \O ^\dag _X) \rightarrow \Gamma (U, \O ^\dag _X)$ est le morphisme canonique
induit par $\Gamma (V, \O  _X) \rightarrow \Gamma (U, \O _X)$ via \cite[1.5]{MonWas68}.
En notant $X ^\dag$ l'espace topologique de $X _0 = X \times _{\Spec \V} \Spec (k)$,
on vérifie que l'espace annelé $(X ^\dag, \O ^\dag _X)$ est un $\V$-schéma formel faible
et que l'on dispose d'un morphisme canonique $(X ^\dag , \O ^\dag _{X}) \rightarrow (X, O _X)$,
qu'il appelle complétion faible de $(X, \O_X)$. Par abus de notation, on écrira
$X ^\dag $ à la place de $(X ^\dag , \O ^\dag _X)$. On vérifie de plus
que l'application $X \mapsto X ^\dag$ induit canoniquement (via \cite[1.5]{MonWas68})
un foncteur de la catégorie des $\V$-schémas de type fini dans celle des
$\V$-schémas formels faibles.

\end{vide}

\begin{defi}
  Soit $U$ un $\V$-schéma formel faible.
  On dit que $U $ a des coordonnées locales s'il existe un morphisme étale
  $U \rightarrow \A ^{d \dag} _\V $.
  Les sections $t _1,\dots, t_d$ induites via ce morphisme par $X _1, \dots,X_d$
  seront appelés {\it coordonnées locales}.
\end{defi}

\begin{rema}
  Soit $U$ un $\V$-schéma formel faible affine et sans torsion.
  Alors, $U _0$ a des coordonnées locales si et seulement si $\widehat{U}$ en a
  si et seulement si $U$ en a.
\end{rema}

\begin{lemm}\label{suireg}
  Soient $U$ un schéma formel faible affine muni de coordonnées locales $t _1,\dots,t_d$ et
  $\tau _1 =1 \otimes t_1 -t_1 \otimes 1, \dots , \tau _d =1 \otimes t_d -t_d \otimes 1$.
  La suite $\tau _1, \dots ,\tau _d$ est une suite
  régulière de générateurs de l'idéal définissant l'immersion fermée $U \hookrightarrow U \times _S U$.
\end{lemm}
\begin{proof}
Cela résulte du {\og cas formel \fg}.
  En effet, soient $A $ la $\V$-algèbre f.c.t.f. de $U$ et $I $ l'idéal de $A\otimes ^\dag _\V A$
  correspondant à l'immersion
fermée $U \hookrightarrow U\times U$. Les images canoniques de $t _1,\dots ,t _d$ dans $\widehat{A}$ sont
des coordonnées locales de $\widehat{U}$.
De plus, par fidèle platitude de $A\otimes ^\dag _\V A \rightarrow A \widehat{\otimes}  _\V A$,
comme $\widehat{I} \riso I \otimes _{A\otimes ^\dag _\V A} A \widehat{\otimes}  _\V A$
est engendré par les images de $\tau _1,\dots ,\tau _d$,
alors $\tau _1, \dots ,\tau _d$ engendre $I$.
De plus,
$I / ( \tau _1,\dots , \tau _r) \rightarrow \widehat{I} / ( \tau _1,\dots , \tau _r)$ est injectif.
Comme la multiplication par $\tau _{r+1}$ est injective dans
$\widehat{I} / ( \tau _1,\dots , \tau _r)$, elle l'est alors aussi
dans $I / ( \tau _1,\dots , \tau _r)$.
\hfill \hfill \qed \end{proof}

\begin{vide} Dans la suite de cette section, $m$ sera un entier positif fixé et
  $U$ un $\V$-schéma formel faible.
  On note $\mathcal{I}$ l'idéal cohérent de l'immersion diagonale :
  $U \hookrightarrow U \times _S U$. On définit l'{\it algèbre
  des parties principales de niveau $m$ et d'ordre $n$ de $U$} et noté
  $\PP ^n _{U, (m)}$, comme étant l'enveloppe à puissance divisée de niveau $m$ et d'ordre $\leq n$ de $\mathcal{I}$.
  Les deux projections
  canoniques $U \times _S U \rightarrow U$ induisent deux structures de $O _U$-algèbres sur $\PP ^n _{U, (m)}$ : la
  structure à gauche et la structure à droite.
  Il résulte de \cite[1.5.3]{Be1} et de \ref{suireg} que, si $U $ est muni de coordonnées locales $t _1,\dots ,t_d$,
  pour chacune des structures de $\O_U$-algèbre, le faisceau $\PP ^n _{U, (m)}$ est un $\O _U$-module
  libre dont les éléments
  $\underline{\tau} ^{\{\underline{k}\}}= \tau _1 ^{\{k_1\}}\cdots \tau _d ^{\{k_d\}}$,
  pour $\underline{k}\leq n$, forment une base.

Le {\it faisceau des opérateurs différentiels de niveau $m$ et d'ordre $\leq n$ sur $U$}, noté $\D _{U,n}^{(m)}$,
est le dual $\O _U$-linéaire pour la structure gauche de $\O _U$-algèbres de $\PP ^n _{U, (m)}$.
Le {\it faisceau des opérateurs différentiels de niveau $m$ sur $U$} est la réunion :
$\D _{U }^{(m)} := \cup _{n\in \N} \; \D _{U,n}^{(m)}.$ De manière analogue à \cite[2.2.1]{Be1}, on munit
ce faisceau d'une structure d'anneau et de deux structures de $\O _U$-algèbre, l'une à droite et l'autre à gauche.

Si $U$ est muni de coordonnées locales, à la base de $\PP ^n _{U, (m)}$ formée des
$\underline{\tau} ^{\{\underline{k}\}}$ pour $\underline{k}\leq n$, correspond la base duale de
$\D _{U,n}^{(m)}$, dont les éléments seront notés $\underline{\partial} ^{<\underline{k}>}$.
Les $\underline{\partial} ^{<\underline{k}>}$ forment donc une base de $\D _{U }^{(m)} $.

Les propositions \cite[2.2.4 et 2.2.5]{Be1} sont encore valables en remplaçant {\og formel\fg} par {\og faiblement formel\fg}.
De même les définitions et résultats sur les $m$-PD-stratifications (\cite[2.3]{Be1}) sont encore valables
en remplaçant {\og formel\fg} par {\og faiblement formel\fg}. Les modules sont par défaut à gauche. Toutefois,
tous les résultats de cette section restent valables pour les modules à droite.
\end{vide}

\begin{lemm}\label{lemm-wildetildeDcoh}
  L'anneau $\Gamma (U, \D ^{(m)} _U)$
  (resp. $\D ^{(m)} _{U ,x}$ pour tout $x \in U$) est noethérien à droite et à gauche.
\end{lemm}
\begin{proof}
  Analogue à \cite[2.2.5]{Be1}.
\hfill \hfill \qed \end{proof}
\begin{prop}\label{wildetildeDcoh}
  Soit $j$ : $U \hookrightarrow P$ une immersion ouverte de $\V$-schémas formels faibles telle que
  $P _0 \setminus U _0$ soit le support d'un diviseur.
  Les faisceaux d'anneaux
  $\D ^{(m)} _P$ et $j _* \D _U ^{(m)}$ sont cohérents à droite et à gauche.
\end{prop}
\begin{proof}
  La cohérence à droite et à gauche de $\D ^{(m)} _P$ est un cas particulier de celle de $j _* \D _U ^{(m)}$.
  Prouvons donc la cohérence à gauche de $j _* \D _U ^{(m)}$.
  La preuve est similaire à celle de \cite[3.1.2]{Be1} :
  il existe une base d'ouverts $\mathcal{B}$ de $P$ telle que
  pour tout $V \in \mathcal{B}$, $V$ et $V \cap U$ sont affines.
  Par \ref{lemm-wildetildeDcoh}, l'anneau $\Gamma ( V ,j _* \D _U ^{(m)}) = \Gamma(U\cap V, \D _U ^{(m)})$
  est noethérien à droite et à gauche. De plus, pour tous ouverts $V '\subset V$ de $\mathcal{B}$,
  $\Gamma ( U \cap V, \O _U) \rightarrow \Gamma ( U \cap V', \O _U)$ est plat (c'est un homomorphisme de $\V$-algèbres
  f.c.t.f. lisses dont la réduction modulo $\pi$ est plate).
  Comme $\Gamma (U \cap V', \O _U) \otimes _{\Gamma(U \cap V, \O _U)} \Gamma ( V ,j _* \D _U ^{(m)})
  \riso \Gamma ( V ',j _* \D _U ^{(m)})$,
  $\Gamma ( V ,j _* \D _U ^{(m)}) \rightarrow \Gamma ( V ',j _* \D _U ^{(m)})$ est plat.
  On conclut via \cite[3.1.1]{Be1}.
\hfill \hfill \qed \end{proof}

\begin{vide}
Pour tous entiers $n$ et $n '$, on notera $\D ^{(m)}_{U , n}\cdot \D ^{(m)}_{U , n'}$
l'image de l'homomorphisme $\O _U$-linéaire à droite et à gauche
$\D ^{(m)}_{U , n}\otimes _{\O _U }  \D ^{(m)}_{U , n'} \rightarrow
\D ^{(m)}_{U , n+ n'}.$

\end{vide}

\begin{prop}
\label{Dabonnefilt}
Pour $n<0$, on pose $\D ^{(m)}_{U , n}:=0$.
Pour tout couple $(r,s)\in \N ^2$,
$\sum _{j=0,...,p ^m -1} \D ^{(m)}_{U , r -j} \cdot \D ^{(m)}_{U , s+j}
= \D ^{(m)}_{U , r+s}  $.
\end{prop}
\begin{proof}
  Même calcul que \cite[7.1.2]{caro}.
\hfill \hfill \qed \end{proof}

\begin{defi}
\label{deffiltration}
Soit $\M$ un $\D ^{(m)} _{U}$-module.
Une {\it filtration} de $\M$ est une famille $(\M _r) _{r\in \N}$ de sous-$\O _U$-modules de $\M$ telle que :

 1. Pour tous $r,\ s\in \N$ : $\M _r \subset \M _{r+1}$, $\D ^{(m)} _{U , r} \cdot \M _s \subset \M _{r+s}$,

 2. $\M = \cup _{r\in \N} \M _r$.
\end{defi}

\begin{defi}
\label{defbonnefilt}
Soit $\M$ un $\D ^{(m)} _{U}$-module à gauche muni d'une filtration. La filtration est {\it bonne}
si et seulement si :
\begin{enumerate}
 \item  Quel que soit $r\in \N $, $\M _r$ est $\O _U$-cohérent ;

 \item Il existe un entier $r_1 \in \N$ tel que pour tout entier $r\geq r_1$, on ait
$$\M _r =\sum_{j=0} ^{p^m -1} \D ^{(m)} _{U, r-r_1 +j}\cdot \M_{r _1 -j} .$$
\end{enumerate}
\end{defi}

D'après la proposition \ref{Dabonnefilt},
la famille $(\D ^{(m)} _{U, r}) _{r\in \N}$ est une filtration de
$\D ^{(m)} _{U}$ vérifiant la condition $2$ pour tout
entier $r_1$. La condition $1$ étant aussi vérifiée, cette filtration est donc bonne.
On l'appellera la {\it filtration par l'ordre}.

\begin{prop}
\label{globpfbfil}
Un $\D ^{(m)} _{U}$-module globalement de présentation finie admet une bonne filtration.
\end{prop}
\begin{proof}
  Analogue à \cite[7.1.5]{caro}.
\hfill \hfill \qed \end{proof}

\begin{theo}[Théorème B]\label{theoB}
  On suppose $U$ affine et soit $\M$ un $\D ^{(m)} _U$-module globalement de présentation finie.
  Alors, pour tout entier $i \neq 0$, $H ^i (U , \M) =0$.
\end{theo}
\begin{proof}
  Cela résulte de \ref{globpfbfil} et du fait que,
  pour tout entier $i \neq 0$, le foncteur $H ^i (U , \M)$ commute aux limites inductives filtrantes et
  du théorème de type $B$ pour les $\O _U$-modules cohérents.
\hfill \hfill \qed \end{proof}

\begin{theo}[Théorème A]\label{theoA}
  Lorsque $U$ est affine, les foncteurs $\M \mapsto \Gamma (U, \M)$ et
  $M \mapsto \widetilde{M}$ sont des équivalences
  quasi-inverses entre la catégorie des $\D ^{(m)} _U $-modules globalement de présentation finie et
  celle des \linebreak g$\Gamma (U,\D ^{(m)} _U )$-modules de type fini.
\end{theo}
\begin{proof}
  Il s'agit de calquer \cite[7.1.8]{caro}.
\hfill \hfill \qed \end{proof}

\begin{rema}
Comme le faisceau d'anneaux $\D ^{(m)} _{U } $ est cohérent
(\ref{wildetildeDcoh}),
un $\D ^{(m)} _{U } $-module $\M$ est cohérent
si et seulement s'il est localement de présentation finie,
i.e., d'après \ref{theoA}, si et seulement s'il est localement de la forme
$\widetilde{M}$, où $M$ est un $\Gamma (U,\D ^{(m)} _{U })$-module de type fini.
\end{rema}

\begin{theo}
Soit $\M$ un $\D ^{(m)} _{U } $-module.
Alors $\M$ est cohérent si et seulement s'il admet localement de bonnes filtrations.
\end{theo}
\begin{proof}
  Similaire à \cite[7.1.10]{caro}.
\hfill \hfill \qed \end{proof}

\begin{prop}\label{j*->Rj*iso}
  Soient $j$ : $U \subset P$ une immersion ouverte de $\V$-schémas formels lisses
  telle que $P _0 \setminus U _0$ soit le support d'un diviseur
  et $\M$ un $\D ^{(m)} _U$-module localement en $P$ de présentation finie
  (i.e. il existe un recouvrement d'ouverts de $P = \cup _\alpha P _\alpha$ tel que, pour tout $\alpha$,
  $\M |_{U \cap P _\alpha}$ soit globalement de présentation finie). Alors :
  \begin{enumerate}
    \item Le morphisme canonique $j _* \M \rightarrow \R j_* \M$ est un isomorphisme ;
    \item Pour tout morphisme lisse $P \rightarrow P'$ de $\V$-schémas formels faibles,
    pour tout ouvert $U'$ de $P'$ tel que $U = f ^{-1} (U')$,
    le morphisme canonique
    $j_* (\Omega ^\bullet _{U/U'} \otimes _{\O _U}  \M) \rightarrow
    \R j_* ( \Omega ^\bullet _{U/U'} \otimes _{\O _U} \M)$
    est un isomorphisme ;
    \item Le $j _* \D ^{(m)} _U$-module $j _* \M$ est cohérent. Plus précisément, pour tout
    $\D ^{(m)} _U$-module globalement de présentation finie $\NN$, $j _* \NN$ est
    un $j _* \D ^{(m)} _U$-module globalement de présentation finie.
  \end{enumerate}
\end{prop}
\begin{proof}
  Prouvons $1)$. Comme $P _0 \setminus U _0$ est un diviseur, l'assertion étant locale, on peut supposer
  $P _\alpha$ et $U \cap P _\alpha$ affines.
  Pour tout entier $i$, $(\mathcal{H} ^i j _* \M) |_{P _\alpha} \riso \mathcal{H} ^i j _{\alpha *} ( \M |_{U \cap P _{\alpha}})$,
  où $j _\alpha$ : $ U \cap P _{\alpha } \hookrightarrow P _\alpha$. Or,
$ \mathcal{H} ^i j _{\alpha *} ( \M _{U \cap P _{\alpha}})$ est le faisceau associé au préfaisceau qui à tout ouvert
principal $P' $ de $P _\alpha$ associe $H ^i ( U \cap P ',\M)$.
Le théorème $B$ \ref{theoB} nous permet de conclure $1)$.
L'assertion $2)$ découle de $1)$ et du fait que $\Omega _{U/U'}$ est un $\O _{U}$-module
localement en $P$ libre de type fini
(on prend une base d'ouverts de $P$ ayant des coordonnées locales au dessus de $P'$).
Traitons à présent l'assertion $3)$. On dispose d'une présentation finie :
$( \D _{U}^{(m)}) ^r
\overset{\phi}{\rightarrow}( \D _{U}^{(m)}) ^s \overset{\epsilon}{\rightarrow}\NN  \rightarrow 0$.
Or, il résulte de $1)$ et de \ref{coro-ker-im-coker}.(i)
que le foncteur $j _*$ de la catégorie des
$\D ^{(m)} _U$-modules localement en $P$ de présentation finie dans celle des
$j _* \D ^{(m)} _U$-modules est exact.
On en déduit la suite exacte
$j _*( \D _{U}^{(m)}) ^r
\overset{j _*\phi}{\rightarrow}j _*( \D _{U}^{(m)}) ^s \overset{j _*\epsilon}{\rightarrow}j _*\NN  \rightarrow 0$.

\hfill \hfill \qed \end{proof}

\begin{nota}
Dans la suite de cette section $U$ sera supposé affine et lisse.
\end{nota}

\begin{defi}\label{defiindfini}
  Soit $M$ un $\Gamma (U ,\O_U)$-module
(resp. $\Gamma (U,\D ^{(m)} _U)$-module).
Le $\O _U$-module (resp. $\D ^{(m)} _U$-module) {\it associé à $M$}
  est $\widetilde{M} := \O _U \otimes _{\Gamma (U, \O _U)} M$
  (resp. $\D ^{(m)} _U \otimes _{\Gamma (U,\D ^{(m)} _U)} M$).
    Un tel module est dit {\it $\mathcal{I} nd $-fini}.

  Si $u$ : $M '\rightarrow M$ est un homomorphisme de
  $\Gamma (U,\O _U)$-modules (resp. $\Gamma (U,\D ^{(m)} _U)$-modules), on note
  $\widetilde{u}$ : $\smash{\widetilde{M}} ' \rightarrow \widetilde{M}$ l'homomorphisme canonique induit.
 On vérifie enfin que l'on obtient un foncteur $\sim$ : $M \mapsto \widetilde{M}$ exact.
\end{defi}

\begin{rema}\label{remapreeqcatindfini}
$\bullet$  Le faisceau $\widehat{\D} ^{(m)} _{\U}$ n'est pas $\mathcal{I} nd$-fini.

$\bullet$ De plus, pour qu'un $\D ^{(m)} _U$-module soit $\mathcal{I}nd$-fini,
  il faut et il suffit que la structure sous-jacente de $\O _U$-module soit $\mathcal{I}nd$-fini.
  En effet, comme $\D ^{(m)} _U$ est limite inductive filtrante de $\O _U$-modules cohérents et
  puisque les foncteurs $\O _U \otimes _{\Gamma (U, \O _U)} -$ et $\Gamma (U,-)$
  commutent aux limites inductives filtrantes,
  on bénéficie de $\D ^{(m)} _U \riso \O _U \otimes _{\Gamma (U, \O _U)} \Gamma (U, \D ^{(m)} _U)$.
\end{rema}

\begin{theo}[Théorème B]\label{theoB'}
Pour tout $\D ^{(m)} _U$-module $\M$ $\mathcal{I} nd$-fini,
pour tout entier $i \neq 0$, $H ^i (U , \M) =0$.
\end{theo}
\begin{proof}
 Cela résulte de \ref{theoB}, \ref{rema-indfil}.(i) et du fait que,
  pour tout entier $i \neq 0$, le foncteur $H ^i (U , \M)$ commute aux limites inductives filtrantes.
\hfill \hfill \qed \end{proof}

\begin{prop}\label{coro-ker-im-coker}

  (i) Soit $u$ : $M \rightarrow N$ un morphisme de $\Gamma (U,\D ^{(m)} _U )$-modules.
  les $\D ^{(m)} _U $-modules associés à $\mathrm{Ker}  u$,
  $\mathrm{Im}  u$, $\mathrm{Coker}  u$, sont respectivement
  $\mathrm{Ker}  \tilde{u}$,
  $\mathrm{Im}  \tilde{u}$, $\mathrm{Coker}  \tilde{u}$.
  En particulier, lorsque $M$ et $N$ sont de type fini, ces $\D ^{(m)} _U$-modules
  sont globalement de présentation finie.

  (ii) Si $M$ est une limite inductive (resp. somme directe) d'une famille de
  $\Gamma (U,\D ^{(m)} _U )$-module ($M _\lambda$), $\widetilde{M}$ est limite inductive
  (resp. somme directe) de la famille ($\widetilde{M _\lambda}$), à isomorphisme canonique près.

  (iii) Si $M$ et $N$ sont deux $\Gamma (U,\D ^{(m)} _U )$-modules,
  les $\D ^{(m)} _U $-modules associés à
  $M \otimes _{\Gamma (U,\O _U)} N$ et $\mathrm{Hom} _{\Gamma (U,\O _U)} (M,N)$ sont respectivement
  $\widetilde{M} \otimes _{\O _U} \widetilde{N}$ et $\mathcal{H} om _{\O _U} (\widetilde{M}, \widetilde{N})$.

  (iv) Le foncteur $\Gamma (U,- )$ est exact sur la catégorie des $\D ^{(m)} _U$-modules $\mathcal{I}nd$-finis.

  On bénéficie de résultats identiques avec $\O _U$ à la place de $\D ^{(m)} _U$.
\end{prop}
\begin{proof}
  La partie (i) se prouve comme \cite[1.3.9]{EGAI}. L'assertion (ii) résulte de la commutation aux limites
  inductives du produit tensoriel tandis que (iii) est aisé.
  Enfin, (iv) résulte de (i) et du théorème de
  type $B$ pour les $\D ^{(m)} _U$-modules $\mathcal{I}nd$-finis (\ref{theoB'}).
\hfill \hfill \qed \end{proof}

\begin{rema}\label{rema-indfil}
  (i) Un $\D ^{(m)} _U $-module $\mathcal{I}nd$-fini
  est limite inductive filtrante de ses sous-$\D ^{(m)} _U $-modules globalement de présentation finie.
  Ainsi, via \ref{coro-ker-im-coker}.(ii), la catégorie des $\D ^{(m)} _U $-modules $\mathcal{I}nd$-finis
  est la plus petite catégorie stable par limite inductive et contenant les
  $\D ^{(m)} _U $-modules globalement de présentation finie.

  (ii) Soient $M$ un $\Gamma ( U ,\D ^{(m)} _U)$-module et $U'$ un ouvert principal de $U$.
  Comme le produit tensoriel et $\Gamma (U', - )$ commutent aux limites inductives filtrantes,
  il résulte de la remarque (i) et du théorème de type $A$ (\ref{theoA}) appliqué à $U$ et $U'$ que
$\Gamma (U', \widetilde{M}) \riso \Gamma (U', \D ^{(m)} _U ) \otimes _{\Gamma (U,\D ^{(m)} _U )} M$.
\end{rema}

\begin{prop}\label{eqcatindfini}
  Les foncteurs $M \mapsto \widetilde{M}$ et $\M \mapsto \Gamma (U, \M)$
  induisent des équivalences
  quasi-inverses entre la catégorie des $\Gamma (U,\D ^{(m)} _U )$-modules
  et celle des $\D ^{(m)} _U $-modules $\mathcal{I}nd$-finis.
  On dispose de résultats identiques en remplaçant $\D ^{(m)} _U $ par $\O _U$.
\end{prop}
\begin{proof}
 Soit $M$ un $\Gamma (U,\D ^{(m)} _U )$-module.
 D'après la remarque ci-dessus \ref{rema-indfil}.(ii),
 $\Gamma ( U, \widetilde{M}) \riso M$. Réciproquement, pour tout
 un $\D ^{(m)} _U $-module $\M $ isomorphe à $\widetilde{M}$,
 $\Gamma (U, \M ) ^{\sim} \riso \Gamma (U, \widetilde{M} ) ^{\sim} \riso \widetilde{M} \riso \M$.

\hfill \hfill \qed \end{proof}

\begin{prop}
  Soit $0 \rightarrow \E '\rightarrow \E \rightarrow \E'' \rightarrow 0$ une suite exacte de
  $\D ^{(m)} _U $-modules. Si deux de ces modules sont $\mathcal{I}nd$-finis, alors le troisième l'est aussi.
\end{prop}
\begin{proof}
  Le cas où $\E$ est $\mathcal{I}nd$-fini se déduit de \ref{coro-ker-im-coker}.(i).

  Supposons à présent $\E '$ et $\E ''$ de la forme $\widetilde{E} '$ et $\widetilde{E} ''$.
  Via le théorème $B$ pour les
 $\D ^{(m)} _U $-modules $\mathcal{I}nd$-finis,
$0 \rightarrow E ' \rightarrow \Gamma (U,\E) \rightarrow E '' \rightarrow 0$ est exacte. En lui
appliquant le foncteur exact $\sim$, on conclut que
l'homomorphisme canonique $\Gamma (U, \E) ^{\sim} \rightarrow  \E$
 est un isomorphisme.
\hfill \hfill \qed \end{proof}

\subsection{Opérations cohomologiques, lien entre les cas formel faible et formel}

\begin{vide}
  Soient $g$ : $U ' \rightarrow U$ et $g' $ : $U'' \rightarrow U'$ deux morphismes de $\V$-schémas formels faibles.
  Le faisceau $ g ^* \D ^{(m)} _{U}$ est muni d'une structure canonique de
  $(\D ^{(m)} _{U'}, g ^{-1} \D ^{(m)} _U)$-bimodule, que l'on notera $\D ^{(m)} _{U ' \rightarrow U}$.
  L'image inverse extraordinaire par $g$ d'un complexe $\E \in D ^- ( \D ^{(m)} _U)$ est définie en posant
  $g ^{!^{(m)}} (\E) := \D ^{(m)} _{U ' \rightarrow U} \otimes ^\L _{g ^{-1} \D ^{(m)} _U} g ^{-1} \E [d _{U'/U}]$,
  où $d _{U'/U}$ est la dimension relative de $U'$ sur $U$.
Si aucune confusion sur le niveau n'est à craindre, on notera $g ^!$ à la place.
  On vérifie que l'on dispose d'un isomorphisme
  $ g ^{\prime !} \circ g ^ ! (\E)\riso (g \circ g ') ^! (\E)$ fonctoriel en $\E$
  et associatif.

  De plus, on notera $\D ^{(m)} _{U  \leftarrow U'} :=
  \omega _{U'} \otimes g ^* _g ( \D ^{(m)} _U \otimes _{\O _U} \omega ^{-1} _U)$, l'indice $g$ signifiant
que l'on prend la structure gauche de $\D ^{(m)} _U $-module à gauche.
  Comme $g _0$ : $U ' _0 \rightarrow U _0 $ est un morphisme quasi-séparé et quasi-compact entre
  schémas noethériens de dimension de Krull finie, le foncteur $g _{0*} = g _*$ est de dimension
  cohomologique finie.
  On définit l'image directe par $g$ de $\E '\in D ^- ( \D ^{(m)} _{U'})$ en posant
  $g _{+^{(m)}} (\E') := \R g _* (\D ^{(m)} _{U  \leftarrow U'} \otimes _{\D ^{(m)} _{U'}} ^\L \E ')$.
On pourra aussi le noter $g _+$.
  De manière analogue au cas des schémas, on vérifie l'isomorphisme
  $g  _+ \circ g ' _+ (\E '') \riso (g \circ g') _+ (\E'')$ fonctoriel en
  $\E '' \in D ^- ( \D ^{(m)} _{U''})$ et associatif.

Soient $\I \subset \O _U$ un idéal cohérent, $Z$ le sous-schéma formel faible fermé défini par $\I$ et $\PP ^n _{(m)}(\I)$
l'enveloppe à puissance divisée de niveau $m$ et d'ordre $n$ de $\I$.
Pour tout $\D ^{(m)} _{U}$-module $\E $, le faisceau, défini en posant
$\underline{\Gamma} ^{(m)} _Z ( \E) :=
\underset{\underset{n}{\longrightarrow}}{\lim} \mathcal{H} om _{\O _U} (\PP ^n _{(m)}(\I), \E)$,
est muni d'une structure canonique de $\D ^{(m)} _{U}$-module. Ce foncteur
s'étend sur $D ^+ ( \D ^{(m)} _{U})$ et sera noté $\R \underline{\Gamma} ^{(m)} _Z$.
Cette construction est fonctorielle
en $Z$, i.e., si $\I ' \subset \I \subset \O _U$ sont des idéaux, $Z'$ et $Z$
les $\V$-schémas formels faibles correspondants,
on dispose d'un morphisme
$\R \underline{\Gamma} ^{(m)} _{Z} (\E) \rightarrow \R \underline{\Gamma} ^{(m)} _{Z'} (\E)$
fonctoriel en $\E \in D ^+ ( \D ^{(m)} _{U})$.

Lorsque $g$ est une immersion fermée, on dispose comme dans le cas des schémas d'un isomorphisme
canonique
\begin{equation}
  \label{g+g!gamma}
  g _+ g ^! (\E) \riso \R \underline{\Gamma} ^{(m)} _{U'} (\E)
\end{equation}
fonctoriel en $\E \in D ^\mathrm{b} ( \D ^{(m)} _{U})$.
De plus, on bénéficie aussi, via un calcul analogue à la situation des schémas, d'un isomorphisme canonique
$\E' \riso g ^! g _+  (\E')$ fonctoriel $\E '\in D ^\mathrm{b} ( \D ^{(m)} _{U'})$.

\end{vide}

\begin{vide}
  [Complexes quasi-cohérents]\label{qcUdag}
  Soient $\E $ un objet de $ D ^-(\overset{^g}{}\D _{U} ^{(m)})$ et $\FF \in D ^-(\D _{U} ^{(m)}\overset{^d}{})$.
  On pose $\E _i := \D ^{(m)} _{U _i} \otimes ^\L _{\D ^{(m)} _U} \E$,
  $\FF _i := \FF \otimes ^\L _{\D ^{(m)} _U} \D ^{(m)} _{U _i}$ et
  $$\FF \smash{\widehat{\otimes}} ^\L_{\D _{U} ^{(m)}} \E :=
  \R \underset{\longleftarrow}{\lim} _i
  \FF _i \otimes ^\L _{\D ^{(m)} _{U _i}} \E _i,$$
  avec $i$ un entier (notations \ref{defschff}).
  De même en remplaçant $\D ^{(m)}_U$ par $\O _U$ ou $\widehat{\D} ^{(m)}_U$ ou plus généralement un
  faisceau d'anneaux sur $U$. En résolvant $\E$ et $\FF$ platement, on construit un morphisme
  canonique $\FF \otimes ^\L_{\D _{U} ^{(m)}} \E  \rightarrow
  \FF \smash{\widehat{\otimes}} ^\L_{\D _{U} ^{(m)}} \E$ bifonctoriel en $\E$ et $\FF$ (i.e., un cube est commutatif).
  Lorsque
  $\FF = \D ^{(m)}_U$, on notera $\widehat{\E} := \D ^{(m)}_U \smash{\widehat{\otimes}} ^\L_{\D _{U} ^{(m)}} \E $.
  On remarque que le morphisme canonique
  $\O  _{U} \smash{\widehat{\otimes}} ^\L_{\O  _{U} } \E \rightarrow
  \D _{U} ^{(m)} \smash{\widehat{\otimes}} ^\L_{\D _{U} ^{(m)}} \E$ est un isomorphisme.
  On définit $D ^\mathrm{b} _\mathrm{qc}(\overset{^*}{}\D _{U} ^{(m)})$, avec $* =g $ ou $*=d$,
la sous-catégorie pleine de $D ^\mathrm{b} (\overset{^*}{}\D _{U} ^{(m)})$, des complexes
  $\E $ tels que, pour tout $i$, $\E _i \in D ^\mathrm{b} _\mathrm{qc} (\overset{^*}{}\D _{U _i} ^{(m)})$.
  Ses objets seront appelés {\it complexes quasi-cohérents}.
  Par exemple, un $\D _{U} ^{(m)}$-module $\mathcal{I}nd$-fini (\ref{defiindfini})
  est un complexe quasi-cohérent.
  Par \cite[3.2.2]{Beintro2},
  on dispose d'un foncteur $D ^\mathrm{b} _\mathrm{qc}(\overset{^*}{}\D _{U} ^{(m)})
  \rightarrow D ^\mathrm{b} _\mathrm{qc}(\overset{^*}{}\widehat{\D} _{\U} ^{(m)})$,
  qui à $\E $ associe $\widehat{\E}$.
  De plus, si $\E \in D ^\mathrm{b} _\mathrm{qc}(\overset{^*}{}\D _{U} ^{(m)})$
  et $\FF \in D ^\mathrm{b} _\mathrm{qc}(\overset{^*}{}\widehat{\D} _{\U} ^{(m)})$,
  tout morphisme $\E \rightarrow \FF$ se factorise de manière unique
  par $\widehat{\E} \rightarrow \FF$.
\end{vide}

\begin{lemm}\label{compDiso}
  Soient $U$ un $\V$-schéma formel faible lisse et
  $\E\in D ^\mathrm{b} _\mathrm{coh}(\D _{U} ^{(m)})$.
  Le morphisme canonique :
  $\widehat{\D} ^{(m)} _{\U }   \otimes _{ \D ^{(m)} _{U}} \E \rightarrow \widehat{\E}$
  est un isomorphisme de $D ^\mathrm{b} _\mathrm{coh}(\widehat{\D} _{\U} ^{(m)})$.
\end{lemm}
\begin{proof}
Cela résulte de \cite[3.2.3]{Beintro2}.
\hfill \hfill \qed \end{proof}

\begin{lemm}\label{globpfinvlis}
  Soit $g$ : $ U ' \rightarrow U$ un morphisme lisse de $\V$-schémas formels faibles affines et lisses.
  Pour tout $\D ^{(m)} _U$-module globalement de présentation finie $\E$,
  $\D ^{(m)} _{U' \rightarrow U} \otimes _{ g ^{-1} \D ^{(m)} _U} g ^{-1} \E$ est un
  $\D ^{(m)} _{U'}$-module globalement de présentation finie.
\end{lemm}
\begin{proof}
Comme $\O _{U'} \otimes _{ g ^{-1} \O _U} g ^{-1} \E
\rightarrow \D ^{(m)} _{U' \rightarrow U} \otimes _{ g ^{-1} \D ^{(m)} _U} g ^{-1} \E$
est un isomorphisme, ceux-ci sont $\mathcal{I}nd$-finis. Via \ref{theoA}, il suffit de prouver
que $\Gamma  (U',\D ^{(m)} _{U' \rightarrow U}) \otimes _{\Gamma(U,\D ^{(m)} _U)} \Gamma(U,\E)$
est de type fini sur $\Gamma(U',\D ^{(m)} _{U'})$. Or, par des calculs en coordonnées locales,
le morphisme
$\D ^{(m)} _{U'}\rightarrow \D ^{(m)} _{U' \rightarrow U}$
est surjectif. Par \ref{coro-ker-im-coker}.(i) et \ref{eqcatindfini},
$\Gamma ( U', \D ^{(m)} _{U'})\twoheadrightarrow \Gamma ( U',  \D ^{(m)} _{U' \rightarrow U})$.

\hfill \hfill \qed \end{proof}

\begin{lemm}\label{prop-comm-^inv}
  Soient $g$ : $U' \rightarrow U$ un morphisme de $\V$-schémas formels faibles lisses et
  $\E\in D ^\mathrm{b} _\mathrm{qc} (\D _{U} ^{(m)})$ (\ref{qcUdag}). Le morphisme canonique
\begin{equation}
  \label{prop-comm-^invdiag}
  \widehat{\D} ^{(m)} _{\U '}   \smash{\widehat{\otimes}} ^\L_{\D _{U'} ^{(m)}} g ^! (\E)
\rightarrow
\hat{g} ^! (\widehat{\D} ^{(m)} _{\U }   \smash{\widehat{\otimes}} ^\L_{\D _{U} ^{(m)}} \E )
\end{equation}
est un isomorphisme. De plus, celui-ci est transitif en $g$.
En particulier, si $\E\in D ^\mathrm{b} _\mathrm{coh}(\D _{U} ^{(m)})$ et
  $g ^! (\E)\in D ^\mathrm{b} _\mathrm{coh}(\D _{U'} ^{(m)})$ (ce qui est le cas si $g$ est lisse \ref{globpfinvlis}),
  le morphisme canonique
$\widehat{\D} ^{(m)} _{\U '}   \otimes _{\D _{U'} ^{(m)}} g ^! (\E)
\rightarrow
\hat{g} ^! (\widehat{\D} ^{(m)} _{\U } \otimes _{\D _{U} ^{(m)}} \E )$
est un isomorphisme.
\end{lemm}
\begin{proof}
  Cela découle de \ref{compDiso} et de la commutation au changement de base de l'image inverse
  extraordinaire.
\hfill \hfill \qed \end{proof}

\begin{prop}\label{prop-comm-^inv1}
  Soient $m _0$ un entier, $g$ : $U' \rightarrow U$ un morphisme de $\V$-schémas formels faibles lisses et
  $\E ^{(m _0)}\in D ^\mathrm{b} _\mathrm{coh} (\D _{U} ^{(m _0)})$.
  Pour tout entier $m \geq m_0 $, on note $\E ^{(m )} : = \D _{U} ^{(m )}\otimes _{\D _{U} ^{(m _0)}} \E ^{(m _0)}$.
  Lorsque $g ^{!(m _0)} (\E^{(m _0)})\in D ^\mathrm{b} _\mathrm{coh} (\D _{U'} ^{(m _0)})$
  et $g ^{!(m )} (\E^{(m )})\in D ^\mathrm{b} _\mathrm{coh} (\D _{U'} ^{(m )})$,
  le morphisme canonique
\begin{equation}
  \label{prop-comm-^inv1diag1}
  \widehat{\D} ^{(m)} _{\U ',\Q}   \otimes _{\widehat{\D} ^{(m _0)} _{\U ',\Q}}
  \hat{g} ^{!(m _0)} (\widehat{\E} ^{(m_0)} _\Q)
\rightarrow
\hat{g} ^{!(m )} ( \widehat{\E} ^{(m)} _\Q )
\end{equation}
est un isomorphisme.

Lorsque, pour tout $m\geq m_0$, $g ^{!(m )} (\E^{(m )})\in D ^\mathrm{b} _\mathrm{coh} (\D _{U'} ^{(m)})$,
le morphisme
\begin{equation}
  \label{prop-comm-^inv1diag2}
  \D ^\dag _{\U ',\Q}   \otimes _{\widehat{\D} ^{(m _0)} _{\U ',\Q}}
  \hat{g} ^{!(m _0)} (\widehat{\E} ^{(m_0)} _\Q)
\rightarrow
\hat{g} ^{!\dag} (\D ^\dag _{\U ,\Q}   \otimes _{\widehat{\D} ^{(m _0)} _{\U ,\Q}}
\widehat{\E} ^{(m _0)} _\Q )
\end{equation}
est un isomorphisme.
De plus, ceux-ci sont transitifs en $g$.
\end{prop}
\begin{proof}
Notons $\E := \E ^{(m)} _\Q$ ($\E$ est indépendant de $m$ à isomorphisme canonique près).
Si $m$ est un entier tel que $g ^! (\E^{(m )})\in D ^\mathrm{b} _\mathrm{coh} (\D _{U'} ^{(m )})$,
il découle de \ref{prop-comm-^inv} que
le morphisme canonique
$\widehat{\D} ^{(m)} _{\U ',\Q}   \otimes _{\D  _{\U ',\Q}} g ^! (\E) \rightarrow
  \hat{g} ^! (\widehat{\E} ^{(m)} _\Q)$
  est un isomorphisme. Il en résulte que \ref{prop-comm-^inv1diag1} est un isomorphisme.
  Par passage à la limite sur le niveau, on en conclut \ref{prop-comm-^inv1diag2}.
\hfill \hfill \qed \end{proof}

\begin{prop}\label{prop-comm-^inv2}
  Soient $g$ : $U' \rightarrow U$ un morphisme lisse de $\V$-schémas formels faibles lisses et
  $\E\in D ^\mathrm{b} _\mathrm{coh}(\D _{U} ^{(m)})$. Le morphisme canonique
\begin{equation}
  \label{prop-comm-^invdiag2}
  \D ^\dag _{\U ',\Q}   \otimes _{ \D ^{(m)} _{U '}} g ^{!^{(m)}} (\E)
\rightarrow
g ^{!\dag} (\D ^\dag _{\U ,\Q}   \otimes _{ \D ^{(m)} _{U }} \E )
\end{equation}
est un isomorphisme. De plus, celui-ci est transitif en $g$.
\end{prop}
\begin{proof}
  Cela découle de \ref{globpfinvlis}, \ref{prop-comm-^inv} et \ref{prop-comm-^inv1}.
\hfill \hfill \qed \end{proof}

\begin{lemm}\label{prop-comm-^}
  Soient $g$ : $U' \rightarrow U$ un morphisme de $\V$-schémas formels faibles lisses et
  $\E'\in D ^\mathrm{b} _\mathrm{qc} (\D _{U'} ^{(m)})$ (\ref{qcUdag}). Le morphisme canonique
\begin{equation}
  \label{prop-comm-^diag}
  \widehat{\D} ^{(m)} _{\U }   \smash{\widehat{\otimes}} ^\L _{ \D ^{(m)} _{U }} g _{+^{(m)}} (\E')
\rightarrow
\hat{g} _{+^{(m)}} (\widehat{\D} ^{(m)} _{\U '}   \smash{\widehat{\otimes}} ^\L _{ \D ^{(m)} _{U '}} \E' )
\end{equation}
est un isomorphisme. De plus, celui-ci est transitif en $g$.
En particulier, si $\E'\in D ^\mathrm{b} _\mathrm{coh}(\D _{U'} ^{(m)})$ et
  $g _{+^{(m)}} (\E')\in D ^\mathrm{b} _\mathrm{coh}(\D _{U} ^{(m)})$
le morphisme canonique
$\widehat{\D} ^{(m)} _{\U }   \otimes _{ \D ^{(m)} _{U }} g _{+^{(m)}} (\E')
\rightarrow
\hat{g} _{+^{(m)}} (\widehat{\D} ^{(m)} _{\U '}   \otimes _{ \D ^{(m)} _{U '}} \E' )$ est un isomorphisme.
\end{lemm}
\begin{proof}
Via la suite exacte
  $0\rightarrow \O _U \overset{\pi ^{i+1}}{\longrightarrow}
  \O _U \rightarrow \O _{U _i} \rightarrow 0$, on voit que $\O _{U _i}$ est de Tor-dimension finie sur $\O_U$.
  On bénéficie donc de l'isomorphisme de projection :
$\O _{U _i} \otimes ^\L _{\O _{U}} \R g _* (\FF')
\riso
\R g _* ( g ^{-1} \O _{U_i} \otimes ^\L _{g ^{-1} \O _U } \FF ')$
fonctoriel en $\FF ' \in D (g ^{-1}  \D ^{(m)} _U)$.
D'où le diagramme commutatif :
\begin{equation}
  \label{prop-comm-^diag2}
  \xymatrix @C=0,2cm @R=0,3cm {
{\O _{U _i} \otimes ^\L _{\O _{U}}
\R g _* ( \D ^{(m)} _{U \leftarrow U '} \otimes ^\L _{\D ^{(m)} _{U '}} \E ') }
\ar[d]
\ar[r] ^-\sim
&
{\R g _* ( g ^{\text{-}1} \O _{U_i}
\smash{\underset{g ^{\text{-}1} \O _U }{\otimes ^\L}}
( \D ^{(m)} _{U \leftarrow U '} \otimes ^\L _{\D ^{(m)} _{U '}} \E '   ))}
\ar[d]
\\
{\O _{U _i} \otimes ^\L _{\O _{U}}
\R g _* ( \widehat{\D} ^{(m)} _{\U \leftarrow \U '} \otimes ^\L _{\D ^{(m)} _{U '}} \E ') }
\ar[r] ^-\sim
\ar[d]
&
{\R g _* ( g ^{\text{-}1} \O _{U_i}
\smash{\underset{g ^{\text{-}1} \O _U }{\otimes ^\L}}
( \widehat{\D} ^{(m)} _{\U \leftarrow \U '} \otimes ^\L _{\D ^{(m)} _{U '}} \E '   ))}
\ar[d]
\\
{\O _{U _i} \otimes ^\L _{\O _{U}}
\R g _* ( \widehat{\D} ^{(m)} _{\U \leftarrow \U '} \smash{\widehat{\otimes}} ^\L _{\D ^{(m)} _{U '}} \E ') }
\ar[r] ^-\sim
&
{\R g _* ( g ^{\text{-}1} \O _{U_i}
\smash{\underset{g ^{\text{-}1} \O _U }{\otimes ^\L}}
(   \widehat{\D} ^{(m)} _{\U \leftarrow \U '} \smash{\widehat{\otimes}} ^\L _{\D ^{(m)} _{U '}} \E '  ))}
}
\end{equation}

De plus, on dispose des morphismes
$\D ^{(m)} _{U \leftarrow U '} \rightarrow \widehat{\D} ^{(m)} _{\U \leftarrow \U '}
\rightarrow
\D ^{(m)} _{U _i \leftarrow U ' _i }$
qui induisent
$g ^{-1} \O _{U_i} \otimes ^\L _{g ^{-1} \O _U }
\D ^{(m)} _{U \leftarrow U '}
\riso g ^{-1} \O _{U_i} \otimes ^\L _{g ^{-1} \O _U }
\widehat{\D} ^{(m)} _{\U \leftarrow \U '}
\riso   \D ^{(m)} _{U _i \leftarrow U ' _i } $.
Il en découle le diagramme commutatif :
\begin{equation}
  \label{prop-comm-^diag3}
  \xymatrix @R=0,3cm {
{\R g _* ( g ^{-1} \O _{U_i} \otimes ^\L _{g ^{-1} \O _U }
( \D ^{(m)} _{U \leftarrow U '} \otimes ^\L _{\D ^{(m)} _{U '}} \E '   ))}
\ar[r] ^-\sim \ar[d]
&
{\R g _* (\D ^{(m)} _{U _i \leftarrow U ' _i } \otimes ^\L _{\D ^{(m)} _{U ' _i}} \E '_i )}
\ar@{=}[d]
\\
{\R g _* ( g ^{-1} \O _{U_i} \otimes ^\L _{g ^{-1} \O _U }
( \widehat{\D} ^{(m)} _{\U \leftarrow \U '} \otimes ^\L _{\D ^{(m)} _{U '}} \E '   ))}
\ar[r] ^-\sim \ar[d]
&
{\R g _* (\D ^{(m)} _{U _i \leftarrow U ' _i } \otimes ^\L _{\D ^{(m)} _{U ' _i}} \E '_i )}
\ar@{=}[d]
\\
{\R g _* ( g ^{-1} \O _{U_i} \otimes ^\L _{g ^{-1} \O _U }
(   \widehat{\D} ^{(m)} _{\U \leftarrow \U '} \smash{\widehat{\otimes}} ^\L _{\D ^{(m)} _{U '}} \E '  ))}
\ar[r] ^-\sim
&
{\R g _* (\D ^{(m)} _{U _i \leftarrow U ' _i } \otimes ^\L _{\D ^{(m)} _{U ' _i}} \E '_i ).}
}
\end{equation}
En composant les isomorphismes du haut de \ref{prop-comm-^diag2} et \ref{prop-comm-^diag3},
on obtient, avec \cite[3.5.2]{Beintro2}, $g _+ (\E') \in D ^\mathrm{b} _\mathrm{qc} (\D ^{(m)} _{U })$.

La flèche \ref{prop-comm-^diag} est alors l'unique morphisme de $D ^\mathrm{b} _\mathrm{qc} (\widehat{\D} ^{(m)} _{\U })$
rendant commutatif le diagramme ci-dessous
  \begin{equation}
  \label{prop-comm-^diag1}
  \xymatrix @R=0,3cm {
{\R g _* ( \D ^{(m)} _{U \leftarrow U '} \otimes ^\L _{\D ^{(m)} _{U '}} \E ') }
\ar[r] \ar[d]
&
{\widehat{\D} ^{(m)} _{\U }   \smash{\widehat{\otimes}} ^\L _{ \D ^{(m)} _{U }} g _+ (\E')}
  \ar@{.>}[dd]
\\
{\R g _* ( \widehat{\D} ^{(m)} _{\U \leftarrow \U '} \otimes ^\L _{\D ^{(m)} _{U '}} \E ') }
  \ar[d]
  &
{ }
\\
{\R g _* ( \widehat{\D} ^{(m)} _{\U \leftarrow \U '}
  \smash{\widehat{\otimes}} ^\L _{\widehat{\D} ^{(m)} _{\U '}} \widehat{\D} ^{(m)} _{\U '}
  \smash{\widehat{\otimes}} ^\L _{\D ^{(m)} _{U '}} \E ') }
  \ar@{=}[r]
&
{\hat{g} _+ (\widehat{\D} ^{(m)} _{\U '}   \smash{\widehat{\otimes}} ^\L _{ \D ^{(m)} _{U '}} \E' ) }
  }
  \end{equation}
  En appliquant $\O _{U _i} \otimes ^\L _{\O _{U}}-$ à \ref{prop-comm-^diag1},
  grâce à \ref{prop-comm-^diag2} et \ref{prop-comm-^diag3}, le composé de gauche devient un isomorphisme.
 Il en résulte par construction que \ref{prop-comm-^diag} est un isomorphisme.
Enfin, la transitivité en $g$ de la construction de \ref{prop-comm-^diag1} est aisée.
\hfill \hfill \qed \end{proof}

\begin{prop}\label{prop-comm-^dir}
  Soient $g$ : $U' \rightarrow U$ un morphisme de $\V$-schémas formels faibles lisses et
  $\E'\in D ^\mathrm{b} _\mathrm{coh}(\D _{U'} ^{(m)})$ tels que
  $g _{+^{(m)}} (\E')\in D ^\mathrm{b} _\mathrm{coh}(\D _{U} ^{(m)})$. Le morphisme canonique
\begin{equation}
  \label{prop-comm-^dirdiag}
  \D ^\dag _{\U }   \otimes _{ \D ^{(m)} _{U }} g _{+^{(m)}} (\E')
\rightarrow
g ^\dag _+ (\D ^\dag _{\U '}   \otimes _{ \D ^{(m)} _{U '}} \E' )
\end{equation}
est un isomorphisme. De plus, celui-ci est transitif en $g$.
\end{prop}
\begin{proof}
  Cela découle de \ref{prop-comm-^} et de \cite[2.4.3]{Beintro2} par complétion et passage à la limite.
\hfill \hfill \qed \end{proof}

\begin{lemm}
  \label{lemmtracecomphat}
  Soient $v$ : $ Y \hookrightarrow U$ une immersion fermée de $\V$-schémas formels faibles lisses
  et $\E \in D ^\mathrm{b} _\mathrm{qc}(\D _{U} ^{(m)})$.
  On dispose d'un isomorphisme canonique
  $\widehat{\D} ^{(m)} _{\U } \smash{\widehat{\otimes}} ^\L _{ \D ^{(m)} _{U }}
  \R \underline{\Gamma} ^{(m)} _{Y} (\E) \riso
  \R \underline{\Gamma} ^{(m)} _{Y _0}
  (\widehat{\D} ^{(m)} _{\U }   \smash{\widehat{\otimes}} ^\L _{ \D ^{(m)} _{U }}\E)$ bifonctoriel
  en $Y$ et $\E $, i.e., si $Y'$ est un sous-schéma fermé de $Y$,
  le diagramme suivant
  $$\xymatrix @R=0,3cm {
  {\widehat{\D} ^{(m)} _{\U }   \smash{\widehat{\otimes}} ^\L _{ \D ^{(m)} _{U }} \R \underline{\Gamma} ^{(m)} _{Y'} (\E)}
  \ar[r]^-\sim \ar[d]
  &
  {\R \underline{\Gamma} ^{(m)} _{Y _0 '} (\widehat{\D} ^{(m)} _{\U }\smash{\widehat{\otimes}} ^\L _{ \D ^{(m)} _{U }}\E)}
  \ar[d]
  \\
  {\widehat{\D} ^{(m)} _{\U }   \smash{\widehat{\otimes}} ^\L _{ \D ^{(m)} _{U }} \R \underline{\Gamma} ^{(m)} _{Y} (\E) }
  \ar[r] ^-\sim
  &
  {\R \underline{\Gamma} ^{(m)} _{Y _0} (\widehat{\D} ^{(m)} _{\U }
\smash{\widehat{\otimes}} ^\L _{ \D ^{(m)} _{U }}\E)}
  }
  $$
 est commutatif et le cube qui s'en déduit par fonctorialité en $\E$ l'est aussi.
De plus, le diagramme canonique
\begin{equation}
  \xymatrix @R=0,3cm {
    {\widehat{\D} ^{(m)} _{\U }   \smash{\widehat{\otimes}} ^\L _{ \D ^{(m)} _{U }} v _+ v ^! (\E)}
    \ar[r] ^-\sim
    \ar[d] ^-\sim
    &
    {\widehat{\D} ^{(m)} _{\U } \smash{\widehat{\otimes}} ^\L _{ \D ^{(m)} _{U }} \R \underline{\Gamma} ^{(m)} _{Y} (\E)}
    \ar[d] ^-\sim
    \\
    {\hat{v} _+ \hat{v} ^! (\widehat{\D} ^{(m)} _{\U }   \smash{\widehat{\otimes}} ^\L _{ \D ^{(m)} _{U }}\E)}
    \ar[r] ^-\sim
    &
    {\R \underline{\Gamma} ^{(m)} _{Y_0} (\widehat{\D} ^{(m)} _{\U }
    \smash{\widehat{\otimes}} ^\L _{ \D ^{(m)} _{U }}\E),}
}
\end{equation}
où l'isomorphisme du haut (resp. du bas) résulte \ref{g+g!gamma} (resp. est de \cite[4.4.5]{Beintro2})
et celui de gauche découle de \ref{prop-comm-^} et \ref{prop-comm-^inv},
 est commutatif.
\end{lemm}
\begin{proof}
  Par \cite[1.5.3]{Be1}, on vérifie que le foncteur $\R \underline{\Gamma} ^{(m)} _{Y}$
  commute aux changements de base. Il en dérive que
  $\R \underline{\Gamma} ^{(m)} _{Y} (\E)\in D ^\mathrm{b} _\mathrm{qc}(\D _{U} ^{(m)})$.
  Le reste de la preuve découle de la commutation aux changements de base de
   l'image directe, de l'image inverse extraordinaire (prouvées précédemment)
   et du foncteur cohomologique local à support dans un sous-$\V$-schéma formel fermé et lisse.
\end{proof}

\begin{lemm}
  \label{tracecomphat}
  Soient $v$ : $Y \hookrightarrow U$ une immersion fermée de $\V$-schémas formels faibles lisses et
  $\E ^{(m _0)} \in D ^\mathrm{b} _\mathrm{coh}(\D _{U} ^{(m_0)})$.
  Pour tout entier $m \geq m_0 $, on note $\E ^{(m )} : = \D _{U} ^{(m )}\otimes _{\D _{U} ^{(m _0)}} \E ^{(m _0)}$.
  Lorsqu'il existe un entier
  $m \geq m _0$ tel que
  $\R \underline{\Gamma} ^{(m _0)} _{Y} (\E ^{(m _0)})\in D ^\mathrm{b} _\mathrm{coh}(\D _{U} ^{(m_0)})$ et
  $\R \underline{\Gamma} ^{(m)} _{Y} (\E^{(m )}) \in D ^\mathrm{b} _\mathrm{coh}(\D _{U} ^{(m)})$,
  $\widehat{\D} ^{(m)} _{\U ,\Q}   \otimes _{\widehat{\D} ^{(m _0)} _{\U ,\Q}}
  \R \underline{\Gamma} ^{(m _0)} _{Y _0 } ( \widehat{\E } ^{(m _0)} _\Q)
  \riso
  \R \underline{\Gamma} ^{(m)} _{Y _0 } ( \widehat{\E } ^{(m)} _\Q)$ canoniquement et
  bifonctoriellement en $Y $ et $\E$,
  et le diagramme canonique qui s'en déduit
\begin{gather}
\notag
  \xymatrix @R=0,4cm {
  {\widehat{\D} ^{(m)} _{\U ,\Q}   \otimes _{\widehat{\D} ^{(m _0)} _{\U ,\Q}}
  v _{+^{(m_0)}} v ^{!^{(m)}}  ( \widehat{\E } ^{(m _0)} _\Q)}
  \ar[r]^-\sim \ar[d] ^-\sim
  &
  {\widehat{\D} ^{(m)} _{\U ,\Q}   \otimes _{\widehat{\D} ^{(m _0)} _{\U ,\Q}}
  \R \underline{\Gamma} ^{(m _0)} _{Y _0} ( \widehat{\E } ^{(m _0)} _\Q)}
  \ar[d] ^-\sim
  \\
  {v _{+^{(m)}} v ^{!^{(m)}}  ( \widehat{\E } ^{(m)} _\Q)}
  \ar[r] ^-\sim
  &
  {\R \underline{\Gamma} ^{(m)} _{Y _0 } ( \widehat{\E } ^{(m)} _\Q)}
  }
  \end{gather}
est commutatif.
Lorsque, pour tout $m\geq m _0$,
$\R \underline{\Gamma} ^{(m)} _{Y} (\E^{(m )}) \in D ^\mathrm{b} _\mathrm{coh}(\D _{U} ^{(m)})$,
on dispose de résultats analogues
en remplaçant respectivement {\og $\widehat{\E } ^{(m)} _\Q$ \fg} par
{\og $\D ^\dag _{\U ,\Q}   \otimes _{\widehat{\D} ^{(m _0)} _{\U ,\Q}}
\widehat{\E} ^{(m _0)} _\Q $ \fg} et
{\og $(m)$ \fg} par {\og $\dag$ \fg}.
\end{lemm}
\begin{proof}
  Cela découle de \ref{lemmtracecomphat} de manière analogue au fait que \ref{prop-comm-^inv1}
  résulte de \ref{prop-comm-^inv}.
\hfill \hfill \qed \end{proof}

\begin{vide}\label{tracecomphat2}
  Soient $v$ : $Y \hookrightarrow U$ une immersion fermée de $\V$-schémas formels faibles lisses et
  $\FF ^{(m _0)} \in D ^\mathrm{b} _\mathrm{coh}(\D _{Y} ^{(m_0)})$.
  Le diagramme
  $$\xymatrix @R=0,4cm {
  {\D ^\dag _{\Y,\Q }   \otimes _{ \D ^{(m)} _{Y }} \FF ^{(m _0)}}
  \ar[r] ^-\sim
  \ar@{=}[d]
  &
  {\D ^\dag _{\Y ,\Q }   \otimes _{ \D ^{(m)} _{Y }} v ^{! ^{(m)}} v _{+^{(m)}} (\FF ^{(m _0)})}
  \ar[d] ^-\sim
  \\
  {\D ^\dag _{\Y ,\Q }   \otimes _{ \D ^{(m)} _{Y }} \FF ^{(m _0)}}
  \ar[r] ^-\sim
  &
  {v ^{! ^{\dag}} v _{+^{\dag}} (\D ^\dag _{\Y ,\Q }   \otimes _{ \D ^{(m)} _{Y }}  \FF ^{(m _0)}),}
}
$$
où l'isomorphisme de droite découle de \ref{prop-comm-^inv2} et \ref{prop-comm-^dir},
est commutatif. En effet, de manière analogue aux preuves de \ref{prop-comm-^inv2} ou \ref{prop-comm-^dir},
cela se vérifie
en complétant, tensorisant par $\Q$ puis passant à la limite l'isomorphisme canonique
$\FF ^{(m _0)} \riso v ^{! ^{(m)}} v _{+^{(m)}} (\FF ^{(m _0)})$.

Avec les notations \ref{tracecomphat}, on dispose du morphisme composé, noté $\mathrm{adj}$,
$ v _{+^{(m)}} v ^{! ^{(m)}} (\E ^{(m _0)}) \riso
\R \underline{\Gamma} ^{(m)} _{Y} (\E ^{(m _0)}) \rightarrow \E ^{(m _0)}$.
Supposons que, pour tout $m\geq m _0$,
$\R \underline{\Gamma} ^{(m)} _{Y} (\E^{(m )}) \in D ^\mathrm{b} _\mathrm{coh}(\D _{U} ^{(m)})$.
Il résulte alors de \ref{lemmtracecomphat} et \ref{tracecomphat} que le diagramme
$$\xymatrix @R=0,4cm {
  {\D ^\dag _{\U,\Q }   \otimes _{ \D ^{(m)} _{U }}  v _{+^{(m)}} v ^{! ^{(m)}} (\E ^{(m _0)})}
  \ar[r]  _-{\mathrm{adj}}
  \ar[d] ^-\sim
  &
  {\D ^\dag _{\U,\Q }   \otimes _{ \D ^{(m)} _{U }} \E ^{(m _0)}}
  \ar@{=}[d]
  \\
  {v _{+^{\dag}} v ^{! ^{\dag}} (\D ^\dag _{\U,\Q }   \otimes _{ \D ^{(m)} _{U }}  \E ^{(m _0)})}
  \ar[r]  _-{\mathrm{adj}}
  &
  {\D ^\dag _{\U,\Q }   \otimes _{ \D ^{(m)} _{U }} \E ^{(m _0)},}
}
$$
où le morphisme du bas est le morphisme d'adjonction (\cite{Vir04}) et celui
de gauche dérive de \ref{prop-comm-^inv2} et \ref{prop-comm-^dir},
est commutatif.
\end{vide}

\begin{prop}\label{globpresenfiniu+}
  Soit $v$ : $ Y \hookrightarrow U$ une immersion fermée de $\V$-schémas formels faibles affines et lisses.
  Pour tout $\D _Y ^{(m)}$-module à droite globalement de présentation finie $\M$, $v _+ (\M)$ est globalement
  de présentation finie et on a un isomorphisme canonique :
  $$ \Gamma ( U, v _+ (\M)) \riso \Gamma (Y ,\M) \otimes _{\Gamma (Y ,\D _Y ^{(m)}) }
  \Gamma (Y ,\D _{Y \hookrightarrow U} ^{(m)}).$$
  De même, en remplaçant {\og module à droite \fg} par {\og module à gauche \fg}.
\end{prop}
\begin{proof}
Comme $\M$ est un $\D _Y ^{(m)}$-module à droite globalement de présentation finie,
$\M \otimes _{\D ^{(m)} _Y} \D _{Y \hookrightarrow U} ^{(m)} \riso
\Gamma (Y, \M) \otimes _{\Gamma (Y, \D ^{(m)} _Y)} \D _{Y \hookrightarrow U} ^{(m)}$ (\ref{theoA}).
Pour tout entier $n$, il découle de la $\O _U$-cohérence de $\D ^{(m)} _{U,n}$ et du théorème de type $A$
pour les $\O _U$-modules cohérents que le morphisme canonique
$\Gamma (Y, \O _Y) \otimes _{\Gamma (U,\O _U) } \Gamma (U, \D ^{(m)} _{U,n}) \rightarrow
\Gamma ( Y, v ^* \D ^{(m)} _{U,n})$ est un isomorphisme.
Le produit tensoriel et le foncteur $\Gamma (U, -)$ commutant aux limites inductives filtrantes,
le morphisme
$\Gamma (Y, \O _Y) \otimes _{\Gamma (U,\O _U) } \Gamma (U, \D ^{(m)} _{U}) \rightarrow
\Gamma ( Y, \D _{Y \hookrightarrow U} ^{(m)})$ est un isomorphisme.
Soient $U'$ un ouvert principal de $U$ et $Y':= Y \cap U'$.
On obtient de même
$\Gamma (Y ', \O _Y) \otimes _{\Gamma (U',\O _U) } \Gamma (U ', \D ^{(m)} _{U}) \riso
\Gamma ( Y ', \D _{Y \hookrightarrow U} ^{(m)})$.
Or,
$\Gamma (Y ', \O _Y)\riso \Gamma (Y, \O _Y) \otimes _{\Gamma (U,\O _U) }\Gamma (U',\O _U)$
(on utilise la remarque de la preuve de \ref{stabimmfer}).
Finalement, on obtient :
$$\Gamma ( Y, \D _{Y \hookrightarrow U} ^{(m)})
\otimes _{\Gamma (U, \D ^{(m)} _{U})} \Gamma (U ', \D ^{(m)} _{U})
\riso
\Gamma ( Y ', \D _{Y \hookrightarrow U} ^{(m)}).$$
Comme le foncteur {\og faisceautisation \fg}  commute à $v _*$ (en effet, $v$ est une
immersion fermée), $v _* (\M \otimes _{\D ^{(m)} _Y} \D _{Y \hookrightarrow U} ^{(m)})$
est le faisceau associé au préfaisceau défini sur les ouverts principaux $U'$ de $U$ par :\newline
$U' \mapsto \Gamma (Y, \M) \otimes _{\Gamma (Y, \D ^{(m)} _Y)}\Gamma ( Y, \D _{Y \hookrightarrow U} ^{(m)})
\otimes _{\Gamma (U, \D ^{(m)} _{U})} \Gamma (U ', \D ^{(m)} _{U}).$
Il reste ainsi à vérifier que
$\Gamma (Y, \M) \otimes _{\Gamma (Y, \D ^{(m)} _Y)}\Gamma ( Y, \D _{Y \hookrightarrow U} ^{(m)})$
est un $\Gamma (U, \D ^{(m)} _{U})$-module à droite de type fini. Le produit tensoriel étant exact à droite
et $\Gamma (Y, \M)$ étant un $\Gamma (Y, \D ^{(m)} _Y)$-module à droite de type fini, il suffit de prouver
que $\Gamma ( Y, \D _{Y \hookrightarrow U} ^{(m)})$ est un $\Gamma (U, \D ^{(m)} _{U})$-module à droite de type fini.
On conclut via
$\Gamma (U, \D ^{(m)} _{U})
\twoheadrightarrow \Gamma (Y, \O _Y) \otimes _{\Gamma (U,\O _U) } \Gamma (U, \D ^{(m)} _{U}) \riso
\Gamma ( Y, \D _{Y \hookrightarrow U} ^{(m)})$.

\hfill \hfill \qed \end{proof}

\begin{coro}
  \label{coro-globpresenfiniu+}
  Soit $v$ : $ Y \hookrightarrow U$ une immersion fermée de $\V$-schémas formels faibles lisses.
  Pour tout $\D _Y ^{(m)}$-module cohérent $\M$, $v _+ (\M)$ est $\D _U ^{(m)}$-cohérent.
\end{coro}

\begin{nota}\label{notafidplatpq}
Soient $f$ : $P ' \rightarrow P$ un morphisme de $\V$-schémas formels faibles lisses,
  $T _0$ (resp. $T ' _0$) un diviseur de $P_0$ (resp. $P ' _0$), $U$ (resp. $U'$)
  l'ouvert de $P$ (resp. $P'$) complémentaire de $T _0$ (resp. $T '_0$),
  $j$ : $U \hookrightarrow P$ et $j'$ : $U' \hookrightarrow P'$ les immersions ouvertes.
  On suppose que $f $ se factorise par $g$ : $U' \rightarrow U$.
  On reprend les notations de Huyghe \cite{huyghe-comparaison} concernant le symbole {\og * \fg},
  apparaissant par exemple avec
les faisceaux $\O _{\PP } (\hdag * T _0)$ et $\D ^\dag _{\PP } (\hdag * T _0)$.
De plus, pour tout entier $m$, on notera
$\D ^{(m)} _{\PP } (\hdag * T _0):=
\O _{\PP } (\hdag * T _0) \otimes _{\O _{\PP } } \D ^{(m)} _{\PP }$
et
$\D _{\PP} (\hdag T _0) _\Q:=
\O _{\PP } (\hdag T _0) _\Q \otimes _{\O _{\PP } } \D ^{(m)} _{\PP }$, ce dernier ne dépendant pas de
$m$ à isomorphisme canonique près. Enfin,
\begin{gather}\notag
  \D ^{(m)} _{\PP '\rightarrow \PP } (\hdag *T '_0  ):=
\O _{\PP '} (\hdag * T '_0 )\otimes  _{f ^{-1}\O _{\PP}(\hdag * T _0 )}
f ^{-1}\D ^{(m)} _{\PP } (\hdag *T _0 ) \\
\tilde{\leftarrow}
\O _{\PP '} (\hdag * T '_0 )\otimes _{\O _{\PP '} }
(\O _{\PP '} \otimes  _{f ^{-1}\O _{\PP}} f ^{-1}\D ^{(m)} _{\PP })
= \O _{\PP '} (\hdag * T '_0 )\otimes _{\O _{\PP '} }
\D ^{(m)} _{\PP '\rightarrow \PP } . \notag
\end{gather}
De même, on pose $\D ^{(m)} _{\PP \leftarrow \PP' } (\hdag *T '_0 ) :=
\O _{\PP '} (\hdag * T '_0 )\otimes  _{\O _{\PP}}\D ^{(m)} _{\PP \leftarrow \PP' } $.

\end{nota}

\begin{prop}\label{fidplatpq}
  Avec les notations de \ref{notafidplatpq}, les homomorphismes ca-\linebreak noniques
 $ \O _{\PP } (\hdag * T _0) \rightarrow j _* \O _{\U}$
et 
$\D ^\dag _{\PP } (\hdag * T _0) \rightarrow j _* \D ^\dag _{\U}$
  sont fidèlement plats à droite et à gauche.
\end{prop}
\begin{proof}
On sait déjà que ces morphismes sans le symbole "$*$" sont plats (cela est prouvé au cours de la preuve
de \cite[4.3.10]{Be1}). Il en découle que les morphismes de \ref{fidplatpq} sont plats.
Par analogie, contentons-nous d'établir la fidèle platitude du deuxième morphisme de \ref{fidplatpq}
dont la preuve est analogue à celle \cite[4.3.10]{Be1} : grâce à \cite[3.3.5 et 4.3.8]{Be1},
  il suffit de prouver que si $\PP '$ est un ouvert affine de $\PP$ sur lequel il existe $f \in \Gamma (\PP ', \O _{\PP})$
  relevant une équation locale de $T _0$, et si $M$ est un
  $\Gamma (\PP ',\D ^\dag _{\PP } (\hdag * T _0))$-module monogène de présentation finie
  tel que $ \Gamma( \PP' \cap \U , \D ^\dag _{\U}) \otimes _{\Gamma (\PP ',\D ^\dag _{\PP } (\hdag * T _0))} M =0$,
  alors $M =0$. Soient $m'\geq m$ assez grands tels qu'il existe un
  $\Gamma (\PP ',\widehat{\B} ^{(m')} _{\PP} (T _0) \widehat{\otimes}\widehat{\D} ^{(m)} _{\PP} (* T _0))$-module
  monogène de présentation finie $M'$
  induisant $M$ par extension des scalaires et vérifiant
  $\Gamma ( \PP' \cap \U , \widehat{\D} ^{(m)} _{\PP})
  \otimes _{ \Gamma (\PP ',\widehat{\B} ^{(m')} _{\PP} (T _0) \widehat{\otimes}\widehat{\D} ^{(m)} _{\PP} (* T _0))}
  M' =0$. Sur $\PP'$, le foncteur $(* T _0)$ est isomorphe au foncteur localisation en $f$.
  Soit $e$ un générateur de $M'$. On note $E$
  le sous-$\Gamma (\PP ',\widehat{\B} ^{(m')} _{\PP} (T _0) \widehat{\otimes}\widehat{\D} ^{(m)} _{\PP})$-module
  de $M'$ engendré par $e$.
  Comme $\Gamma (\PP ',\widehat{\B} ^{(m')} _{\PP} (T _0) \widehat{\otimes}\widehat{\D} ^{(m)} _{\PP})$
  est noethérien, $E$ est un
  $\Gamma (\PP ',\widehat{\B} ^{(m')} _{\PP} (T _0) \widehat{\otimes}\widehat{\D} ^{(m)} _{\PP})$-module monogène,
   de présentation finie, sans $f$-torsion et tel que $E [ 1 /f] \riso M'$.
   Comme $M' / \pi  M' =0$, on obtient l'égalité $(E / \pi E )_{\overline{f}}=0$, où $\overline{f}$ est l'image canonique
   de $f$ sur $A /\pi A$. On termine la preuve de façon similaire en reprenant la preuve de \cite[4.3.10]{Be1}
   à partir de {\og Si $E _0 =E /pE$, on voit ainsi que $(E _0) _f=0$. \fg}.

\hfill \hfill \qed \end{proof}

\begin{coro}\label{corofidplatpq}
  Avec les notations de \ref{notafidplatpq}, pour tout entier $m$,
  les homomorphismes ,
  $j _* \O _U \rightarrow \O _{\PP } (\hdag * T _0)$,
  $j _* \D ^{(m)} _U \rightarrow \D ^{(m)} _{\PP } (\hdag * T _0)$ et
$j _*  \D _{U ,\Q} \rightarrow  \D _{\PP} (\hdag T _0) _\Q$
  sont fidèlement plats à droite et à gauche.
De plus, $j _* \D  _U \rightarrow \D ^\dag _{\PP } (\hdag * T _0)$
et $j _* \D ^{(m)} _U \rightarrow \D ^\dag _{\PP } (\hdag T _0) _\Q$
  sont plats à droite et à gauche.
\end{coro}
\begin{proof}
Traitons d'abord la pleine fidélité.
Comme $T _0$ est un diviseur,
il existe une base de voisinages affines de $P$ telle que celle induite sur $U$ soit constituée
d'ouverts affines.
Or, pour tout ouvert affine $V $ de $U $,
l'extension $\Gamma (V , \O _{U }) \rightarrow
\Gamma (V  , \O _{\U })$ est fidèlement plate à droite et à gauche.
L'extension $j _* \O _U \rightarrow j _* \O _\U$ est donc fidèlement plate à droite et à gauche.
Grâce à \ref{fidplatpq}, on l'obtient pour $j _* \O _U \rightarrow \O _{\PP } (\hdag * T _0)$.

Celle du deuxième morphisme du corollaire découle de l'isomorphisme canonique :
$\O _{\PP } (\hdag * T _0) \otimes _{j _* \O _U }  j _* \D ^{(m)} _U
\riso \D ^{(m)} _{\PP } (\hdag * T _0)$.
Le troisième s'en déduit en lui appliquant $\otimes \Q$.

Passons à la platitude. Pour tout $m$, les extensions
$j _* \D ^{(m)} _U \rightarrow  j _* \widehat{\D} ^{(m)} _\U $ sont plates.
Par passage à la limite sur le niveau,
$j _* \D _U \rightarrow  j _* \D ^\dag _\U $ est plate.
La proposition \ref{fidplatpq} nous permet de conclure le premier cas.
La platitude du dernier morphisme résulte, en plus des arguments précédents,
de la platitude de
$j _* \widehat{\D} ^{(m)} _\U \rightarrow j _* \D ^\dag  _{\U,\Q}$
(on utilise \cite[3.5.4]{Be1}).
\hfill \hfill \qed \end{proof}

\begin{rema}\label{remathABdag*}
  Soient $\PP$ un $\V$-schéma formel lisse, $T _0$ un diviseur de $P_0$ et $\E$ un
  $\D ^\dag _{\PP} (\hdag T _0)$-module. On suppose qu'il existe
  une section $f \in \Gamma (\PP, \O _{\PP})$ relevant une équation locale de $T _0$. Comme $\PP$ est noethérien,
  le préfaisceau qui a un ouvert $\PP '$ associe $\Gamma (\PP ', \E) [1 /f]$ est un faisceau.
  Pour tout entier $q \geq 0$, il en dérive $H ^q (\PP, \E) [1 /f] \riso H ^q (\PP, \E (*T _0))$.

  De plus, si $\FF $ est un $\D ^\dag _{\PP} (\hdag * T _0)$-module globalement de présentation finie, alors il
  existe un $\D ^\dag _{\PP} (\hdag T _0)$-module de présentation finie $\E$ tel que
  $\E (* T _0) \riso \FF$.

  Supposons à présent $\PP$ affine. On obtient alors les théorèmes $A$ et $B$ :

  A) Pour qu'un $\D ^\dag _{\PP} (\hdag * T _0)$-module $\FF$ soit globalement de présentation finie,
  il faut et il suffit que $\Gamma (\PP, \FF)$ soit un $\Gamma ( \PP,\D ^\dag _{\PP} (\hdag * T _0))$-module
  de présentation finie, et que l'homomorphisme
  $\D ^\dag _{\PP} (\hdag * T _0)\otimes _{\Gamma (\PP, \D ^\dag _{\PP} (\hdag * T _0))} \Gamma (\PP, \FF)$
  soit un isomorphisme.

  B) Pour tout $\D ^\dag _{\PP} (\hdag * T _0)$-module $\FF$ globalement de présentation finie,
  pour tout entier $q \geq 1$, $H ^q (\PP, \FF)=0$.

  En outre, les foncteurs $\Gamma (\PP, -)$ et
  $\D ^\dag _{\PP} (\hdag * T _0)\otimes _{\Gamma (\PP, \D ^\dag _{\PP} (\hdag * T _0))}- $ sont
  des équivalences quasi-inverses entre la catégorie des
  $\Gamma ( \PP,\D ^\dag _{\PP} (\hdag * T _0))$-modules
  de présentation finie et celle des
  $\D ^\dag _{\PP} (\hdag * T _0)$-modules de présentation finie.

  En effet, en vertu de \cite[3.6]{Be1}, on le sait sans le symbole "$*$". Les remarques précédentes
  nous permettent alors de conclure.
\end{rema}

\begin{prop}\label{4.3.12.Be1}
  Avec les notations de \ref{fidplatpq}, pour tout $\D ^\dag _{\PP } (\hdag * T _0)$-module localement
  de présentation finie $\E$, l'homomorphisme canonique
  \begin{equation}
    \label{4.3.12.Be1-iso}
  j _* \D ^\dag _{\U} \otimes _{\D ^\dag _{\PP } (\hdag * T _0)} \E \rightarrow j _* j^* \E
  \end{equation}
  est un
  isomorphisme. De plus, pour qu'un morphisme de $\D ^\dag _{\PP } (\hdag * T _0)$-modules localement
  de présentation finie $\E$ soit injectif, surjectif ou bijectif,
  il faut et il suffit que sa restriction à $\U$ le soit.

  On bénéficie d'un résultat analogue en remplaçant {\og $\D ^\dag _{\PP } (\hdag * T _0)$ \fg}
  par {\og $\D ^{(m)} _{\PP } (\hdag * T _0)$ \fg} et {\og $\D ^\dag _{\U}$ \fg}
  par {\og $\D ^{(m)} _{\U}$ \fg}.
\end{prop}
\begin{proof}
  On reprend les arguments de \cite[4.3.12]{Be1} :
  l'exactitude à droite des foncteurs en $\E$ (grâce au théorème $B$ pour les $\D ^\dag _{\U}$-modules
  cohérents ou $\D ^{(m)} _{\U}$-modules cohérents impliquent l'exactitude de $j _*$ dans nos cas)
  nous ramène, afin d'établir l'isomorphisme \ref{4.3.12.Be1-iso},
  au cas immédiat où $\E=\D ^\dag _{\PP } (\hdag * T _0)$ ou
  $\E=\D ^{(m)} _{\PP } (\hdag * T _0)$. La proposition \ref{fidplatpq} permet ensuite d'obtenir le critère
  d'injectivité, surjectivité ou bijectivité.
\hfill \hfill \qed \end{proof}

\begin{lemm}\label{j*D->D}
Il existe un morphisme canonique
  $j ' _* \D ^{(m)} _{U '\rightarrow U } \rightarrow \D ^{(m)} _{\PP '\rightarrow \PP } (\hdag *T '_0 )$
  de $(j ' _*\D ^{(m)} _{U'},j ' _* g ^{-1} \D ^{(m)} _{U} )$-bimodules.
  En outre, le morphisme induit par extension des scalaires,
  $\D ^{(m)} _{\PP '} (\hdag *T '_0 ) \otimes _{j ' _* \D ^{(m)} _{U '}} j ' _* \D ^{(m)} _{U '\rightarrow U }
  \rightarrow \D ^{(m)} _{\PP '\rightarrow \PP } (\hdag *T '_0 )$,
est un isomorphisme.
\end{lemm}
\begin{proof}
  L'injection canonique $g ^* \D ^{(m)} _U \hookrightarrow \hat{g} ^* \widehat{\D} ^{(m)} _{\U}$ induit la suivante
$j ' _* g ^* \D ^{(m)} _U \hookrightarrow j ' _* \hat{g} ^* \widehat{\D} ^{(m)} _{\U}$.
Puisque le morphisme $\O _{\PP '} (\hdag * T '_0 ) \rightarrow j ' _* \O _{\U '}$
est injectif, par un calcul en coordonnées locales, on obtient les injections (adjonction de $j '_*$ puis complétion)
$$\D ^{(m)} _{\PP '\rightarrow \PP } (\hdag *T _0 )
\hookrightarrow
j '_* \D ^{(m)} _{\U ' \rightarrow \U }
\hookrightarrow
j '_* \widehat{\D} ^{(m)} _{\U ' \rightarrow \U }.
$$
Via un calcul en coordonnées locales et l'injection canonique
$j '_* \O _{U '} \hookrightarrow \O _{\PP'}(\hdag *T '_0 )$,
on établit la factorisation :
$$\xymatrix @R=0,4cm {
{j ' _* \hat{g} ^* \widehat{\D} ^{(m)} _{\U}}
&
{}
\\
{j ' _* g ^* \D ^{(m)} _U}
\ar@{.>}[r] \ar@{^{(}->}[u]
&
{\D ^{(m)} _{\PP '\rightarrow \PP } (\hdag *T '_0 ).} \ar@{^{(}->}[ul]
}$$
Par extension des scalaires, il dérive de cette factorisation
$\O  _{\PP '} (\hdag *T '_0 ) \otimes _{j ' _* \O  _{U '}} j ' _* \D ^{(m)} _{U '\rightarrow U }
  \rightarrow  \D ^{(m)} _{\PP '\rightarrow \PP } (\hdag *T '_0 )$.
Par un calcul en coordonnées locales, on vérifie
que celui-ci est un isomorphisme.
  On conclut grâce à l'isomorphisme canonique
  $\O  _{\PP '} (\hdag *T '_0 ) \otimes _{j ' _* \O  _{U '}} j ' _* \D ^{(m)} _{U '\rightarrow U }
  \riso \D ^{(m)} _{\PP '} (\hdag *T '_0 ) \otimes _{j ' _* \D ^{(m)} _{U '}} j ' _* \D ^{(m)} _{U '\rightarrow U }$.
\hfill \hfill \qed \end{proof}

\begin{prop}\label{compiminv}
Avec les notations \ref{notafidplatpq}, on suppose $f$ lisse.
Pour tout $\D ^{(m)} _U$-module localement en $P$ de présentation
  finie $\E$, on dispose d'un isomorphisme canonique
  \begin{gather}\notag
{\D ^\dag _{\PP '} (\hdag T  '_0 )_\Q \otimes _{j ' _* \D ^{(m)} _{U'}}
  j '_*(\D ^{(m)} _{U' \rightarrow U} \otimes _{ g ^{-1} \D ^{(m)} _U} g ^{-1} \E )}
  \\
\riso
  {\D ^\dag _{\PP '\rightarrow \PP } (\hdag T ' _0 ) _\Q  \otimes _{f ^{-1} \D ^\dag _{\PP} (\hdag T _0) _\Q}
  f ^{-1} ( \D ^\dag _{\PP} (\hdag  T _0) _\Q \otimes _{j _* \D ^{(m)} _U} j _* \E ).} \notag
\end{gather}
  En outre, ceux-ci sont transitifs en $f$ et $g$.
\end{prop}
\begin{proof}
  L'homomorphisme
  $j '_*\D ^{(m)} _{U' \rightarrow U} \otimes _{ j '_* g ^{-1} \D ^{(m)} _U} j '_* g ^{-1} \E
  \rightarrow \linebreak
  j '_*(\D ^{(m)} _{U' \rightarrow U} \otimes _{ g ^{-1} \D ^{(m)} _U} g ^{-1} \E )$
  induit par extension des scalaires
\begin{gather}
  \notag {\D ^{(m)} _{\PP '} (\hdag *T '_0 ) \otimes _{j ' _* \D ^{(m)} _{U '}}
j '_*\D ^{(m)} _{U' \rightarrow U} \otimes _{ j '_* g ^{-1} \D ^{(m)} _U} j '_* g ^{-1} \E}
 \\
\rightarrow
{\D ^{(m)} _{\PP '} (\hdag *T '_0 ) \otimes _{j ' _* \D ^{(m)} _{U '}}
j '_*(\D ^{(m)} _{U' \rightarrow U} \otimes _{ g ^{-1} \D ^{(m)} _U} g ^{-1} \E ).}
\notag
\end{gather}
Par \ref{j*D->D}
  et via le morphisme canonique $f ^{-1} j _* \rightarrow  j ' _* g ^{-1}$, construit par adjonction,
   on obtient alors le suivant
  $\D ^{(m)} _{\PP '\rightarrow \PP } (\hdag *T '_0 )
  \otimes _{ f ^{-1} j _* \D ^{(m)} _U} f ^{-1} j _* \E
\rightarrow
\D ^{(m)} _{\PP '} (\hdag *T '_0 ) \otimes _{j ' _* \D ^{(m)} _{U '}}
j '_*(\D ^{(m)} _{U' \rightarrow U} \otimes _{ g ^{-1} \D ^{(m)} _U} g ^{-1} \E ).$
  D'où
  \begin{gather}\notag
    \D ^{(m)} _{\PP '\rightarrow \PP } (\hdag *T '_0 ) \otimes _{f ^{-1} \D ^{(m)} _{\PP } (\hdag *T _0 )}
  f ^{-1} (\D ^{(m)} _{\PP } (\hdag *T _0 )
  \otimes _{  j _* \D ^{(m)} _U}  j _* \E)
  \rightarrow \\
  \D ^{(m)} _{\PP '} (\hdag *T '_0 ) \otimes _{j ' _* \D ^{(m)} _{U '}}
  j '_*(\D ^{(m)} _{U' \rightarrow U} \otimes _{ g ^{-1} \D ^{(m)} _U} g ^{-1} \E ).\label{compiminv-equ3}
  \end{gather}
Comme $\E$ est localement en $P$ de présentation finie, il découle de \ref{globpfinvlis} que
$\D ^{(m)} _{U' \rightarrow U} \otimes _{ g ^{-1} \D ^{(m)} _U} g ^{-1} \E$ est localement en $P'$ de
présentation finie.
Par \ref{j*->Rj*iso}, le terme du bas de \ref{compiminv-equ3} est
localement de présentation finie sur $\D ^{(m)} _{\PP '} (\hdag *T '_0 )$.
De plus, il en est de même de celui du haut.
Comme \ref{compiminv-equ3} est un isomorphisme au-dessus de $\U'$, il découle de \ref{4.3.12.Be1}
que \ref{compiminv-equ3} est un isomorphisme.
Pour en conclure la construction de l'isomorphisme de \ref{compiminv}, il suffit
de vérifier que l'homomorphisme canonique
$\D ^\dag _{\PP ' } (\hdag T '_0 ) _\Q \otimes _{\D ^{(m)} _{\PP ' } (\hdag *T '_0 )}
\D ^{(m)} _{\PP '\rightarrow \PP } (\hdag *T '_0 ) \rightarrow
\D ^\dag _{\PP '\rightarrow \PP } (\hdag T '_0 ) _\Q$, induit par extension des scalaires, est un isomorphisme.
Comme $f$ est lisse, ce dernier est un morphisme de
$\D ^\dag _{\PP ' } (\hdag T '_0 ) _\Q$-modules cohérents. Par \cite[4.3.12]{Be1},
il suffit d'établir que sa restriction à $\U '$ est un isomorphisme.
Or, comme $g$ est lisse, $\D ^{(m)} _{\U ' \rightarrow \U}$ est $\D ^{(m)} _{\U'}$-cohérent.
Le morphisme $\widehat{\D} ^{(m)} _{\U '} \otimes _{\D ^{(m)} _{\U'}} \D ^{(m)} _{\U ' \rightarrow \U}
\rightarrow \widehat{\D} ^{(m)} _{\U ' \rightarrow \U}$ est donc un isomorphisme.
D'où $\widehat{\D} ^{(m)} _{\U ',\Q} \otimes _{\D ^{(m)} _{\U'}} \D ^{(m)} _{\U ' \rightarrow \U}
\riso \widehat{\D} ^{(m)} _{\U ' \rightarrow \U,\Q}$. On termine grâce à la commutation à isogénie près de l'image
inverse extraordinaire par un morphisme lisse au niveau (\cite[4.3.3]{Beintro2}).

En ce qui concerne la commutation aux isomorphismes de composition
  des images inverses extraordinaires, il s'agit de vérifier la commutativité d'un diagramme canonique.
Par \cite[4.3.12]{Be1},
il suffit de l'établir pour sa restriction à $\U ''$, ce qui est aisé.
\hfill \hfill \qed \end{proof}

\begin{prop}\label{g-+commdag*f+}
Avec les notations \ref{notafidplatpq}, on suppose $T ' _0=f ^{-1} (T_0)$.
  Soit $\E '$ un $\D ^{(0)} _{U'}$-module localement en $P '$ de présentation finie
  tel que $g _+ (\E')$ soit un $\D ^{(0)} _U$-module localement en $P$ de présentation finie.
  On dispose alors d'un morphisme canonique :
  \begin{gather}
\notag
  {\D ^\dag _{\PP } (\hdag * T _0) \otimes _{j _* \D ^{(0)} _U} j _* g _+ (\E ')}
  \\ \rightarrow
  {\R f _* (
  \D ^\dag _{\PP \leftarrow \PP '} (\hdag *T _0 )
  \otimes _{\D ^\dag _{\PP '} (\hdag * T '_0 )} ^\L
  ( \D ^\dag _{\PP '} (\hdag{* T '_0 }) \otimes ^\L  _{j '_* \D ^{(0)} _{U'}}   j '_* \E ' ))}.\label{gathg-+commdag*f+}
  \end{gather}
  Celui-ci est un isomorphisme au dessus de $\U$.
De plus, ceux-ci sont transitifs en $f$ et $g$.
\end{prop}
\begin{proof}
Par \ref{globpresenfiniu+},
on remarque que lorsque $f$ est une immersion fermée,
$g _+ (\E')$ est forcément un $\D ^{(0)} _U$-module localement en $P$ de présentation finie.
De plus, le morphisme $f$ est le composé de son graphe $P '\hookrightarrow P ' \times P$ suivant
de la projection $P' \times P \rightarrow P$. On se ramène ainsi à supposer
que $f$ est soit lisse soit une immersion fermée.

Grâce à \ref{j*->Rj*iso},
$j _* \R g _* ( \D ^{(0)} _{U \leftarrow U '} \otimes ^\L _{\D ^{(0)} _{U '}} \E ')
\riso \R f _* \R j ' _* ( \D ^{(0)} _{U \leftarrow U '} \otimes ^\L _{\D ^{(0)} _{U '}} \E ')$.
Lorsque $f$ est lisse, de manière analogue à \cite[2.4.6.2]{Beintro2},
on dispose d'un quasi-isomorphisme canonique
$\Omega ^\bullet _{U ' /U} \otimes _{\O _{U'}} \D _{U '} ^{(0)} [d _{U'/U}]
\rightarrow \D ^{(0)} _{U \leftarrow U'}$.
En particulier,
$\D ^{(0)} _{U \leftarrow U '}$ admet une résolution finie par des
$\D ^{(0)} _{ U '}$-modules libres de type fini.
Lorsque $f$ est une immersion fermée, $\D ^{(0)} _{U \leftarrow U '}$ est
un $\D ^{(0)} _{ U '}$-module libre.
Dans les deux cas, il en découle que le morphisme
$j ' _*  \D ^{(0)} _{U \leftarrow U '} \otimes ^\L _{j ' _* \D ^{(0)} _{U '}} j ' _* \E '
\rightarrow
\R j ' _* ( \D ^{(0)} _{U \leftarrow U '} \otimes ^\L _{\D ^{(0)} _{U '}} \E ')$
est un isomorphisme.
On obtient par composition :
$j _* g _+ (\E ')
\riso
\R f _*( j ' _*  \D ^{(0)} _{U \leftarrow U '} \otimes ^\L _{j ' _* \D ^{(0)} _{U '}} j ' _* \E ')$.
Nous aurons besoin du lemme ci-après.
\begin{lemm}\label{j*D->Dtordu}
  Pour tout entier $m$, il existe un homomorphisme canonique \linebreak
  $j ' _* \D ^{(m)} _{U \leftarrow U '} \rightarrow \D ^{(m)} _{\PP \leftarrow \PP '} (\hdag *T _0 )$
  de $(j ' _* g ^{-1} \D ^{(m)} _{U}, j ' _*\D ^{(m)} _{U'})$-bimodules.
\end{lemm}
\begin{proof}
En rajoutant les faisceaux de formes différentielles de degré maximum,
cela se vérifie de manière analogue à \ref{j*D->D}.
\hfill \hfill \qed \end{proof}
On dispose d'un morphisme canonique
$\D ^{(m)} _{\PP \leftarrow \PP '} (\hdag *T _0 )\rightarrow
\D ^\dag _{\PP \leftarrow \PP '} (\hdag *T _0 )$. Avec \ref{j*D->Dtordu}, on obtient
$j ' _* \D ^{(m)} _{U \leftarrow U '} \rightarrow
\D ^\dag _{\PP \leftarrow \PP '} (\hdag *T _0 )$.
D'où : \linebreak
$\R f _*( j ' _*  \D ^{(0)} _{U \leftarrow U '} \otimes ^\L _{j ' _* \D ^{(0)} _{U '}} j ' _* \E ')
\rightarrow
\R f _*( \D ^\dag _{\PP \leftarrow \PP '} (\hdag *T _0 ) \otimes ^\L _{j ' _* \D ^{(0)} _{U '}} j ' _* \E ')
\linebreak
\riso \R f _* (
  \D ^\dag _{\PP \leftarrow \PP '} (\hdag *T _0 )
  \otimes _{\D ^\dag _{\PP '} (\hdag * T '_0 )} ^\L
  \D ^\dag _{\PP '} (\hdag{* T '_0 }) \otimes ^\L _{j '_* \D ^{(0)} _{U'}} j '_* \E ').$
On obtient par composition : \newline
$\theta $ : $ j _* g _+ (\E ') \rightarrow
\R f _* (
  \D ^\dag _{\PP \leftarrow \PP '} (\hdag *T _0 )
  \otimes _{\D ^\dag _{\PP '} (\hdag * T '_0 )} ^\L
  \D ^\dag _{\PP '} (\hdag{* T '_0 }) \otimes ^\L _{j '_* \D ^{(0)} _{U'}} j '_* \E ').$
  Il en résulte par extension la construction du morphisme de \ref{gathg-+commdag*f+}.

La restriction à $\U$ de $\theta$ correspond au morphisme composé :
$$\R g _* ( \D ^{(0)} _{U \leftarrow U '} \otimes ^\L _{\D ^{(0)} _{U '}} \E ')
\rightarrow
\R g _* ( \widehat{\D} ^{(0)} _{\U \leftarrow \U '} \otimes ^\L _{\D ^{(0)} _{U '}} \E ')
\rightarrow
\R g _* ( \D ^\dag _{\U \leftarrow \U '} \otimes ^\L _{\D ^{(0)} _{U '}} \E ').$$
D'après \ref{prop-comm-^},
$\widehat{\D} ^{(0)} _{\U }   \otimes _{ \D ^{(0)} _{U }}
\R g _* ( \D ^{(0)} _{U \leftarrow U '} \otimes ^\L _{\D ^{(0)} _{U '}} \E ')
\rightarrow
\R g _* ( \widehat{\D} ^{(0)} _{\U \leftarrow \U '} \otimes ^\L _{\D ^{(0)} _{U '}} \E ')$
est un isomorphisme. De plus, il résulte de \cite[3.5.3.(ii)]{Beintro2} (en passant à la limite),
que $ \D ^\dag _{\U}\otimes _{ \widehat{\D} ^{(0)} _{\U}}
\R g _* ( \widehat{\D} ^{(0)} _{\U \leftarrow \U '} \otimes ^\L _{\D ^{(0)} _{U '}} \E ')
\rightarrow
\R g _* ( \D ^\dag _{\U \leftarrow \U '} \otimes ^\L _{\D ^{(0)} _{U '}} \E ')$
est aussi un isomorphisme. Il en dérive que la restriction à $\U$ de \ref{gathg-+commdag*f+}
est un isomorphisme.
\hfill \hfill \qed \end{proof}

\subsection{Description des isocristaux surconvergents sur les schémas affines et lisses}

  Soit $Y ^\dag $ un $\V$-schéma formel faible affine et lisse. On remarque que grâce à Elkik (\cite{elkik}),
  il existe un
$\V$-schéma affine et lisse $Y$ dont le complété faible est isomorphe à $Y ^\dag$.

\begin{vide}
Il résulte de \cite[2.5.2]{Berig}, qu'il existe un foncteur pleinement fidèle de la catégorie
des isocristaux surconvergents sur $Y _0$ dans celle des
$\Gamma (Y ^\dag , \O _{Y ^\dag ,\Q})$-modules projectifs de type fini
munis d'une connexion intégrable.
Or, pour tout $\Gamma (Y ^\dag , \O _{Y ^\dag ,\Q})$-module de type fini $E$
muni d'une connexion intégrable, le morphisme canonique
$\O _{Y^\dag \Q} \otimes _{\Gamma (Y ^\dag ,\O _{Y^\dag \Q})}E \rightarrow
\D _{Y^\dag \Q} \otimes _{\Gamma (Y ^\dag ,\D _{Y^\dag \Q})}E $ est un isomorphisme.
La catégorie des
$\Gamma (Y ^\dag , \O _{Y ^\dag \Q})$-modules de type fini
munis d'une connexion intégrable est ainsi équivalente à celle des
$\D _{Y ^\dag , \Q}$-modules globalement de présentation finie et qui soient en outre
$\O _{Y ^\dag ,\Q}$-cohérents.

Il en résulte un foncteur pleinement fidèle,
noté $\sp _*$, de la catégorie
des isocristaux surconvergents sur $Y _0$ dans celle des
$\D _{Y ^\dag , \Q}$-modules globalement de présentation finie et
$\O _{Y ^\dag ,\Q}$-cohérents.
\end{vide}

\begin{prop}\label{propmodeisocaffi}
  Soient $E$ un isocristal surconvergent sur $Y_0$ et $\E := \sp _* (E)$.
  Il existe $\E ^{(m)}$, un $\D ^{(m)} _{Y ^\dag }$-module globalement de présentation et
  $\O _{Y ^\dag }$-cohérent, et un isomorphisme $\D _{Y ^\dag , \Q}$-linéaire $\E ^{(m)} _\Q \riso \E$.
\end{prop}
\begin{proof}
  Notons $A ^\dag  :=\Gamma (Y ^\dag , \O _{Y ^\dag})$ et $E:= \Gamma (Y^\dag , \E)$.
  Comme $E$ est un $A ^\dag _K $-module projectif de type fini,
  il existe un $A ^\dag _K $-module $F$, un entier $r$ et un isomorphisme $A ^\dag _K $-linéaire :
  $E \oplus F \riso (A ^\dag _K ) ^r$.
  En notant $\widehat{E} := \widehat{A} _K \otimes _{A ^\dag _K} E$ et
  $\widehat{F} := \widehat{A} _K \otimes _{A ^\dag _K} F$, on obtient l'isomorphisme
  $\widehat{E} \oplus \widehat{F} \riso(\widehat{A} _K)^r$ ainsi que les injections canoniques
  $E \hookrightarrow \widehat{E}$ et $F \hookrightarrow \widehat{F}$. Par abus de notations, on considérera
  tous ces ensembles inclus dans $(\widehat{A} _K)^r$.

  On pose $E ^{(m)}:=\{ e \in E/\ \forall \underline{k},\
  \underline{\partial} ^{<\underline{k}> _{(m)}} e \in E \cap (A ^\dag ) ^r\}$.
  Il résulte de \cite[2.2.4.(iii)]{Be1},
  que $E ^{(m)}$ est un sous-$A^\dag $-module de $E \cap (A ^\dag ) ^r$.
Il en découle que l'ensemble $E ^{(m)}$ est un sous-$\Gamma(Y^\dag,  \D ^{(m)} _{Y ^\dag})$-module de $E$,
de type fini sur $A ^\dag$.

  Notons $\parallel -\parallel$ la norme spectrale sur $\widehat{A} _K$. Pour tout $a \in \widehat{A} _K$,
  $\parallel  a \parallel \leq 1$ si et seulement si $a \in \widehat{A}$ (\cite[2.4.2]{Be1}).
  Si $\smash{\widehat{A} } ^N \twoheadrightarrow \widehat{E} \cap (\widehat{A})^r$
  est un homomorphisme surjectif, celui-ci induit (on applique $\otimes _\V K$) le suivant
  $\smash{\widehat{A} _K} ^N \twoheadrightarrow \widehat{E} $.
 On notera $\parallel -\parallel$ la norme quotient (qui est une norme de Banach)
 induite par celui-ci. Pour tout
 $e \in \widehat{E}$, si $\parallel e \parallel \leq 1$ alors $e \in \widehat{E} \cap (\widehat{A})^r$.

  Pour $k \in \N$, posons $k = p ^m q _k ^{(m)} + r _k ^{(m)}$. D'après \cite[2.4.3.(i)]{Be1},
  il existe $\eta <1$, $c \in \R$ tels que $|q _k ^{(m)}!|\leq c \eta ^k$ pour tout $k$.
  Pour tous $e \in \widehat{E}$ et $\underline{k}\in \N ^d$,
  $\parallel \underline{\partial} ^{<\underline{k}> _{(m)}} e \parallel=
  |q _{\underline{k}} ^{(m)} !|\parallel \underline{\partial} ^{[\underline{k}]} e \parallel
  \leq c ^d \eta ^{|\underline{k}|} \parallel \underline{\partial} ^{[\underline{k}]} e \parallel$.
  Comme $\widehat{E}$ est un isocristal convergent sur $Y _0$, ce dernier terme tend vers $0$ lorsque
  $|\underline{k}|$ tend vers l'infini.
  Il existe donc un entier $N$
  tel que 
  $|\underline{k}|\geq N$ implique
  $\underline{\partial} ^{<\underline{k}>_{(m)}}  e \in \widehat{E} \cap (\widehat{A})^r$.
  Comme $E$ est stable par l'action des $\underline{\partial} ^{<\underline{k}>_{(m)}}$ et puisque
  que $E \cap \widehat{E} \cap (\widehat{A})^r = E \cap (A ^\dag _K )^r \cap (\widehat{A})^r =
  E \cap (A ^\dag  )^r $ (la dernière égalité résultant de
  \cite[Corollaire de la Proposition 2]{Etesse-descente-etale}), on en déduit que
  $\forall e \in E$, $\exists N \in \N$, $\forall \underline{k}$ tel que $|\underline{k}|\geq N$,
  $\underline{\partial} ^{<\underline{k}>_{(m)}}  e \in E \cap (A ^\dag  )^r$. Il en dérive que le morphisme
  canonique $E ^{(m)} _\Q \rightarrow E$ est un isomorphisme.

  Comme le morphisme
$ \O _{Y^\dag } \otimes _{\Gamma (Y ^\dag ,\O _{Y^\dag })}E ^{(m)} \rightarrow
\D ^{(m)} _{Y ^\dag } \otimes _{\Gamma (Y ^\dag ,\D ^{(m)}_{Y^\dag })}E ^{(m)}$ est un isomorphisme,
le faisceau $\E ^{(m)} := \D ^{(m)}_{Y^\dag } \otimes _{\Gamma (Y ^\dag ,\D ^{(m)}_{Y^\dag })}E ^{(m)}$ répond
à la question.
\hfill \hfill \qed \end{proof}

\begin{rema}
  De manière analogue à \cite[3.1.3]{Be0}, on voit qu'un
  $\D ^{(m)} _{Y ^\dag}$-module, cohérent en tant que $\O _{Y ^\dag}$-module
  est alors cohérent en tant que $\D ^{(m)} _{Y ^\dag}$-module.
\end{rema}

\begin{coro}\label{propmodeisocaffi-cor}
  Avec les notations de \ref{propmodeisocaffi},
  le morphisme $\O _{\Y,\Q} \otimes _{\O _{Y ^\dag ,\Q}}\E
  \rightarrow \D ^\dag _{\Y,\Q} \otimes _{\D _{Y ^\dag ,\Q}}\E$ est un isomorphisme.

  De plus, soient $\widehat{E}$ l'isocristal convergent sur $Y _0$ induit par $E$
  et $\widehat{\E} :=\D ^\dag _{\Y,\Q} \otimes _{\D _{Y ^\dag ,\Q}}\E$. Alors,
  $\widehat{\E} \riso \sp _* (\widehat{E})$,
  où $\mathrm {sp}$ : $\Y _K \rightarrow \Y$ est le morphisme de spécialisation.
\end{coro}
\begin{proof}
 En notant $\smash{\widehat{\E}} ^{(m)}$ le complété $p$-adique
  de $\E ^{(m)}$, on dispose du diagramme commutatif
  $$\xymatrix @R=0,1cm {
  {\O _{\Y} \otimes _{\O _{Y ^\dag }}\E ^{(m)}}
  \ar[rr]  \ar[rd] _-{\widetilde{{}\hspace{0,3cm}  }}
  &&
  {\smash{\widehat{\D}} ^{(m)} _{\Y} \otimes _{\D _{Y ^\dag }} \E ^{(m)}}
  \ar[ld] ^-{\widetilde{{}\hspace{0,3cm}  }}
  \\
  &
  {\smash{\widehat{\E}} ^{(m)},}
  }$$
dont les flèches obliques sont des isomorphismes. Celle horizontale l'est donc aussi.
Il en dérive que le morphisme canonique $\O _{\Y,\Q} \otimes _{\O _{Y ^\dag ,\Q}}\E
  \rightarrow \smash{\widehat{\D}} ^{(m)} _{\Y,\Q} \otimes _{\D _{Y ^\dag ,\Q}}\E$ est un isomorphisme.
  On obtient la première partie du corollaire en passant à la limite sur le niveau.
 La seconde partie s'en déduit.
\hfill \hfill \qed \end{proof}

\subsection{$\D$-modules arithmétiques associés aux isocristaux surconvergents sur les schémas affines et lisses}

  Soient $P ^\dag$ un $\V$-schéma formel faible lisse, $T _0$ un diviseur de $P _0$, $U ^\dag$ l'ouvert de
$P ^\dag $ complémentaire de $T _0$, $j$ : $U ^\dag \hookrightarrow P ^\dag$ l'immersion
  ouverte et $v$ : $Y ^\dag \hookrightarrow U ^\dag$ une immersion fermée
de $\V$-schémas formels faibles.
On suppose en outre $Y ^\dag$ affine et lisse et
on note $X _0$ l'adhérence schématique de
$Y _0$ dans $P _0$ et
$(F\text{-})\mathrm{Coh} ( \PP, T _0 ,X _0)$, la catégorie des
$(F\text{-}) \D ^\dag _{\PP} (\hdag T _0) _\Q$-modules cohérents à support dans $X _0$.

\begin{vide}\label{notations-isocr-sub}
Soient $E$ un isocristal surconvergent sur $Y _0$ et $\E := \sp _* (E)$.
Choisissons $\E ^{(0)}$, un $\D ^{(0)} _{Y ^\dag }$-module globalement de présentation,
$\O _{Y ^\dag }$-cohérent et vérifiant $\E ^{(0)} _\Q \riso \E$ (voir \ref{propmodeisocaffi}).
  On garde les notations de \ref{propmodeisocaffi-cor}
concernant l'isocristal convergent $\widehat{E}$ sur $Y _0$ induit par
  $E$ et $\widehat{\E} :=\D ^\dag _{\Y,\Q} \otimes _{\D _{Y ^\dag ,\Q}}\E$.

Grâce à la proposition \ref{globpresenfiniu+},
le $\D ^{(0)} _{U ^\dag }$-module
$v _+ (\E ^{(0)})$ est globalement de présentation finie.
D'après \ref{j*->Rj*iso}.3, il en dérive que
$j _* v _+ (\E ^{(0)})$ est un $j _* \D _{U^\dag} ^{(0)}$-module globalement de présentation finie.
On en déduit que le faisceau
$$ \sp _{Y ^{\dag}  \hookrightarrow U ^{\dag},T _0 , + }(E) :=
\D ^\dag _{\PP} (\hdag T _0) _\Q \otimes _{j _* \D _{U^\dag} ^{(0)}} j_* v _+ (\E ^{(0)}),$$
est un $ \D ^\dag _{\PP} (\hdag T _0) _\Q$-module globalement de présentation finie à support dans $X _0$.
On remarque que celui-ci est indépendant du choix de $\E ^{(0)}$.
On dispose ainsi d'un foncteur
$\sp _{Y ^{\dag}  \hookrightarrow U ^{\dag},T _0 , + }$ :
$\mathrm{Isoc} ^\dag (Y_0/K) \rightarrow \mathrm{Coh} ( \PP, T _0 ,X _0)$.
Lorsque $T _0 $ est vide, on omettra comme d'habitude de l'indiquer.

\end{vide}

\begin{vide}
  Avec les notations de \ref{notations-isocr-sub}, soient $T' _0 \supset T _0$ un diviseur de $P _0$,
  $U ^{\prime \dag} := P ^{\dag} \setminus T '_0$, $Y ^{\prime \dag} := Y ^{\dag} \setminus T ' _0$
  et $j ' $ : $U ^{\prime \dag} \subset P ^{\dag}$.
  Grâce à \ref{compiminv} (utilisé lorsque $f$ est l'identité), on vérifie que l'on dispose de l'isomorphisme
    fonctoriel en $E\in \mathrm{Isoc} ^\dag (Y_0/K)$
  $(\hdag T' _0 ) \sp _{Y ^{\dag}  \hookrightarrow U ^{\dag},T _0 , + } (E)
  \riso \sp _{Y ^{\prime \dag}  \hookrightarrow U ^{\prime \dag},T '_0 , + } (j ^{\prime \dag} E)$.
\end{vide}

\begin{rema}
  Soit $Y ^\dag $ un $\V$-schéma formel faible affine et lisse.
Grâce à Elkik (\cite{elkik}), il existe un
$\V$-schéma affine et lisse $Y$ dont le complété faible est isomorphe à $Y ^\dag$.
Il existe alors une immersion fermée $Y \hookrightarrow \A _\V ^r$.
En notant $P := \P^r _\V$, $U := \A _\V ^r$ et $T$ le diviseur $P \setminus U$, on obtient
une immersion fermée $v$ : $Y ^\dag\hookrightarrow U ^\dag$ de $\V$-schémas formels faibles lisses
et l'inclusion $j $ : $U ^\dag \subset P ^\dag$. Cette situation géométrique est un cas particulier
de celle de la section.
\end{rema}

\begin{lemm}\label{spyu+rest}
  Avec les notations de \ref{notations-isocr-sub}, on dispose des isomorphismes canoniques
$\sp _{Y ^{\dag}  \hookrightarrow U ^{\dag},T _0 , + } (E) |_{\U}
\riso \sp _{Y _0  \hookrightarrow \U  + } (\widehat{E})
\riso v ^\dag _+ \widehat{\E}$.
\end{lemm}
\begin{proof}
  Par \ref{prop-comm-^},
$\sp _{Y ^{\dag}  \hookrightarrow U ^{\dag},T _0 , + } (E) |_{\U} =
\D ^\dag _{\U,\Q }\otimes _{ \D _{U^\dag} ^{(0)}} v _+ (\E ^{(0)}) \linebreak \riso
\D ^\dag _{\U,\Q } \otimes _{\widehat{\D} ^{(0)} _{\U }}
\hat{v} _+ (\widehat{\D} ^{(0)} _{\Y }   \otimes _{ \D ^{(0)} _{Y ^\dag }} \E ^{(0)})$.
Par \cite[4.3.8]{Beintro2}, il en découle l'isomorphisme
$\sp _{Y ^{\dag}  \hookrightarrow U ^{\dag},T _0 , + }(E) |_{\U}
\riso v ^\dag _+ (\D ^\dag _{\Y ,\Q}   \otimes _{ \D ^{(0)} _{Y ^\dag }} \E ^{(0)})$.
On conclut via $\widehat{\E} \riso \D ^\dag _{\Y ,\Q}   \otimes _{ \D ^{(0)} _{Y ^\dag }} \E ^{(0)}$
(\ref{propmodeisocaffi-cor}).
\hfill \hfill \qed \end{proof}
\begin{prop}\label{sp+exacfidel}
Le foncteur
  $\sp _{Y ^{\dag}  \hookrightarrow U ^{\dag},T _0 , + }$ est exact et fidèle.
\end{prop}
\begin{proof}
D'après \ref{j*->Rj*iso}, le morphisme canonique
  $j _* v _+ (\E ^{(0)}) \rightarrow \R j_* v_+ (\E ^{(0)})$ est un isomorphisme. Le foncteur
  $E \mapsto j _* \D _{U^\dag, \Q}   \otimes _{j _* \D _{U^\dag} ^{(0)}} j_* v _+ (\E ^{(0)})$
  (qui ne dépend pas du choix de $\E ^{(0)}$) est donc exact.
  Comme l'extension $\D _{\PP} (\hdag T _0) _\Q \rightarrow \D ^\dag _{\PP} (\hdag T _0) _\Q$ est plate
à droite et à gauche,
il résulte de \ref{corofidplatpq} que
$\D ^\dag _{\PP} (\hdag T _0) _\Q$ est $j _* \D _{U^\dag,\Q}$-plat à droite et à gauche.
Le foncteur $\sp _{Y ^{\dag}  \hookrightarrow U ^{\dag},T _0 , + }$ est donc exact.

Traitons à présent la fidélité.
Via \cite[4.3.12]{Be1}, il suffit d'établir celle de
$E \mapsto \sp _{Y ^{\dag}  \hookrightarrow U ^{\dag},T _0 , + } (E)|_{\U}\riso v ^\dag _+ \widehat{\E}$.
Or, le foncteur $E \mapsto \widehat{E}$ est fidèle (cela résulte de \cite[2.1.11]{Berig}),
$\sp _*$ est pleinement fidèle sur la catégorie des isocristaux convergents sur $Y _0$ et
$v ^\dag _+$ est pleinement fidèle sur la catégorie des
$\D ^\dag _{\Y ,\Q}$-modules cohérents. D'où le résultat.
\hfill \hfill \qed \end{proof}

\begin{vide}\label{plfidked}
Lorsque $X _0 $ est lisse,
 on dispose (voir \cite[2.2.17]{caro_surholonome})
  du foncteur canonique $\sp _{X _0\hookrightarrow \PP,T_0 +}$ :
  $(F\text{-})\mathrm{Isoc} ^\dag (Y_0,X_0/K) \rightarrow (F\text{-})\mathrm{Coh} (\PP , T_0,X _0)$.
  Ce dernier est pleinement fidèle avec et sans Frobenius et on notera
  ($F$-)$\mathrm{Isoc} ^{\dag}( \PP, T_0, X_0/K)$ son image essentielle (on omet d'indiquer $T _0$ lorsque
  celui-ci est vide).
  On dispose du diagramme essentiellement commutatif (i.e., les deux foncteurs composés sont
  canoniquement isomorphes) :
  $$\xymatrix @R=0,4cm {
  { F\text{-}\mathrm{Isoc} ^{\dag}( \PP, T_0, X_0/K)}
  \ar[r] ^{| _{\U}}
  &
  { F\text{-}\mathrm{Isoc} ^{\dag}( \U , Y _0/K)}
  \\
 {F \text{-} \mathrm{Isoc} ^{\dag}(Y _0 /X _0/K)}
 \ar[r] \ar[u] ^{\sp _{X _0 \hookrightarrow \PP , T _0,+}} _\cong
 &
  {F \text{-} \mathrm{Isoc} (Y _0 /K),}
  \ar[u] ^{\sp _{Y _0 \hookrightarrow \U,+}} _\cong
  }$$
  dont les foncteurs verticaux sont des équivalences de catégorie
  et où les foncteurs horizontaux sont les restrictions canoniques.
  D'après Kedlaya, 
  $F \text{-} \mathrm{Isoc} ^{\dag}(Y _0 /X _0/K) \rightarrow
  F \text{-} \mathrm{Isoc} (Y _0 /K)$ est pleinement fidèle (\cite{kedlaya_full_faithfull}).
  Il en est donc de même de $|_\U$ :
  $F\text{-}\mathrm{Isoc} ^{\dag}( \PP, T_0, X_0/K)
  \rightarrow F\text{-}\mathrm{Isoc} ^{\dag}( \U , Y _0/K)$.
\end{vide}

\begin{vide}\label{notaplfidfrobgen}
Lorsque l'on travaille avec les $F$-isocristaux surconvergents, le théo-\linebreak rème de pleine fidélité de Kedlaya
(\cite{kedlaya_full_faithfull}) est un outil très puissant. Le théorème ci-après en fournit une illustration.
Définissons d'abord la catégorie $F\text{-} \mathrm{Isoc} ^{\dag *} (\PP, T _0, X _0/K)$.
Ses objets sont les couples $(\E,\phi)$, où
  $\E \in \mathrm{Coh} ( \PP, T _0 ,X _0)$ tel que $\E |_\U \in \mathrm{Isoc} ^{\dag } (\U,Y _0)$
et $\phi$ : $F ^* \E |_\U \riso \E |_\U$ dans $\mathrm{Isoc} ^{\dag } (\U,Y _0)$.
  Ses flèches $(\E _1 ,\phi _1)\rightarrow (\E _2 ,\phi _2)$ sont les morphismes
  $\E _1 \rightarrow \E _2 $
$\D ^\dag _{\PP} (\hdag T _0) _\Q$-linéaires
   dont la restriction à $\U$ commute à $\phi _1$ et $\phi _2$.
  Par \ref{spyu+rest}, le foncteur $\sp _{Y ^{\dag}  \hookrightarrow U ^{\dag},T _0 , + }$ induit 
  $F\text{-}\mathrm{Isoc} ^\dag (Y_0/K) \rightarrow F\text{-} \mathrm{Isoc} ^{\dag *} (\PP, T _0, X _0)$,
celui-ci étant fidèle grâce à \ref{sp+exacfidel}.

\end{vide}

\begin{theo}\label{plfidfrobgen}
  Avec les notations de \ref{notaplfidfrobgen}, le foncteur $\sp _{Y ^{\dag}  \hookrightarrow U ^{\dag},T _0 , + }$ :
  $F\text{-}\mathrm{Isoc} ^\dag (Y_0/K) \rightarrow  F\text{-} \mathrm{Isoc} ^{\dag *} (\PP, T _0, X _0/K)$
  est pleinement fidèle.
\end{theo}
\begin{proof}
  On dispose du diagramme canonique essentiellement commutatif
  $$\xymatrix @C=2cm @R=0,3cm {
  {F \text{-} \mathrm{Isoc} ^{\dag}(Y _0 /K)}
  \ar[d] \ar[r] ^{\sp _{Y ^{\dag}  \hookrightarrow U ^{\dag},T _0 , + }}
  &
  { F\text{-} \mathrm{Isoc} ^{\dag *} (\PP, T _0, X _0)}
  \ar[d] ^{|_{\U}}
  \\
 {F \text{-} \mathrm{Isoc} (Y _0 /K)}
 \ar[r] ^{\sp _{Y _0 \hookrightarrow \U,+}} _\cong
 &
 { F\text{-} \mathrm{Isoc} ^{\dag } (\U,Y _0).}
  }$$
Comme $\sp _{Y ^{\dag}  \hookrightarrow U ^{\dag},T _0 , + }$ et $|_{\U}$ sont des foncteurs
fidèles et que les deux autres sont pleinement fidèles, il en résulte que
$\sp _{Y ^{\dag}  \hookrightarrow U ^{\dag},T _0 , + }$ est pleinement fidèle.
\hfill \hfill \qed \end{proof}

\begin{rema}
\label{isocanfrob}
   Si $E\in (F\text{-})\mathrm{Isoc} ^\dag (Y_0/K)$, le théorème \ref{plfidfrobgen}
  n'implique pas directement que l'on dispose d'un isomorphisme canonique de commutation à Frobenius
  $\sp _{Y ^{\dag}  \hookrightarrow U ^{\dag},T _0 , + } (F ^* E)
  \riso F ^* \sp _{Y ^{\dag}  \hookrightarrow U ^{\dag},T _0 , + } (E)$.

  On dit que
  $\phi$ : $\sp _{Y ^{\dag}  \hookrightarrow U ^{\dag},T _0 , + } (F ^* E)
  \riso F ^* \sp _{Y ^{\dag}  \hookrightarrow U ^{\dag},T _0 , + } (E)$
 est
  {\og l'isomorphisme canonique de commutation à Frobenius \fg} s'il induit le diagramme commutatif suivant
    \begin{equation}
    \label{isocanfrobdiag}
\xymatrix @R=0,3cm {
{\sp _{Y ^{\dag}  \hookrightarrow U ^{\dag},T _0 , + } (F ^* E) |_{\U}}
\ar[d] _-\sim ^{\phi |_{\U }} \ar[r] ^-\sim
&
{\sp _{Y _0 \hookrightarrow \U + } (F ^* \widehat{E}) }
\ar[r] ^-\sim
&
{v ^\dag _+ F ^* \widehat{\E}}
\ar[d] ^-\sim
\\
{ F ^* (\sp _{Y ^{\dag}  \hookrightarrow U ^{\dag},T _0 , + } (E) |_{\U})}
\ar[r] ^-\sim
&
{ F ^* \sp _{Y _0 \hookrightarrow \U + } (\widehat{E}) }
\ar[r] ^-\sim
&
{F ^* v ^\dag _+ \widehat{\E},}
}
  \end{equation}
  où les isomorphismes horizontaux sont \ref{spyu+rest} et celui-ci de droite est l'isomorphisme
  canonique de commutation à Frobenius de l'image directe (voir \cite{Beintro2}).
Le diagramme \ref{isocanfrobdiag} implique l'unicité de cet isomorphisme. Tout le problème est
d'établir son existence. On remarque, via cette unicité, que son existence est locale en $\PP$.
Nous la prouverons lorsque $X _0$ est lisse (voir \ref{sp+froblisse}) et surtout
  lorsque $Y _0$ se désingularise localement idéalement (voir \ref{theo-ideal->oK}).
\end{rema}

\begin{conj}\label{conj}
Notons $\rho _{Y _0,X _0}$ : $\mathrm{Isoc} ^\dag (Y_0/K) \rightarrow
\mathrm{Isoc} ^\dag (Y_0,X _0/K)$
et $\mathrm{cv}$ : $\mathrm{Isoc} ^\dag (Y_0,X _0/K) \rightarrow \mathrm{Isoc} (Y_0/K)$
les foncteurs restrictions (de même en rajoutant $F\text{-}$).
Il est conjectural que ceux-ci soient pleinement fidèles (voir \cite{tsumono}).
On conjecture alors les deux suivantes :
\begin{enumerate} \label{sp+essent}
  \item [(a)]Lorsque $X _0$ est lisse, pour tout isocristal surconvergent $E$ sur $Y _0$,
on dispose d'un isomorphisme fonctoriel en $E$
$$\sp _{Y ^{\dag}  \hookrightarrow U ^{\dag},T _0 , + } (E)\riso
\sp _{X _0\hookrightarrow \PP,T _0 +} \circ \rho _{Y _0,X _0} (E).$$
\item [(b)] \label{conjpf}
  Le foncteur $\sp _{Y ^{\dag}  \hookrightarrow U ^{\dag},T _0 , + }$ est pleinement fidèle.
\end{enumerate}

Le théorème \ref{conjvalF} valide la conjecture $(a)$ lorsque l'on dispose d'une structure de Frobenius.
\end{conj}

\begin{lemm} \label{lemmindXlisse}
  Soient $\PP$ un $\V$-schéma formel lisse, $T _0 \subset T _0'$ deux diviseurs de $P$ et
$\E$ un ($F$-)$\D ^\dag _{\PP} (\hdag T _0) _\Q$-module cohérent.
Il existe au plus une structure de ($F$-)$\D ^\dag _{\PP} (\hdag T _0') _\Q$-module cohérent sur $\E$
prolongeant sa structure canonique de ($F$-)$\D ^\dag _{\PP} (\hdag T _0) _\Q$-module.
De plus, l'existence d'une telle structure est locale en $\PP$.
\end{lemm}
\begin{proof}
  On dispose du morphisme
  $\E \rightarrow \D ^\dag _{\PP} (\hdag T _0') _\Q \otimes _{\D ^\dag _{\PP} (\hdag T _0) _\Q} \E$.
  Par \cite[4.3.12]{Be1}, celui-ci est un isomorphisme si et seulement si
$\E$ est muni d'une structure de ($F$-)$\D ^\dag _{\PP} (\hdag T _0') _\Q$-module cohérent prolongeant
  sa structure canonique de ($F$-)$\D ^\dag _{\PP} (\hdag T _0) _\Q$-module.
\hfill \hfill \qed \end{proof}
\begin{lemm}
 \label{propindXlisse}
  Soient $\PP$ un $\V$-schéma formel lisse, $ X_0$, $ X_0'$ des sous-schémas fermés lisses de $P$, $T _0 $ et $T _0'$ deux
  diviseurs de $P$ tels que $ X_0 \setminus  T _0 = X_0'\setminus T _0'$.

Soit $\E$ un ($F$-)$\D ^\dag _{\PP} (\hdag T _0) _\Q$-module 
tel qu'il existe un ($F$-) isocristal $E$ sur $ X_0 \setminus T _0$ surconvergent le long $T _0 \cap  X_0$ et
un isomorphisme ($F$-)$\D ^\dag _{\PP} (\hdag T _0) _\Q$-linéaire $\E \riso \sp _{ X_0 \hookrightarrow \PP, T _0,+} (E)$.

La structure de ($F$-)$\D ^\dag _{\PP} (\hdag T _0) _\Q$-module de $\E$ se prolonge en
une structure de ($F$-)$\D ^\dag _{\PP} (\hdag T _0 \cup T _0') _\Q$-module cohérent.
De plus, il existe un ($F$-) isocristal $E'$ sur $ X_0 '\setminus T _0'$ surconvergent le long $T _0 '\cap  X_0 '$
et un isomorphisme ($F$-) $\D ^\dag _{\PP} (\hdag T _0 ') _\Q$-linéaire
$\E \riso \sp _{ X_0 '\hookrightarrow \PP, T _0',+} (E')$.
\end{lemm}
  \begin{proof}

    Traitons d'abord le cas où $T _0'=T _0$.
    D'après la caractérisation de l'image essentielle de $\sp _{ X_0 '\hookrightarrow \PP, T _0',+}$,
    il suffit de prouver que $\E $ est à support dans $ X_0'$, i.e., $(\hdag  X_0 ') (\E)=0$.
    Grâce à \cite[2.2.9]{caro_surcoherent},
    comme $\R \underline{\Gamma} ^\dag _{T _0} ((\hdag  X_0') \E)=0$,
    $(\hdag  X_0') \E =0 $ si et seulement si $((\hdag  X_0') \E ) |_{\PP \setminus T _0} =0$. 
    Or, $[(\hdag  X_0') (\E)] |_{\PP \setminus T _0} \riso (\hdag  X_0' \setminus T _0) (\E |_{\PP \setminus T _0 })
    = (\hdag  X_0  \setminus T _0) (\E |_{\PP \setminus T _0 }) =0$.

    Supposons maintenant $ X_0 =  X_0'$. Par \ref{lemmindXlisse}, le fait que
    $\E$ soit de surcroît un ($F$-)$\D ^\dag _{\PP} (\hdag T _0\cup T _0') _\Q$-module cohérent est local
    en $\PP$. Supposons donc que $ X_0 \hookrightarrow P$ se relève en un
    morphisme $u$ : $\ X_0 \rightarrow \PP$ de $\V$-schémas formels lisses. Comme
    les foncteurs $u _+$, $u _{T _0,+}$ et $u _{T _0\cup T _0',+}$ sont isomorphes et
    puisque $(T _0 \cup T _0') \cap  X_0 = T _0 \cap  X_0$,
    $\sp _{ X_0 \hookrightarrow \PP,T _0,+} (E) \riso u _{T _0\cup T _0',+} \sp _* (E)$.
    La cohérence se préservant par image directe par un morphisme propre,
    $\sp _{ X_0 \hookrightarrow \PP,T _0,+} (E)$ est donc un ($F$-)$\D ^\dag _{\PP} (\hdag T _0\cup T _0') _\Q$-module cohérent.
    De même, $\sp _{ X_0 \hookrightarrow \PP,T _0,+} (E)$ est
    dans l'image essentielle de $\sp _{ X_0 \hookrightarrow \PP,T _0',+}$.

    En remarquant que $ X_0 \setminus T _0 =  X_0 \setminus (T _0 \cup T _0') =  X_0 '\setminus (T _0 \cup T _0')
    =  X_0 '\setminus T _0'$, on vérifie que le cas général se déduit des deux précédents.
\hfill \hfill \qed \end{proof}

\begin{prop}\label{sp+froblisse}
  On suppose $X _0 $ lisse, $P _0$ séparé et on se donne
  $E$ un isocristal sur $Y _0$ surconvergent le long de $T _0 \cap X _0$.
  L'isomorphisme canonique
  $\phi$ : $\sp _{Y ^{\dag}  \hookrightarrow U ^{\dag},T _0 , + } (F ^* E)
  \riso F ^* \sp _{Y ^{\dag}  \hookrightarrow U ^{\dag},T _0 , + } (E)$
  de commutation à Frobenius existe (voir \ref{isocanfrob}),
  i.e., il induit le diagramme commutatif ci-après
  \begin{equation}
    \label{sp+froblissediagp}
\xymatrix @R =0,3cm {
{\sp _{Y ^{\dag}  \hookrightarrow U ^{\dag},T _0 , + } (F ^* E) |_{\U}}
\ar[r] ^{\phi |_{\U }}  _-\sim \ar[d] ^-\sim
&
{ F ^* (\sp _{Y ^{\dag}  \hookrightarrow U ^{\dag},T _0 , + } (E) |_{\U})}
\ar[d] ^-\sim
\\
{v ^\dag _+ F ^* \widehat{\E}}
\ar[r]^-\sim
&
{F ^* v ^\dag _+ \widehat{\E},}
}
  \end{equation}
  où les isomorphismes verticaux sont \ref{spyu+rest} et celui du bas est l'isomorphisme de commutation
  à Frobenius de l'image directe.
  Lorsque $E$ est un $F$-isocristal surconvergent,
  $\sp _{Y ^{\dag}  \hookrightarrow U ^{\dag},T _0 , + }(E)$ est donc muni d'une structure canonique
  de $F\text{-} \D ^\dag _{\PP} (\hdag T _0) _\Q$-module.

\end{prop}
\begin{proof}
  Fixons $F $ : $Y ^\dag \rightarrow Y ^\dag $ un relèvement de Frobenius.
  Notons $P ^{\prime \dag} := P ^\dag \times P ^\dag$, $U ^{\prime \dag} := U ^\dag \times U ^\dag$,
  $j '$ : $U ^{\prime \dag} \subset P ^{\prime \dag}$,
  $T '_0:= P ' _0 \setminus U ' _0$,
  $f _1$ et $f _2$ : $P ^{\prime \dag} \rightarrow P ^\dag$
  (resp. $g _1$ et $g _2$ : $U ^{\prime \dag} \rightarrow U ^\dag$),
  les projections respectives à gauche et à droite, $h _1$ : $U ^\dag \times P ^\dag \rightarrow U ^\dag$ et
  $v ' $ l'immersion fermée $(v, v\circ F)$ : $Y ^\dag \hookrightarrow U ^{\prime \dag}$.
  D'après \ref{compiminv},
  $\D ^\dag _{\PP '} (\hdag T  '_0 )_\Q \otimes _{j ' _* \D ^{(0)} _{U'}}
  j '_*g _2 ^! (v _+ (\E ^{(0)}))
  \riso
  f ^! _{2,T'_0,T_0} (\D ^\dag _{\PP} (\hdag T _0) _\Q \otimes _{j _* \D _{U^\dag} ^{(0)}} j_* v _+ (\E ^{(0)}))$.
  Or, on dispose des morphismes :
  \begin{equation}
    \label{sp+froblissediag1}
  v ' _+ F ^! ( \E ^{(0)}) \riso v ' _+ F ^! v ^! v _+  ( \E ^{(0)}) \riso
  v ' _+  v ^{\prime !} g _2 ^! v _+  ( \E ^{(0)} ) \rightarrow g _2 ^! v _+  ( \E ^{(0)}).
  \end{equation}
  Il en dérive par fonctorialité et composition
  \begin{equation}
    \label{sp+froblissediag2}
  \D ^\dag _{\PP '} (\hdag T  '_0 )_\Q
 \underset{j ' _* \D ^{(0)} _{U'}}{\otimes}
  j '_*(v ' _+ F ^! ( \E ^{(0)} ))
  \rightarrow
  f ^! _{2,T_0',T_0} (\D ^\dag _{\PP} (\hdag T _0) _\Q \otimes _{j _* \D _{U^\dag} ^{(0)}} j_* v _+ (\E ^{(0)})).
  \end{equation}

Notons $\delta _{F_0} $ le morphisme $(id, F _0)$ : $P _0 \hookrightarrow P _0 \times P _0$,
où $F_0$ est Frobenius.
On identifiera $P _0$, $ U_0$, $X _0$, $Y_0$ à des sous-schémas de $P _0 \times P _0$
via $\delta _{F_0}$ (on prendra garde de ne pas les identifier via l'immersion diagonale par exemple).

Comme il existe un isocristal $\widetilde{E}$ sur $Y _0$ surconvergent le long de $X _0 \cap T _0$ tel que
$\D ^\dag _{\PP} (\hdag T _0) _\Q \otimes _{j _* \D _{U^\dag} ^{(0)}} j_* v _+ (\E ^{(0)})
\riso \sp _{X _0 \hookrightarrow \PP, T_0,+} (\widetilde{E})$, celui-ci
est $\D ^\dag _{\PP} (\hdag T _0) _\Q$ surcohérent (\ref{theo-surcohlisse}).
Puisque
$\D ^\dag _{\PP '} (\hdag T  '_0 )_\Q \otimes _{j ' _* \D ^{(0)} _{U'}}  j '_*(v ' _+ F ^! ( \E ^{(0)} ))$
est un $\D ^\dag _{\PP '} (\hdag T  '_0 )_\Q $-module cohérent à support dans $X _0$ (et donc dans $P _0$),
il dérive alors de \ref{sp+froblissediag2} le morphisme
dans $D ^\mathrm{b} _\mathrm{coh} (\D ^\dag _{\PP '} (\hdag T  '_0 )_\Q)$ :
\begin{gather}
\notag
  \D ^\dag _{\PP '} (\hdag T  '_0 )_\Q \otimes _{j ' _* \D ^{(0)} _{U'}}
  j '_*(v ' _+ F ^! ( \E ^{(0)} ))
\\
\rightarrow
 \R \underline{\Gamma} _{P _0} ^\dag \circ
 f ^! _{2,T'_0,T_0} (\D ^\dag _{\PP} (\hdag T _0) _\Q \otimes _{j _* \D _{U^\dag} ^{(0)}} j_* v _+ (\E ^{(0)})).
\label{sp+froblissediag3}
  \end{gather}

Comme $\delta _{F_0} ^{-1} (f_1 ^{-1} (T _0))= T _0$,
$\delta _{F_0} ^{-1} (f_2 ^{-1} (T _0))= T _0 ^{(s)}$
et $T '_0 =f_1 ^{-1} (T _0) \cup f_2 ^{-1} (T _0)$,
$\R \underline{\Gamma} _{P _0} ^\dag \circ ( \hdag f_2 ^{-1} (T _0)) \riso
\R \underline{\Gamma} _{P _0} ^\dag \circ ( \hdag T'_0)$.
On en déduit $\R \underline{\Gamma} _{P _0} ^\dag \circ f ^! _{2,T'_0,T_0} \linebreak \riso
 \R \underline{\Gamma} _{P _0} ^\dag \circ  f ^! _{2,T _0}$.
Comme $\delta _{F_0,f_1 ^{-1} (T _0),+}\riso \delta _{F_0,+} \riso \delta _{F_0,f_2 ^{-1} (T _0),+}$,
on obtient
\begin{gather}\notag
  \R \underline{\Gamma} _{P _0} ^\dag \circ f ^! _{2,T'_0,T_0} \riso
 \R \underline{\Gamma} _{P _0} ^\dag \circ  f ^! _{2,T _0} \riso
  \delta _{F_0,f_2 ^{-1} (T _0),+} \delta _{F_0,f_2 ^{-1} (T _0)} ^! f ^! _{2,T _0}
\\
  \riso
  \delta _{F_0,f_2 ^{-1} (T _0),+} F  ^* _{0,T _0}
  \riso \delta _{F_0,f_1 ^{-1} (T _0),+} F  ^* _0 .\notag
\end{gather}
  Le morphisme \ref{sp+froblissediag3} est donc canoniquement isomorphe à
\begin{gather}
    \notag
  \D ^\dag _{\PP '} (\hdag T  '_0 )_\Q \otimes _{j ' _* \D ^{(0)} _{U'}} j '_*(v ' _+ F ^! ( \E ^{(0)} ))
   \\ \rightarrow
 \delta _{F_0,f_1 ^{-1} (T _0),+} F  ^* _0
 (\D ^\dag _{\PP} (\hdag T _0) _\Q \otimes _{j _* \D _{U^\dag} ^{(0)}} j_* v _+ (\E ^{(0)})).
\label{sp+froblissediag4}
  \end{gather}

Notons $j ''$ l'immersion ouverte $U ^\dag \times P ^\dag \hookrightarrow P ^{\prime \dag}$
et $v ''$ : $Y ^\dag \hookrightarrow U ^\dag \times P ^\dag$
l'immersion fermée composée du graphe de $j\circ v\circ F$ (qui est une immersion fermée
car $P _0$ est séparé) suivi de
$v \times id$ : $Y ^\dag \times P ^\dag \hookrightarrow U ^\dag \times P ^\dag$.
Comme $j '' \circ v''= j ' \circ v'$, $j ''_* \circ v '' _+   \riso j '_* \circ v ' _+ $.
Il en découle que le morphisme d'adjonction
$j '' _* \D ^{(0)} _{U ^\dag \times P ^\dag} \rightarrow j ' _* \D ^{(0)} _{U ^{\prime \dag}}$
induit la flèche
\begin{equation}
    \label{sp+froblissediag5}
    \D ^\dag _{\PP '} (\hdag T  '_0 )_\Q
\underset{j '' _* \D ^{(0)} _{U ^\dag \times P ^\dag}}{\otimes}
    j ''_*(v '' _+ F ^! ( \E ^{(0)} ))
    \rightarrow
    \D ^\dag _{\PP '} (\hdag T  '_0 )_\Q \otimes _{j ' _* \D ^{(0)} _{U ^{\prime \dag}}}
    j '_*(v ' _+ F ^! ( \E ^{(0)} )).
  \end{equation}

Comme $f _1\circ \delta _{F_0} =id$, on vérifie $Y _0 = X _0 \times _{P _0 \times P _0} (U _0\times P _0)$.
$F ^* \E ^{(0)}$ est un
  $\D ^{(0)} _{Y ^\dag }$-module globalement de présentation,
  $\O _{Y ^\dag }$-cohérent et tel que $F ^* \E ^{(0)} _\Q \riso F ^*\E$.
  Il représente donc un modèle de l'isocristal surconvergent $F ^* E$.
 Or, $\D ^\dag _{\PP '} (\hdag f _1 ^{-1} (T_0) )_\Q \otimes _{j '' _* \D ^{(0)} _{U ^\dag \times P ^\dag}}
    j ''_*(v '' _+ F ^* ( \E ^{(0)} ))$ est
    dans l'image essentielle de $\sp _{X _0\hookrightarrow P' _0, f _1 ^{-1} (T_0),+}$.
Par \ref{propindXlisse},
celui-ci est même un
$\D ^\dag _{\PP '} (\hdag T'_0 )_\Q$-module cohérent et donc
le morphisme canonique
\begin{gather}
    \notag
    \D ^\dag _{\PP '} (\hdag f _1 ^{-1} (T_0) )_\Q \otimes _{j '' _* \D ^{(0)} _{U ^\dag \times P ^\dag}}
    j ''_*(v '' _+ F ^! ( \E ^{(0)} ))
    \\ \rightarrow
    \D ^\dag _{\PP '} (\hdag T  '_0 )_\Q \otimes _{j '' _* \D ^{(0)} _{U ^\dag \times P ^\dag}}
    j ''_*(v '' _+ F ^! ( \E ^{(0)} )).\label{sp+froblissediag6}
  \end{gather}
est un isomorphisme.
En composant \ref{sp+froblissediag4}, \ref{sp+froblissediag5} et \ref{sp+froblissediag6},
on obtient
\begin{gather}
\notag
  \D ^\dag _{\PP '} (\hdag f _1 ^{-1} (T_0) )_\Q \otimes _{j '' _* \D ^{(0)} _{U ^\dag \times P ^\dag}}
    j ''_*(v '' _+ F ^! ( \E ^{(0)} ))
\\  \rightarrow
 \delta _{F_0,f_1 ^{-1} (T _0),+} F  ^* _0
 (\D ^\dag _{\PP} (\hdag T _0) _\Q \otimes _{j _* \D _{U^\dag} ^{(0)}} j_* v _+ (\E ^{(0)})).
\label{sp+froblissediag7}
  \end{gather}
Or, d'après \ref{g-+commdag*f+}, on dispose du morphisme
\begin{gather}
\notag
\D ^\dag _{\PP } (\hdag T_0 )_\Q \otimes _{j  _* \D ^{(0)} _{U ^\dag }}
    j _*(h _{1,+} v '' _+ F ^* ( \E ^{(0)} ))
    \\ \rightarrow
      f _{1, T _0,+ }( \D ^\dag _{\PP '} (\hdag f _1 ^{-1} (T_0) )_\Q
      \otimes _{j '' _* \D ^{(0)} _{U ^\dag \times P ^\dag}}
    j ''_*(v '' _+ F ^* ( \E ^{(0)} ))),\label{sp+froblissediag8}
\end{gather}
qui est un isomorphisme grâce à \cite[4.3.12]{Be1}.
Comme $h _{1,+} v '' _+ \riso v _+$ et
$f _{1, T _0,+ } \circ \delta _{F_0,f_1 ^{-1} (T _0),+}\riso id _+$, en utilisant \ref{sp+froblissediag8} et
en appliquant $f _{1, T _0,+ }$ à \ref{sp+froblissediag7}, on obtient alors
\begin{gather}
    \notag
  \D ^\dag _{\PP } (\hdag T_0 )_\Q \otimes _{j  _* \D ^{(0)} _{U ^\dag }}
    j _*(v  _+ F ^* ( \E ^{(0)} ))
  \\ \rightarrow
  F  ^* _0 (\D ^\dag _{\PP} (\hdag T _0) _\Q \otimes _{j _* \D _{U^\dag} ^{(0)}} j_* v _+ (\E ^{(0)})).
\label{sp+froblissediag9}
  \end{gather}

Par \cite[4.3.12]{Be1}, pour conclure la démonstration, il suffit de prouver que
le diagramme \ref{sp+froblissediagp} est commutatif.
Cela est local en $\U$ et on se ramène au cas où $U ^\dag $ est affine.
On notera $F$ : $U ^\dag \rightarrow U ^\dag $ un relèvement de Frobenius
tel que $F \circ v = v \circ F$ (un tel choix est possible grâce à \cite[3.3.2]{arabia-treslisse})
et $\delta _{F}:s= (id,F)$ : $ U ^\dag \hookrightarrow  U ^{\prime \dag}$.
On aura besoin du lemme ci-après.
\begin{lemm*}\label{lef1h1}
Avec les notations ci-dessus,
  pour tout $\D ^\dag _{\PP '} (\hdag f _1 ^{-1} (T_0) )_\Q$-module (ou complexe) cohérent $\E$
  à support dans $P _0$, on dispose d'un isomorphisme canonique
  $(f _{1, T _0,+} (\E) | _{\U} \riso g _{1,+} (\E |_{\U '})$ fonctoriel en $\E$.
\end{lemm*}
\begin{proof}
   On a toujours
  $(f _{1, T _0,+} (\E) | _{\U} \riso h _{1,+} (\E | _{\U \times \PP})$.
  On termine alors la preuve en remarquant que $\E | _{\U \times \PP}$ est à support dans $U _0$
  et l'immersion fermée $U _0 \hookrightarrow U _0 \times P _0$ induite par $\delta _{F _0}$
  se factorise par $U _0 \hookrightarrow U '_0 $.
\hfill \hfill \qed \end{proof}
Grâce au lemme \ref{lef1h1}, en appliquant $g _{1,+} ^\dag$ au morphisme $\theta$ de
\begin{equation}
  \label{deftheta}
\xymatrix @R=0,3cm{
{\D ^\dag _{\U',\Q} \otimes _{\D ^{(0)} _{U'}} (v ' _+ F ^! ( \E ^{(0)}))}
\ar@{.>}[d] ^-\theta \ar@{=}[r]
&
{\D ^\dag _{\U',\Q} \otimes _{\D ^{(0)} _{U'}} (v ' _+ F ^! ( \E ^{(0)}))}
\ar[d]
\\
{\delta _{F +^{\dag}} F ^* (\D ^\dag _{\U,\Q} \otimes _{\D ^{(0)} _{U}} v _+  ( \E ^{(0)}))}
&
{\D ^\dag _{\U',\Q} \otimes _{\D ^{(0)} _{U'}}(  g _2 ^! v _+  ( \E ^{(0)}))}
\ar[d] ^-\sim
\\
{\delta _{F +^{\dag}} \delta _F ^{!^{\dag}} g _2 ^{!^{\dag}} (\D ^\dag _{\U,\Q} \otimes _{\D ^{(0)} _{U}} v _+  ( \E ^{(0)})),}
\ar[r] \ar[u] ^-\sim
&
{g _2 ^{!^{\dag}} (\D ^\dag _{\U,\Q} \otimes _{\D ^{(0)} _{U}} v _+  ( \E ^{(0)}))}
}
\end{equation}
où la flèche de droite du haut dérive de \ref{sp+froblissediag1},
on obtient
(via \ref{prop-comm-^dir} et les isomorphismes $g _{1,+ ^{\dag}} \delta _{F +^{\dag}}  \riso id _+$ et
$g _{1,+}  v ' _{+}  \riso v _{+}$)
la restriction de \ref{sp+froblissediag9} à $\U$, i.e.,
$\D ^\dag _{\U ,\Q} \otimes _{\D ^{(0)} _{U ^\dag }} v  _+ F ^* ( \E ^{(0)} )
  \rightarrow
  F  ^*  (\D ^\dag _{\U,\Q} \otimes _{\D _{U^\dag} ^{(0)}}v _+ (\E ^{(0)}))$.
  En appliquant $\delta _{F +^{\dag}}$
  à l'isomorphisme canonique $v _{+ ^{\dag}} F ^* \widehat{\E} \riso F ^* v _{+ ^{\dag}}  \widehat{\E}$,
  on obtient $v ^{\prime \dag} _+ F ^* \widehat{\E} \riso \delta _{F + ^\dag} F ^* v _{+ ^{\dag}}  \widehat{\E}$,
  que l'on notera $\theta ^\prime$. Il s'agit de vérifier que $\theta$ et
  $\theta ^\prime$ se correspondent modulo les isomorphismes \ref{prop-comm-^dir},
  \ref{prop-comm-^inv2}.
  Or, d'après \cite{caro_courbe}, le morphisme d'adjonction
   $\mathrm{adj}$ : $id \rightarrow v ^{!\dag} v _{+ ^\dag} $ est compatible à Frobenius, i.e.,
  le diagramme canonique
  \begin{equation}
    \xymatrix @R=0,4cm {
  {F ^* (\widehat{\E}) } \ar[r] _-\sim ^-{\mathrm{adj}} \ar[d] ^{\mathrm{adj}} _-\sim &
  { v ^{!\dag} v _{+ ^\dag} F ^*(\widehat{\E})}
  \ar[d] ^-\sim
  \\
  {F ^* v ^{!\dag} v _{+ ^\dag} (\widehat{\E})}
  \ar[r] ^-\sim
  &
  {  v ^{!\dag} F ^*  v _{+ ^\dag}(\widehat{\E}),}
  }
    \end{equation}
  est commutatif. En lui appliquant $v _{+ ^\dag}$, on obtient la commutativité du carré de gauche du diagramme
  \begin{equation}
  \label{v+frobdiag1}
    \xymatrix @R =0,4cm {
  {v _{+ ^\dag} F ^* (\widehat{\E}) } \ar[r] _-{\mathrm{adj}} ^-\sim \ar[d] ^-{\mathrm{adj}} _-\sim
  \ar@/^1,2pc/@{=}[rr]
  &
  {v _{+ ^\dag} v ^{!\dag} v _{+ ^\dag} F ^*(\widehat{\E})}
  \ar[d] ^-\sim   \ar[r] _-{\mathrm{adj}} ^-\sim
  &
  { v _{+ ^\dag} F ^*(\widehat{\E})}
  \ar[d] _-\sim
  \\
  {v _{+ ^\dag} F ^* v ^{!\dag} v _{+ ^\dag} (\widehat{\E})}
  \ar[r] ^-\sim
  &
  { v _{+ ^\dag} v ^{!\dag} F ^*  v _{+ ^\dag}(\widehat{\E})}
  \ar[r] ^-{\mathrm{adj}} _-\sim
  &
  {  F ^*  v _{+ ^\dag}(\widehat{\E}).}
  }
  \end{equation}
  Celle du carré de droite se vérifie par fonctorialité. D'où la commutativité de \ref{v+frobdiag1}.
  Il en résulte alors celle du diagramme
\begin{equation}
  \label{v+frobdiag2}
  \xymatrix @R=0,2cm {
&
{ \delta _{F + ^\dag }  v ^{\dag} _+ F ^*(\widehat{\E})}
  \ar[d] ^-{\mathrm{adj}} _-\sim
&
{ v ^{\prime ^{\dag}} _+ F ^*(\widehat{\E})}
  \ar[d] ^-{\mathrm{adj}} _-\sim
  \ar[l] ^-\sim
\\
{ \delta _{F + ^\dag }  v ^{\dag} _+  v ^{!\dag} F ^* v _{+ ^\dag} (\widehat{\E})}
  \ar@{=}[d]
&
{ \delta _{F + ^\dag }  v ^{\dag} _+ F ^* v ^{!\dag} v _{+ ^\dag} (\widehat{\E})}
  \ar[l] ^-\sim \ar[d] ^-\sim
&
{ v ^{\prime \dag} _+ F ^* v ^{!\dag} v _{+ ^\dag} (\widehat{\E})}
    \ar[d] ^-\sim
    \ar[l] ^-\sim
\\
{ \delta _{F + ^\dag }  v ^{\dag} _+  v ^{!\dag} F ^* v _{+ ^\dag} (\widehat{\E})}
  \ar[d] ^-{\mathrm{adj}}
&
{ \delta _{F + ^\dag }  v ^{\dag} _+ v ^{  !\dag} \delta _{F } ^{! \dag}  g _2 ^{!\dag} v _{+ ^\dag} (\widehat{\E})}
  \ar[d] ^-{\mathrm{adj}}
  \ar[l] ^-\sim
  \ar[r] _-\sim
&
{ v ^{\prime \dag} _+  v ^{\prime  !\dag}  g _2 ^{!\dag} v _{+ ^\dag} (\widehat{\E})}
  \ar[d] ^-{\mathrm{adj}}
\\
{ \delta _{F + ^\dag }  F ^* v _{+ ^\dag} (\widehat{\E})}
&
{ \delta _{F + ^\dag }  \delta _{F } ^{! \dag}  g _2 ^{!\dag} v _{+ ^\dag} (\widehat{\E})}
  \ar[l] ^-\sim
  \ar[r] _-\sim ^-{\mathrm{adj}}
&
{ g _2 ^{!\dag} v _{+ ^\dag} (\widehat{\E})}
}
\end{equation}
En effet, le carré en bas à droite est commutatif par transitivité des morphismes d'adjonction.
Les deux du milieu le sont
par transitivité de l'isomorphisme de commutation à la composition
des images inverses extraordinaires. Enfin, les deux derniers le sont par fonctorialité.
En utilisant la commutativité du diagramme déduit de \ref{v+frobdiag1} par application de
$ \delta _{F + ^{\dag} }$, on vérifie que
le morphisme
$  v ^{\prime ^{\dag}} _+ F ^*(\widehat{\E})
\rightarrow \delta _{F + ^\dag }  F ^* v _{+ ^\dag} (\widehat{\E})$
de \ref{v+frobdiag2} passant par le haut puis la gauche est $\theta ^\prime$.
Or, grâce à \ref{tracecomphat2}, \ref{prop-comm-^dir} et \ref{prop-comm-^inv2}, on vérifie que la flèche
$  v ^{\prime ^{\dag}} _+ F ^*(\widehat{\E})
\rightarrow \delta _{F + ^\dag }  F ^* v _{+ ^\dag} (\widehat{\E})$
de \ref{v+frobdiag2} passant par la droite puis le bas,
correspond à $\theta$ (voir \ref{deftheta}).

\hfill \hfill \qed \end{proof}

\begin{prop}\label{conjvalF}
On suppose $X _0 $ lisse et $P _0$ séparé.
  Pour tout $F$-isocristal $E$ surconvergent sur $Y _0$,
  on dispose d'un isomorphisme
$F\text{-}\D ^\dag _{\PP} (\hdag T _0) _{\Q}$-linéaire
$  \sp _{Y ^{\dag}  \hookrightarrow U ^{\dag},T _0 , + } (E)\riso
\sp _{X _0\hookrightarrow \PP,T _0 +} \circ \rho _{Y _0,X _0} (E)$
et fonctoriel en $E$.
\end{prop}
\begin{proof}
Grâce à \ref{sp+froblisse}, 
$\sp _{Y ^{\dag}  \hookrightarrow U ^{\dag},T _0 , + } (E)$
est muni d'une structure de $F\text{-}\D ^\dag _{\PP} (\hdag T _0) _\Q$-module telle que
  $\sp _{Y ^{\dag}  \hookrightarrow U ^{\dag},T _0 , + } (E) |_{\U}
\riso 
\sp _{Y _0 \hookrightarrow \U,+} (\widehat{E})$
commute à Frobenius.
On conclut grâce à la pleine fidélité du foncteur restriction
$|_\U$ :
$F\text{-}\mathrm{Isoc} ^{\dag}( \PP, T_0, X_0/K)
  \rightarrow F\text{-}\mathrm{Isoc} ^{\dag}( \U , Y _0/K)$ (voir \ref{plfidked}).
\hfill \hfill \qed \end{proof}

Terminons cette section par la proposition suivante.
\begin{prop}\label{sp+plssstabiminv}
Soient $f$ : $P ^{\prime \dag} \rightarrow P ^\dag$ un morphisme lisse de $\V$-schémas formels faibles lisses,
$T _0$ un diviseur de $P _0$, $T ' _0:= f ^{-1} (T _0)$,
$U ^\dag $ (resp. $U ^{\prime \dag}$) l'ouvert de $P ^\dag $ (resp. $P ^{\prime \dag}$)
complémentaire de $T _0$ (resp. $T ' _0$),
$j$ (resp. $j'$) l'immersion ouverte correspondante.
On se donne de plus un morphisme
$b$ : $Y ^{\prime \dag} \rightarrow Y ^{\dag}$ de $\V$-schémas formels faibles affines et lisses,
des immersions fermées
$v$ : $ Y ^\dag \hookrightarrow U ^\dag$
et $v'$ : $ Y ^{\prime \dag} \hookrightarrow U ^{\prime \dag}$,
tels que
$g \circ v' = v \circ b$, où $g$ : $U ^{\prime \dag} \rightarrow U ^\dag$
est le morphisme induit par $f$. Enfin, $X _0$ (resp. $X ' _0$) désigne l'adhérence
de $Y _0$ (resp. $Y ' _0$) dans $P _0$ (resp. $P ' _0$).

Pour tout objet $E \in \mathrm{Isoc} ^\dag (Y_0/K)$ tel que
$\R \underline{\Gamma} ^\dag _{X ' _0} f _{T _0} ^! (\sp _{Y ^{\dag}  \hookrightarrow U ^{\dag},T _0 , + }(E))$
soit à cohomologie $\D ^\dag _{\PP'} (\hdag T '_0) _\Q$-cohérente, il existe un isomorphisme canonique
$\sp _{Y ^{\prime \dag}  \hookrightarrow U ^{\prime \dag},T '_0 , + } (b ^*E)[d_{X '_0/X _0}]
\riso
\R \underline{\Gamma} ^\dag _{X ' _0} f _{T _0} ^!
(\sp _{Y ^{\dag}  \hookrightarrow U ^{\dag},T _0 , + }(E) ).$

En outre, celui-ci commute aux isomorphismes de commutation des images inverses (resp. images
inverses extraordinaires).

En particulier, lorsque $X ' _0 = f ^{-1} (X _0)$,
pour tout $E \in \mathrm{Isoc} ^\dag (Y_0/K)$, on obtient
$\sp _{Y ^{\prime \dag}  \hookrightarrow U ^{\prime \dag},T '_0 , + } (b ^*E)[d_{X '_0/X _0}]
\riso
f _{T _0} ^! (\sp _{Y ^{\dag}  \hookrightarrow U ^{\dag},T _0 , + }(E) )$.
\end{prop}
\begin{proof}
  Il résulte de la \ref{compiminv},
  $\D ^\dag _{\PP'} (\hdag T '_0) _\Q \otimes _{j ' _* \D _{U^{\prime \dag}}} ^{(0)} j '_* g ^! v _+ (\E ^{(0)})
  \linebreak
  \riso
  f _{T _0} ^! (\sp _{Y ^{\dag}  \hookrightarrow U ^{\dag},T _0 , + }(E))$.
  Via le morphisme $v ' _+ v ^{\prime !} \rightarrow id $, il en résulte
  $\sp _{Y ^{\prime \dag}  \hookrightarrow U ^{\prime \dag},T '_0 , + }(b ^*E) [d_{X '_0/X _0}]\rightarrow
  f _{T _0} ^! (\sp _{Y ^{\dag}  \hookrightarrow U ^{\dag},T _0 , + }(E))$.
  En lui appliquant $\R \underline{\Gamma} ^\dag _{X ' _0}$, comme
  $\sp _{Y ^{\prime \dag}  \hookrightarrow U ^{\prime \dag},T '_0 , + }(b ^*E) \in
  \mathrm{Coh} ( \PP', T' _0 ,X' _0)$, on obtient la flèche
  $\sp _{Y ^{\prime \dag}  \hookrightarrow U ^{\prime \dag},T '_0 , + } (b ^*E)[d_{X '_0/X _0}]
\rightarrow
\R \underline{\Gamma} ^\dag _{X ' _0} f _{T _0} ^!
(\sp _{Y ^{\dag}  \hookrightarrow U ^{\dag},T _0 , + }(E) ).$
Comme celle-ci est un isomorphisme en dehors de $T ' _0$, on conclut avec
\cite[4.3.12]{Be1}.
\hfill \hfill \qed \end{proof}

\section{$F$-isocristaux surcohérents sur les schémas lisses}
\subsection{Surcohérence différentielle des $F$-isocristaux surconvergents dans le cas d'une compactification lisse}
\begin{vide}[Propriété $P _{\PP,T}$]
\label{PPPT}
  Soient $\PP$ un $\V$-schéma formel lisse, $T$ un diviseur de $P$ et
  $\E \in D^\mathrm{b} _\mathrm{coh} (\D ^\dag _{\PP} (\hdag T  ) _\Q )$.
  On dira que $\E$ vérifie $P _{\PP,T}$ si $\E$ est $\D ^\dag _{\PP} (\hdag T  ) _\Q$-surcohérent
  et si pour tout morphisme
  lisse $\alpha$ : $\QQ \rightarrow \PP$, pour tout sous-schéma fermé $Z  \subset Q $, en notant
  $U := \alpha ^{-1} (T)$, $\DD _{\QQ, U} \R \underline{\Gamma} ^\dag _{Z} (\alpha _T ^! (\E))$ est
  $\D ^\dag _{\QQ} (\hdag U) _\Q$-surcohérent.
\end{vide}

\begin{rema}\label{remaPPPT}
 (i)  Avec les notations de \ref{PPPT}, comme la $\D ^\dag _{\QQ} (\hdag U) _\Q$-surcohérence est
 locale en $\QQ$, la propriété $P _{\PP,T}$ est locale en $\PP$, i.e., si $(\PP _\alpha) _{\alpha \in \Lambda}$
 est un recouvrement ouvert de $\PP$, alors $\E $ vérifie $P _{\PP,T}$
 si et seulement si $\E |_{\PP _\alpha}$ vérifie $P _{\PP _\alpha,T\cap P _\alpha}$ pour tout $\alpha \in \Lambda$.

  (ii)  Il n'est pas évident que la propriété $P _{\PP,T}$ se préserve par foncteur dual.

  (iii) $\E$ vérifie $P _{\PP,T}$ si et seulement si
  $\E$ est $\D ^\dag _{\PP} (\hdag T  ) _\Q$-surcohérent
  et si pour tout morphisme
  $\alpha$ : $\QQ \rightarrow \PP$ tel que $U :=\alpha ^{-1} (T)$ soit un diviseur de $Q$,
  pour tous sous-schémas fermés $Z,Z'  \subset Q $,
  $\DD _{\QQ, U} \R \underline{\Gamma} ^\dag _{Z}(\hdag Z') (\alpha _T ^! (\E))$ est
  $\D ^\dag _{\QQ} (\hdag U) _\Q$-surcohérent.

  En effet, la condition est suffisante. Prouvons qu'elle est nécessaire.
  Via un triangle de localisation, on peut supposer $Z' =\emptyset$.
  De plus, $\alpha $ est le composé de son graphe (que l'on suppose fermé pour simplifier les notations)
  $\gamma$ : $\QQ \rightarrow  \QQ \times \PP$ suivi
  de la projection canonique $p$ : $\QQ \times \PP \rightarrow \PP$.
  Or, grâce au théorème de Kashiwara, au théorème de dualité relative et à la commutation
  de la cohomologie locale à l'image directe, en posant $T' = p ^{-1} (T)$, on a
 \begin{gather}\notag
   \DD _{\QQ, U} \R \underline{\Gamma} ^\dag _{Z} (\alpha _T ^! (\E)) \riso
  \DD _{\QQ, U} \R \underline{\Gamma} ^\dag _{Z} \gamma _{T' } ^! p _T ^! (\E)
  \riso
  \gamma _{T' } ^! \gamma _{T' +} \DD _{\QQ, U} \R \underline{\Gamma} ^\dag _{Z} \gamma _{T' } ^! p _T ^! (\E)
  \\
  \riso
  \gamma _{T' } ^!  \DD _{\QQ, T'} \gamma _{T' +} \R \underline{\Gamma} ^\dag _{Z} \gamma _{T' } ^! p _T ^! (\E)
  \notag
  \riso
\gamma _{T' } ^!  \DD _{\QQ, T'}  \R \underline{\Gamma} ^\dag _{Z} \gamma _{T' +} \gamma _{T' } ^! p _T ^! (\E)
\\ \notag
\riso
\gamma _{T' } ^!  \DD _{\QQ, T'}  \R \underline{\Gamma} ^\dag _{Z}  \R \underline{\Gamma} ^\dag _{Q} p _T ^! (\E)
\riso
\gamma _{T' } ^!  \DD _{\QQ, T'}  \R \underline{\Gamma} ^\dag _{Z} p _T ^! (\E).
 \end{gather}
 On conclut la remarque grâce à la stabilité de la surcohérence par image inverse extraordinaire.
\end{rema}

\begin{prop}\label{PPPTstable+propre}
  Soient $f$ : $ \PP ' \rightarrow \PP$ un morphisme de $\V$-schémas formels lisses,
  $T$ un diviseur de $P$ tel que $T' :=f ^{-1} (T)$ soit un diviseur de $Q$.
  Pour tous sous-schémas fermés $Z'_1$ et $Z ' _2$ de $P'$, pour tout
  $\E \in D^\mathrm{b} _\mathrm{coh} (\D ^\dag _{\PP} (\hdag T  ) _\Q )$ vérifiant $P _{\PP ,T}$,
  $\R \underline{\Gamma} ^\dag _{Z'_2} (\hdag Z ' _1) f _T ^! (\E)$ satisfait $P _{\PP ' ,T'}$.
  Si $f$ est propre alors, pour tout
  $\FF \in D^\mathrm{b} _\mathrm{coh} (\D ^\dag _{\PP '} (\hdag T'  ) _\Q )$
  vérifiant $P _{\PP ' ,T'}$, $f _{T,+} (\E )$ satisfait $P _{\PP,T}$.
\end{prop}
\begin{proof}
La première assertion résulte de \ref{remaPPPT}.(iii) et
de la commutation de l'image inverse extraordinaire à la cohomologie locale
(voir \cite[2.2.18.1]{caro_surcoherent}) et à la composition.
Prouvons à présent la deuxième.
  Soient $\alpha$ : $\QQ\rightarrow \PP$ un morphisme lisse et $Z $
  un sous-schéma fermé de $Q$. On note $\QQ ' := \PP ' \times _{\PP} \QQ$,
  $\alpha '$ : $\QQ ' \rightarrow \PP '$ et $f'$ : $ \QQ ' \rightarrow \QQ$ les projections, $U := \alpha ^{-1} (T)$,
  $U' := \alpha ^{\prime -1} (T')$ et $Z' : =f ^{\prime -1} (Z)$.
  On obtient :
  $\DD _{\QQ, U} \R \underline{\Gamma} ^\dag _{Z} \alpha _T ^! (f _{T,+} (\FF ))
  \riso
  f '_{T,+} \DD _{\QQ ', U'} \R \underline{\Gamma} ^\dag _{Z'}  \alpha _T ^{\prime !} (\FF )$
grâce à \cite[1.2.5]{caro_courbe} et \cite[2.2.18, 3.1.8]{caro_surcoherent}.
Par hypothèse, $\DD _{\QQ ', U'} \R \underline{\Gamma} ^\dag _{Z'}  \alpha _T ^{\prime !} (\FF )$ est
$\D ^\dag _{\QQ '} (\hdag U'  ) _\Q$-surcohérent.
  De même que \cite[1.1.15]{caro_surholonome},
 on remarque que \cite[3.1.9]{caro_surcoherent} est valable lorsque le morphisme $f _0$ se relève.
Il en dérive la $\D ^\dag _{\QQ} (\hdag U  ) _\Q$-surcohérence de
$f '_{T,+} \DD _{\QQ ', U'} \R \underline{\Gamma} ^\dag _{Z'}  \alpha _T ^{\prime !} (\FF )$.
\hfill \hfill \qed \end{proof}

\begin{theo}
\label{theo-surcohlisse}
  Soient $\PP$ un $\V$-schéma formel séparé et lisse, $X$ un sous-schéma fermé lisse de $P$ et $T$ un diviseur de $P$
  tel que $T _X := T \cap X$ soit un diviseur de $X$. On note $\U$ (resp. $Y$) l'ouvert de $\PP$ (resp. $X$)
  complémentaire de $T$ (resp. $T _X$).

  Pour tout isocristal $E$ sur $Y$ surconvergent le long de $T _X $, le $\D ^\dag _{\PP} (\hdag T  ) _\Q$-module
  cohérent à support dans $X$ associé, $\sp _{X \hookrightarrow \PP, T,+} (E)$ (\cite[1.5]{caro_unite}),
  vérifie $P _{\PP,T}$.
\end{theo}

\begin{proof}
Notons $\E :=\sp _{X \hookrightarrow \PP, T,+} (E)$.
  Comme le théorème est local en $\PP$, on peut supposer $P$ affine et irréductible. L'immersion fermée
  $X \hookrightarrow P$ se relève en une immersion fermée de $\V$-schémas formels lisses.
  Grâce à \ref{PPPTstable+propre},
  on se ramène donc à traiter le cas où $X=P$. On notera alors $T $ pour $T _X$ et $\X$ pour $\PP$.

  Tout d'abord, il est suffisant de prouver que pour tout sous-schéma fermé $Z$ de $X$,
  $\R \underline{\Gamma} ^\dag _Z (\E) $ est $\D ^\dag _{\X} (\hdag T ) _\Q$-cohérent (ce qui implique la surcohérence
  de $\E$) puis que
  $\DD _{\X , T} \R \underline{\Gamma} ^\dag _{Z} (\E)$ est
  $\D ^\dag _{\X} (\hdag T ) _\Q$-surcohérent.
En effet, l'image inverse extraordinaire d'un isocristal surconvergent est, à un décalage près,
un isocristal surconvergent. \'Etablissons donc ces deux propriétés.

  Grâce au théorème de désingularisation de de Jong (\cite{dejong}), il existe un morphisme
  projectif, surjectif, génériquement fini et étale $a$ : $X '\rightarrow X$
  tel que $X'$ soit irréductible et $k$-lisse,
  et tel que $a ^{-1} (Z)$ soit un diviseur à croisements normaux de $X'$.
  Comme $a$ est projectif, il existe (\cite{Beintro2})
  donc un $\V$-schéma formel lisse $\PP '$, une immersion fermée $u '$ : $ X' \hookrightarrow \PP'$, un morphisme
  propre et lisse $f$ : $\PP' \rightarrow \X$ tels que $f \circ u' = u \circ a$.
  On notera $T'=f ^{-1} (T)$ et $Y':= a ^{-1} (Y)$.

  Par adjonction (\cite[1.2.6]{caro_courbe}), on a  :
  $f _{T+} \R \underline{\Gamma} ^\dag _{X'} f _T ^! (\E ) \rightarrow \E $.
  En dualisant celui-ci et en utilisant les théorèmes de dualité relative et de bidualité (\cite{virrion}),
  on obtient le suivant :
  $\E \rightarrow f _{T+} \DD _{\PP ', T'} \R \underline{\Gamma} ^\dag _{X'} f _T ^! \DD _{\X, T }( \E )$.

  Or, comme $X'$ est irréductible et $T'$ est un diviseur de $P'$,
  $X '\cap T'$ est alors un diviseur de $X'$ ou est égal à $X'$. Puisque
  $a$ est génériquement fini et étale, on ne peut avoir $X '\cap T'=X'$. On en conclut que
  $X '\cap T'$ est un diviseur de $X'$. D'où les isomorphismes :
  $\R \underline{\Gamma} ^\dag _{X'} f_T ^! (\E )\riso \sp _{X' \hookrightarrow \PP', T',+}(a ^*(E)),$
  $\DD _{\PP ', T'} \R \underline{\Gamma} ^\dag _{X'} f_T ^! \DD _{\X, T }( \E ) \riso
  \sp _{X' \hookrightarrow \PP', T'+} (a ^*(E ^\vee )^\vee) \riso \sp _{X' \hookrightarrow \PP', T'+} (a ^*(E))$
  (\cite{caro_unite}).
  Les morphismes d'adjonction induisent alors la suite de morphismes :
   $ \E \rightarrow f _{T+} \R \underline{\Gamma} ^\dag _{X'} f_T ^! (\E ) \rightarrow \E.$
On vérifie de manière analogue à \cite{caro_unite}, que ce morphisme composé est un isomorphisme, i.e.,
que $\E$ est un facteur direct de $f _{T+} \R \underline{\Gamma} ^\dag _{X'} f_T ^! (\E )$.
Par \ref{PPPTstable+propre}, on se ramène alors au cas où $Z$ est un diviseur à croisements normaux de $X$.

Grâce aux suites spectrales de Mayer-Vietoris (\cite{caro_surcoherent}),
il nous reste à traiter le cas où $Z$ est un sous-schéma fermé
intègre et lisse de $X$. On a alors deux cas : soit $Z \subset T$, soit
$Z \cap T$ est un diviseur de $Z$. Le premier cas donne les égalités
$\R \underline{\Gamma} ^\dag _Z (\E) =0$ et
$\DD _{\X , T} \R \underline{\Gamma} ^\dag _{Z} (\E)=0$. 
Traitons à présent le second cas. Les deux propriétés à vérifier étant locales en $\X$,
il ne coûte rien de supposer que l'immersion fermée $Z \hookrightarrow X$ se
relève en une immersion fermée $u$ : $\ZZ\hookrightarrow \X$ de $\V$-schémas formels lisses.
L'isomorphisme $\R \underline{\Gamma} ^\dag _Z (\E) \riso u _+ u ^! (\E) $ et le fait que $u ^! (\E)$ soit,
au décalage $[d_{Z/X}]$ près,
un isocristal sur $Z \setminus T$ surconvergent le long de $Z \cap T$ nous permet de conclure que
$\R \underline{\Gamma} ^\dag _Z (\E)$ est $\D ^\dag _{\X} (\hdag T ) _\Q$-cohérent.
Enfin, comme $\DD _{\X , T} \R \underline{\Gamma} ^\dag _{Z} (\E) \riso u _+ \DD _{\ZZ, T} u^!(\E)$ et puisque
$\DD _{\ZZ, T} u^!(\E)$ est, au décalage $[d_{Z/X}]$ près,
un isocristal sur $Z \setminus T$ surconvergent le long de $Z \cap T$,
on conclut que $\DD _{\X , T} \R \underline{\Gamma} ^\dag _{Z} (\E)$
est $\D ^\dag _{\X} (\hdag T ) _\Q$-surcohérent.
\hfill \hfill \qed \end{proof}

\begin{prop}\label{lissestableotimes}
  Avec les notations de \ref{theo-surcohlisse},
  pour tous isocristaux $E_1$ et $E _2$ sur $Y$ surconvergents le long de $T _X $,
  en notant $\sp _+ = \sp _{X \hookrightarrow \PP, T,+}$, on dispose
  d'un isomorphisme canonique :
  $$\sp _+ (E _1 \otimes _{j ^\dag \O _{\PP _K}} E _2) \riso
  \sp _+ (E _1 )
  \smash{\overset{\L}{\otimes}} ^\dag _{\O _{\PP} (\hdag T) _\Q}
  \sp _+ (E _2 )[d _{X/P}] .$$
\end{prop}
\begin{proof}
  Supposons d'abord $\PP$ affine. Il existe alors un relèvement
  $u$ : $\X \hookrightarrow \PP$ de $X \hookrightarrow P$.
  En notant $\E _1 := \sp _* (E _1)$ et $\E _2 := \sp _* (E _2)$, on obtient
  $\sp _+ (E _1 ) \riso u _{T+} (\E _1)$,
  $\sp _+ (E _2 ) \riso u _{T+} (\E _2)$  et
  $\sp _+ (E _1 \otimes E _2) \linebreak \riso
  u _{T+}  ( \E _1 \otimes _{\O _{\X} ( \hdag T _X) _\Q} \E _2)$, ces isomorphismes
  ne dépendant pas (modulo les isomorphismes induits par \cite[2.1.5]{Be2})
  du relèvement $u$ choisi.

  Or, $\E _2 \riso u ^! _T u_{T +} (\E _2)$ et,
  de manière analogue à \cite[2.1.4]{caro_surcoherent},
  on dispose d'un isomorphisme canonique
 $$ u _{T+}  ( \E _1 \otimes _{\O _{\X} ( \hdag T _X) _\Q} u ^! _T (u_{T +} (\E _2)))
 \riso
 u _{T+}  ( \E _1)
 \smash{\overset{\L}{\otimes}} ^\dag _{\O _{\PP} (\hdag T) _\Q}
 u_{T +} (\E _2))[d _{X/P}]$$
  qui ne dépend pas du relèvement $u$ choisi.
  Il en dérive l'isomorphisme de \ref{lissestableotimes}.

  Lorsque $\PP$ est quelconque, via la construction précédente, on obtient
  localement en $\PP$ des isomorphismes qui se recollent.
\hfill \hfill \qed \end{proof}
\subsection{Construction des $F$-isocristaux surcohérents}

Sauf mention contraire, on gardera les notations suivantes :
soient $\PP$ un $\V$-schéma formel séparé et lisse, $u$ : $X \hookrightarrow P$ une immersion fermée, $T$ un diviseur de $P$
  et $\U$ l'ouvert de $\PP$ complémentaire de $T$.
  Enfin, on suppose $Y :=X \setminus T$ lisse.

\begin{defi}\label{defidagdagpxt}

  On définit la catégorie ($F$-)$\mathrm{Isoc} ^{\dag \dag}( \PP, T, X/K)$ de la manière suivante :
  les objets sont les
  ($F$-)$\D ^\dag _{\PP} (\hdag T) _\Q$-modules cohérents $\E$ à support dans $X$ tels que
  \begin{enumerate}
    \item $\E |_{\U} \in \mathrm{Isoc} ^{\dag}( \U ,Y /K)$ (notations de \ref{plfidked}), i.e.,
    il existe un isocristal $G$ convergent sur $Y$ tel que
    $\E |_{\U} \riso \sp _{Y \hookrightarrow \U,+}(G)$;
   \item $\E$ et $\DD_{\PP,T}  (\E)$ sont $\D ^\dag _{\PP} (\hdag T) _\Q$-surcohérents.
  \end{enumerate}
  Enfin, les flèches sont les morphismes ($F$-)$\D ^\dag _{\PP} (\hdag T) _\Q$-linéaires.
  Lorsque le diviseur $T$ est vide, on omettra de l'indiquer.

  Comme le foncteur $\sp _{Y \hookrightarrow \U,+}$ est pleinement fidèle, l'isocristal $G$
  est défini à isomorphisme près. Avec les notations de \cite{caro_surholonome},
  on a en fait \linebreak $G \riso  \mathcal{R}ecol \circ \sp ^* \circ \mathcal{L} oc ( \E |_{\U})$.
  De plus, lorsque $\E \in F \text{-}\mathrm{Isoc} ^{\dag \dag}( \PP, T, X/K)$, on a
  automatiquement $\E  |_{\U} \in F \text{-} \mathrm{Isoc} ^{\dag}( \U ,Y /K)$.
\end{defi}

\begin{vide}\label{dag2=daglisse}
Il découle du théorème \ref{theo-surcohlisse} que
l'égalité
$F \text{-}\mathrm{Isoc} ^{\dag }( \PP, T, X/K) =
F \text{-}\mathrm{Isoc} ^{\dag \dag }( \PP, T, X/K)$
est validée si $X$ est lisse.
Le foncteur $\sp _{X \hookrightarrow \PP,T,+}$ induit alors
l'équivalence de catégorie : \newline
$F \text{-}\mathrm{Isoc} ^{\dag }( Y, X/K) \cong
F \text{-}\mathrm{Isoc} ^{\dag }( \PP, T, X/K)
=F \text{-}\mathrm{Isoc} ^{\dag \dag }( \PP, T, X/K)$.

Via la proposition ci-après,
on se ramène par descente au cas où $X$ est lisse.
\end{vide}

\begin{prop}\label{defidagdagpxt?}
On suppose $X$ irréductible et $T \not \supset X$.
  Il existe alors un diagramme commutatif de la forme
  $$\xymatrix @R=0,3cm {
  { Y '} \ar[r] \ar[d] ^b & {X'} \ar[r] ^{u'} \ar[d] ^a & {\PP'} \ar[d] ^f \\
  {Y} \ar[r] & {X } \ar[r] ^{u} & {\PP,}}$$
  où $f$ est un morphisme propre et lisse de $\V$-schémas formels lisses,
  le carré de gauche est cartésien, $X'$ est lisse, $u'$ est une immersion fermée,
  et $a$ est un morphisme
  projectif, génériquement fini et étale, tel que
  $a ^{-1} (T \cap X)$ soit un diviseur à croisement normaux de $X'$.

  Supposons choisi un tel diagramme.
  Soit $\E\in F\text{-}\mathrm{Isoc} ^{\dag \dag}( \PP, T, X/K)$.
 Il existe un
  (unique à isomorphisme près) $F$-isocristal $E'$ sur $Y'$ surconvergent le long de $a ^{-1} (T \cap X)$ tel que
  $\R \underline{\Gamma} ^\dag _{X'} f _T ^! (\E ) \riso \sp _{X' \hookrightarrow \PP', f ^{-1} (T)} (E')$.
  De plus, le faisceau $\E$ est un facteur direct de
  $f _{T,+} \sp _{X' \hookrightarrow \PP', f ^{-1} (T)} (E')$.
\end{prop}

\begin{proof}
  Grâce au théorème de désingularisation de de Jong (\cite{dejong}), il existe un morphisme
  projectif, surjectif, génériquement fini et étale $a$ : $X '\rightarrow X$ tel que $X'$ soit irréductible et $k$-lisse,
  et $a ^{-1} (T\cap X)$ soit un diviseur à croisements normaux de $X'$.
  Comme $a$ est projectif, il existe
  un $\V$-schéma formel lisse $\PP '$, une immersion fermée $u '$ : $ X' \hookrightarrow \PP'$, un morphisme
  propre et lisse $f$ : $\PP' \rightarrow \PP$ tels que $f \circ u' = u \circ a$.
  En posant $ Y ' :=X ' \setminus a ^{-1} (T \cap X) $,
  on obtient l'existence du diagramme de la proposition.

  On note $\U '$ l'ouvert de
  $\PP'$ complémentaire de $f ^{-1} (T)$ et $g$ : $ \U' \rightarrow \U$ le morphisme induit par $f$.
  Par adjonction $f _{T,+} \R \underline{\Gamma} ^\dag _{X'} f_T ^! (\E ) \rightarrow \E$.
  Il en dérive par dualité
  $\E \rightarrow f _{T,+} \DD ^* _{\PP ', f^{-1} (T) } \R \underline{\Gamma} ^\dag _{X'} f_T ^! \DD ^* _{\PP, T}(\E )$.
  Comme le foncteur $\sp _+ ^\prime :=\sp _{X' \hookrightarrow \PP',f ^{-1} (T )+}$ est pleinement fidèle,
  via la caractérisation de son image essentielle,
  il existe de manière unique des $F$-isocristaux $E'_1$, $E _2 '$ sur $ Y '$ surconvergents
  le long de $a ^{\text{-}1} (T \cap X)$ tels que
$\R \underline{\Gamma} ^\dag _{X'} f_T ^! (\E ) \riso \sp _+ ^{\prime}(E'_1) $ et
 $ \DD ^* _{\PP ', f^{-1} (T) } \R \underline{\Gamma} ^\dag _{X'} f_T ^! \DD^* _{\PP, T}(\E )
  \riso \sp _+ ^{\prime}(E _2 ')$. 
 En notant $\widehat{E} ' _1$, $\widehat{E} ' _2$ leur image par 
$F\text{-}\mathrm{Isoc} ^{\dag}(Y',X'/K)\rightarrow F\text{-}\mathrm{Isoc} (Y',X'/K)$,
on obtient les isomorphismes
$\smash{\widehat{\sp}} ' _+ (\widehat{E} ' _1) \riso
\R \underline{\Gamma} ^\dag _{Y'} g ^! (\E | _{\U})$
et
$\smash{\widehat{\sp}} ' _+ (\widehat{E} ' _2) \riso
\DD _{\U '}^*  \R \underline{\Gamma} ^\dag _{Y'}  g ^! \DD ^* _{\U} (\E | _{\U})$,
où $\smash{\widehat{\sp}} ' _+ :=\sp _{Y' \hookrightarrow \U '+}$.

Or, $\E |_{\U} \riso \sp _{Y \hookrightarrow \U,+}(G)$, où
$G$ est un $F$-isocristal convergent sur $Y$.
D'où
$\R \underline{\Gamma} ^\dag _{Y'} g ^! (\E | _{\U}) \riso
\smash{\widehat{\sp}} ' _+ b ^* (G)$
et
$\DD _{\U '}^*  \R \underline{\Gamma} ^\dag _{Y'}  g ^! \DD ^* _{\U} (\E | _{\U}) \riso
\smash{\widehat{\sp}} ' _+ (b ^* (G ^\vee ) ^\vee)$.
Puisque $b ^* (G) \riso b ^* (G ^\vee ) ^\vee$
et comme $\smash{\widehat{\sp}} ' _+$ est pleinement fidèle,
$\widehat{E} ' _1 \riso \widehat{E} ' _2$.
  Par pleine fidélité de $F\text{-}\mathrm{Isoc} ^{\dag}(Y',X'/K)\rightarrow F\text{-}\mathrm{Isoc} (Y',X'/K)$
  (\cite{kedlaya_full_faithfull}), on en tire $E ' _1 \riso E' _2$.
  D'où :
  $\E \rightarrow f _{T,+} \R \underline{\Gamma} ^\dag _{X'} f _T ^! (\E ) \rightarrow \E$.
  Or, ce morphisme composé est un isomorphisme puisqu'il l'est en dehors de $T$. Ainsi,
  $\E$ est un facteur direct de $f _{T,+} \R \underline{\Gamma} ^\dag _{X'} f_T ^! (\E )$.
\hfill \hfill \qed \end{proof}

On remarque que la proposition \ref{defidagdagpxt?}, très utile par la suite,
nécessite une structure de Frobenius (en effet, on a utilisé le théorème de pleine fidélité
de Kedlaya \cite{kedlaya_full_faithfull} dont la généralisation sans structure de Frobenius est
à l'heure actuelle une conjecture).
Intéressons-nous maintenant à la stabilité de
$F\text{-} \mathrm{Isoc} ^{\dag \dag}( \PP, T, X/K)$.

\begin{prop}\label{daddagstDjdag}
    La catégorie $F\text{-}\mathrm{Isoc} ^{\dag \dag}( \PP, T, X/K)$ est stable par
  $\DD_{\PP,T}$.
  Soit $T '\supset T$ un deuxième diviseur de $P$. Le foncteur localisation
$(\hdag T')$ induit
  $(\hdag T')\ : \ F\text{-} \mathrm{Isoc} ^{\dag \dag}( \PP, T, X/K)
  \rightarrow F\text{-} \mathrm{Isoc} ^{\dag \dag}( \PP, T', X/K)$.
\end{prop}
\begin{proof}
  La première assertion découle de l'isomorphisme canonique
$\DD _{\U} (\E |_{\U}) \riso \sp _{Y \hookrightarrow \U,+}(G ^\vee )$
  (\cite{caro_surholonome}). La deuxième résulte de la commutation
  du foncteur dual à l'extension des scalaires $(\hdag T')$ et de la stabilité de la surcohérence par $(\hdag T')$.
\hfill \hfill \qed \end{proof}

\begin{vide}\label{isosurcohpf0}
 Comme $Y$ est lisse,
  $Y$ est la somme directe de ses composantes connexes $Y _r$, pour $r=1,\dots ,N$.
  Comme, pour tout $r$, les immersions $Y _r \hookrightarrow Y$ et $Y \hookrightarrow P \setminus T$
sont fermées,
  on a $\overline{Y} \setminus T =Y$ et $\overline{Y} _r \setminus T = Y _r$, où
  $\overline{Y}$ (resp. $\overline{Y} _r$) est l'adhérence schématique de $Y$ (resp. $Y _r$) dans $P$.
  Pour tout objet $\E$ de $F\text{-}\mathrm{Isoc} ^{\dag \dag}( \PP, T , \overline{Y}/K)$,
  le morphisme canonique
  $\oplus _r \R \underline{\Gamma} ^\dag  _{\overline{Y} _r} \E \rightarrow \E$
  est alors un isomorphisme en dehors de $T$. Par \cite[4.3.12]{Be1}, ce dernier est un isomorphisme.
  Ainsi, on obtient des foncteurs $\R \underline{\Gamma} ^\dag  _{\overline{Y} _r}$ :
  $F\text{-}\mathrm{Isoc} ^{\dag \dag}( \PP, T , \overline{Y}/K) \rightarrow F\text{-}\mathrm{Isoc} ^{\dag \dag}( \PP, T , \overline{Y} _r/K)$.
  De plus, il en dérive une équivalence de catégorie
  $\oplus _r \R \underline{\Gamma} ^\dag  _{\overline{Y} _r}$ :
  $F\text{-}\mathrm{Isoc} ^{\dag \dag}( \PP, T , \overline{Y}/K) \cong
  \oplus _r F\text{-}\mathrm{Isoc} ^{\dag \dag}( \PP, T , \overline{Y} _r/K)$.
\end{vide}

\begin{prop}\label{propindX}
  Soient $\PP$ un $\V$-schéma formel séparé et lisse, $X$, $X'$ deux sous-schémas fermés de $P$, $T $ et $T'$ deux
  diviseurs de $P$ tels que $Y := X \setminus T$ soit lisse et $X \setminus  T =X'\setminus T'$.
  On a alors les égalités $F\text{-}\mathrm{Isoc} ^{\dag \dag} ( \PP, T, X/K)=
  F\text{-}\mathrm{Isoc} ^{\dag \dag}  (\PP, T', X'/K)$.
  En particulier, en notant $\overline{Y}$ l'adhérence schématique de $Y$ dans $P$,
  $F\text{-}\mathrm{Isoc} ^{\dag \dag} ( \PP, T, X/K)=
  F\text{-}\mathrm{Isoc} ^{\dag \dag} (\PP, T, \overline{Y} /K)$.

  En outre, pour tout $\E \in F\text{-}\mathrm{Isoc} ^{\dag \dag} ( \PP, T, X/K)$,
  $\E$ vérifie $P _{\PP ,T}$ (voir \ref{PPPT}).
\end{prop}
    \begin{proof}
    Le cas où $X \setminus T$ est vide implique, en notant $0$ le faisceau nul,
    $F\text{-}\mathrm{Isoc} ^{\dag \dag} (\PP, T, X/K)=F\text{-}\mathrm{Isoc} ^{\dag \dag} ( \PP, T', X'/K)=\{0 \}$.
    Le cas où $T'=T$ se traite de manière analogue à \ref{propindXlisse}.

    Comme $X \setminus T = \overline{Y} \setminus T$ et $X '\setminus T '= \overline{Y} \setminus T'$,
    on se ramène ainsi au cas où $X = X' =\overline{Y}$ et $X \setminus T$ est non vide.
    Par \ref{isosurcohpf0},
    on peut en outre supposer $X$ irréductible.

    Soit $\E$ un objet de $F\text{-}\mathrm{Isoc} ^{\dag \dag} ( \PP, T, X/K)$.
    Par symétrie, pour terminer la preuve, il suffit de prouver que $\E$ vérifie $P _{\PP ,T'}$.

  Mais, comme $X \setminus T =X \setminus T '$, il
  découle de \ref{propindXlisse} et \ref{defidagdagpxt?} (et avec ses notations) que
  $\R \underline{\Gamma} ^\dag _{X'} f_T ^! (\E ) $ est associé
  via $\sp _{X' \hookrightarrow \PP',f ^{-1} (T ')+}$ à un $F$-isocristal
   sur $ X ' \setminus a ^{-1} (T '\cap X) $ surconvergent le long de $a ^{-1} (T '\cap X)$.
   Grâce à \ref{theo-surcohlisse}, $\R \underline{\Gamma} ^\dag _{X'} f _T^! (\E ) $
   vérifie $P _{\PP',f ^{-1} T '}$.
   Par \ref{PPPTstable+propre}, $f _{T,+} \R \underline{\Gamma} ^\dag _{X'} f _T ^! (\E )$ vérifie
   $P _{\PP ,T'}$. Comme $\E$ est un facteur direct de $f _{T,+} \R \underline{\Gamma} ^\dag _{X'} f_T ^! (\E )$,
   il en est de même de $\E$.

\hfill \hfill \qed \end{proof}

Grâce à \ref{propindX}, la catégorie $F\text{-}\mathrm{Isoc} ^{\dag \dag}( \PP , T , Z /K)$
est stable par images inverses extraordinaires :

\begin{prop}\label{isoctilstaiinvlis}
   Soient $f$ : $\PP ' \rightarrow \PP $ un morphisme de $\V$-schémas formels séparés et lisses,
   $X $ et $X '$ deux sous-schémas fermés respectifs de $P$ et $P'$, $T $ un diviseur de $P $ tel que
   $T ':= f ^{-1} (T )$ soit un diviseur.
   On suppose les $k$-schémas $Y :=X \setminus T $ et $Y ' : = X ' \setminus T'$
   lisses et $f (Y') \subset Y$.

   On dispose de la factorisation
   $$\R \underline{\Gamma} ^\dag _{X '} f _{T } ^! [-d _{Y '/Y}]\ : \
   F\text{-}\mathrm{Isoc} ^{\dag \dag}( \PP , T , X /K)
   \rightarrow
   F\text{-}\mathrm{Isoc} ^{\dag \dag}( \PP ', T ', X '/K).$$
\end{prop}
\begin{proof}
Soient $\E \in F\text{-}\mathrm{Isoc} ^{\dag \dag}( \PP, T, X/K)$, $\U :=\PP \setminus T$,
$\U ' := \PP' \setminus T'$, $g$ : $\U' \rightarrow \U$ et
$b$ : $Y '\rightarrow Y$ les morphismes induits par $f$.
Par \ref{propindX} et \ref{PPPTstable+propre},
$\R \underline{\Gamma} ^\dag _{Z '} f _{T } ^! (\E )[-d _{Y '/Y}]$ vérifie la propriété $P _{\PP', T'}$
et donc \ref{defidagdagpxt}.2.
De plus, $\R \underline{\Gamma} ^\dag _{X'} f _T ^! [-d_{Y'/Y}] (\E) |_{\U'} \riso
\R \underline{\Gamma} ^\dag _{Y'} g ^! [-d_{Y'/Y}] (\E |_{\U})$.
D'après \cite[2.2.27]{caro_surholonome}, on obtient,
$\R \underline{\Gamma} ^\dag _{Y'} g ^! [-d_{Y'/Y}] (\E |_{\U}) \riso
\sp _{Y '\hookrightarrow \U '} (b ^*(G))$, où
$G$ est le $F$-isocristal convergent sur $Y$ tel que
$\E |_{\U} \riso \sp _{Y \hookrightarrow \U} (G)$.
  D'où le résultat.
\hfill \hfill \qed \end{proof}

\begin{prop}\label{dagdagstblotimes}
Soient $\E _1$ et
$\E _2 \in F \text{-}\mathrm{Isoc} ^{\dag \dag}( \PP, T, X/K)$.

On a $\E _1  \smash{\overset{\L}{\otimes}} ^\dag _{\O _{\PP} (\hdag T) _\Q} \E _2[d _{X/P}]
\in F \text{-}\mathrm{Isoc} ^{\dag \dag}( \PP, T, X/K)$.
\end{prop}
\begin{proof}
Supposons $X = \overline{Y}$. Avec les notations de \ref{isosurcohpf0}, comme
$\R \underline{\Gamma} ^\dag  _{\overline{Y} _r}
(\E _1  \smash{\overset{\L}{\otimes}} ^\dag _{\O _{\PP} (\hdag T) _\Q} \E _2)
\riso
\R \underline{\Gamma} ^\dag  _{\overline{Y} _r} (\E _1 )
\smash{\overset{\L}{\otimes}} ^\dag _{\O _{\PP} (\hdag T) _\Q}
\R \underline{\Gamma} ^\dag  _{\overline{Y} _r} (\E _2 )$,
on se ramène au cas où $X$ est irréductible.
Il découle de \ref{lissestableotimes} que
$\E _1  \smash{\overset{\L}{\otimes}} ^\dag _{\O _{\PP} (\hdag T) _\Q} \E _2[d _{X/P}] |_{\U}
= \E _1  | _{\U} \smash{\overset{\L}{\otimes}} ^\dag _{\O _{\U,\Q}} \E _2 |_{\U} [d _{Y/U}]
\in F\text{-}\mathrm{Isoc} ^\dag ( \U, Y/K)$.
De plus, grâce à \ref{defidagdagpxt?} et avec ses notations,
  $\E _1$ est un facteur direct de
  $f _{T+} \R \underline{\Gamma} ^\dag _{X'} f _T ^! (\E _1)$.
  Il en résulte que
$\E _1  \smash{\overset{\L}{\otimes}} ^\dag _{\O _{\PP} (\hdag T) _\Q} \E _2$
est un facteur direct de
$f _{T+} (\R \underline{\Gamma} ^\dag _{X'} f _T ^! (\E _1))
\smash{\overset{\L}{\otimes}} ^\dag _{\O _{\PP} (\hdag T) _\Q} \E _2$. 
D'après \cite[2.1.4]{caro_surcoherent} (toujours valable en rajoutant des diviseurs),
on dispose de l'isomorphisme :
$$f _{T+} ( \R \underline{\Gamma} ^\dag _{X'} f _T ^! (\E _1)
\smash{\overset{\L}{\otimes}} ^\dag _{\O '}   f _T ^! (\E _2))
\riso
f _{T+} (\R \underline{\Gamma} ^\dag _{X'} f _T ^! (\E _1) )
\smash{\overset{\L}{\otimes}} ^\dag _{\O '} \E _2  [ d _{_{P'/P}}] ,$$
où $\O ':=\O _{\PP '} (\hdag T') _\Q$.
Or, il existe des
$F$-isocristaux $E' _1$ et $E ' _2$
sur $Y'$ surconvergents le long de $a ^{-1} (T \cap X)$ tels que,
  $\R \underline{\Gamma} ^\dag _{X'} f _T ^! (\E _r) \riso \sp _+ ^\prime (E' _r)$,
  où $r =1,2$ et $\sp _+ ^\prime :=\sp _{X' \hookrightarrow \PP', f ^{-1} (T)+}$.
Avec \ref{lissestableotimes},
on obtient :
\begin{gather} \notag
\sp _+ ^\prime    (E '_1 \otimes _{j ^\dag \O _{\PP' _K}} E '_2)
\riso
\R \underline{\Gamma} ^\dag _{X'} f _T ^! (\E _1)
\smash{\overset{\L}{\otimes}} ^\dag _{\O _{\PP '} (\hdag T') _\Q}
 \R \underline{\Gamma} ^\dag _{X'}  f _T ^! (\E _2) [d _{X '/P'}]
 \\  \riso
 \notag
\R \underline{\Gamma} ^\dag _{X'} f _T ^! (\E _1)
\smash{\overset{\L}{\otimes}} ^\dag _{\O _{\PP '} (\hdag T') _\Q} f _T ^! (\E _2)[d _{X '/P'}] .
\end{gather}
Les propositions \ref{propindX} et \ref{PPPTstable+propre} nous permettent de conclure.
\hfill \hfill \qed \end{proof}

\begin{prop}\label{isosurcohpfpre}
  On reprend la construction et les notations de \ref{defidagdagpxt?}.
Le foncteur
\begin{gather}\notag
 (\R \underline{\Gamma} ^\dag _{X'} f ^!, |_\U ) \ :\
 F \text{-}\mathrm{Isoc} ^{\dag \dag}( \PP, T,X /K) \\
\rightarrow
F \text{-}\mathrm{Isoc} ^{\dag \dag}( \PP', T',X' /K) \times _{F \text{-}\mathrm{Isoc} ^\dag (\U ', Y' /K) }
F \text{-}\mathrm{Isoc}^\dag  (\U,  Y /K)\notag 
\end{gather}
est pleinement fidèle.
\end{prop}
\begin{proof}
Par \ref{isoctilstaiinvlis} et \ref{defidagdagpxt},
$\R \underline{\Gamma} ^\dag _{X'} f ^!$ :
$F \text{-}\mathrm{Isoc} ^{\dag \dag}( \PP, T,X /K)
\rightarrow \linebreak
F \text{-}\mathrm{Isoc} ^{\dag \dag}( \PP',T',X' /K)$
et
$| _\U$ : $F \text{-}\mathrm{Isoc} ^{\dag \dag}( \PP, T,X /K) \rightarrow F \text{-}\mathrm{Isoc}^\dag  (\U,  Y /K)$
sont bien définis.
De plus, la fidélité de $(\R \underline{\Gamma} ^\dag _{X'} f ^!, |_\U )$
résulte de celle de $|_\U$. Il reste à prouver que cette fidélité est pleine.

Soient $\E _1, \E _2\in F\text{-}\mathrm{Isoc} ^{\dag \dag}( \PP,T , X/K)$,
$\phi$ :
$\R \underline{\Gamma} ^\dag _{X'} f ^! (\E _1) \rightarrow \R \underline{\Gamma} ^\dag _{X'} f ^! (\E _2)$
et $\psi$ : $\E _1 |_{\U} \rightarrow \E _2 |_{\U}$ induisant le même
morphisme $\R \underline{\Gamma} ^\dag _{Y'}g ^! (\E _1 |_{\U} )
\rightarrow \R \underline{\Gamma} ^\dag _{Y'} g ^! (\E _2 |_{\U} )$.
Il existe un et un seul morphisme $\theta $ rendant commutatif le diagramme de gauche ci-dessous
\begin{equation}
\label{isosurcohpfdiag1}
  \xymatrix @C=0,5cm @R=0,4cm{
{ \E _1} \ar[r] ^-{\rho _1}  \ar@{.>}[d] ^\theta
&
{f _{T, +}\R \underline{\Gamma} ^\dag _{X'} f _{T} ^! (\E _1)}
\ar[r] ^-{\mathrm{adj}} \ar[d] ^{f _{T, +} \phi}
&
{\E _1}
\ar@{.>}[d] ^\theta
\\
{ \E _2} \ar[r] ^-{\rho _2}
&
{f _{T, +}\R \underline{\Gamma} ^\dag _{X'} f _{T} ^! (\E _2)}
\ar[r] ^-{\mathrm{adj}}
&
{\E _2,}
}
\
\xymatrix @C=0,5cm  @R=0,4cm {
{ \E _1 |_{\U}} \ar[r]^-{\rho _1 |_{\U}} \ar[d] ^{\theta |_{\U}}
&
{g _{ +}\R \underline{\Gamma} ^\dag _{Y'} g ^! (\E _1)}
\ar[d] ^{g _+ (\phi |_{\U} )}
\ar[r] ^-{\mathrm{adj}}
&
{ \E _1 |_{\U}}
\ar[d] ^{\psi}
\\
{ \E _2 |_{\U}} \ar[r] ^-{\rho _2 |_{\U}}
&
{g _{ +}\R \underline{\Gamma} ^\dag _{Y'} g  ^! (\E _2 |_{\U})}
\ar[r] ^-{\mathrm{adj}}
&
{\E _2 |_{\U},}
}
\end{equation}
où, $\mathrm{adj}$ est le morphisme induit par adjonction,
$\rho _1$ et $\rho _2$ désignent une section de $\mathrm{adj}$.
Le diagramme de gauche de \ref{isosurcohpfdiag1}
s'en déduit par restriction.
Comme $\R \underline{\Gamma} ^\dag _{Y'} g ^! \psi  = \phi |_{\U}$,
le diagramme de droite de \ref{isosurcohpfdiag1} est commutatif.
D'où $\theta |_{\U} = \psi$.

De plus, en composant les diagrammes commutatifs suivants
\begin{gather} \notag
\xymatrix @R=0,3cm {
{  f _{T} ^!}
\ar[rr] ^-{\mathrm{adj}}
&&
{ f _{T} ^!f _{T +}  f _{T} ^! }
\\
{\R \underline{\Gamma} ^\dag _{X'}  f _{T} ^!}
\ar[r] ^-{\mathrm{adj}}
\ar[u]
&
{\R \underline{\Gamma} ^\dag _{X'}  f _{T} ^!f _{T +} f _{T} ^! }
\ar[ur]
&
{\R \underline{\Gamma} ^\dag _{X'}  f _{T} ^!f _{T +}\R \underline{\Gamma} ^\dag _{X'} f _{T} ^! ,}
\ar[u]
\ar[l]
}
\\
\xymatrix @R=0,3cm{
{ f _{T} ^! f _{T +}  f _{T} ^! }
\ar[rr] ^-{\mathrm{adj}}
& &
{  f _{T} ^!}
\\
{\R \underline{\Gamma} ^\dag _{X'}  f _{T} ^! f _{T +}\R \underline{\Gamma} ^\dag _{X'} f _{T} ^!}
\ar[r] ^-{\mathrm{adj}}
\ar[u]
&
{\R \underline{\Gamma} ^\dag _{X'} \R \underline{\Gamma} ^\dag _{X'} f _{T} ^! }
\ar[r]
&
{\R \underline{\Gamma} ^\dag _{X'}  f _{T} ^!,}
\ar[u]
}\notag
\end{gather}
on vérifie que le morphisme composé du haut est l'identité. Il en est donc de même de celui du bas.
Il en résulte, en notant $f$ pour $f _T$, que les composés horizontaux de
\begin{equation}
\label{isosurcohpfdiag2}
\xymatrix @R=0,5cm {
{\R \underline{\Gamma} ^\dag _{X'}  f  ^! f _{+}\R \underline{\Gamma} ^\dag _{X'} f  ^! (\E _1)}
\ar[r] ^-{\mathrm{adj}} \ar[d] _{\R \underline{\Gamma} ^\dag _{X'}  f  ^! f _{+} \phi}
&
{\R \underline{\Gamma} ^\dag _{X'}  f ^!\E _1}
\ar[d] _{\R \underline{\Gamma} ^\dag _{X'}  f  ^! \theta}
\ar[r] ^-{\mathrm{adj}}
&
{\R \underline{\Gamma} ^\dag _{X'}  f  ^!f _{ +}\R \underline{\Gamma} ^\dag _{X'} f  ^! (\E _1)}
\ar[d] _{\R \underline{\Gamma} ^\dag _{X'}  f  ^!f _{ +}\R \underline{\Gamma} ^\dag _{X'}  f  ^! \theta}
\\
{\R \underline{\Gamma} ^\dag _{X'}  f  ^! f _{ +}\R \underline{\Gamma} ^\dag _{X'} f  ^! (\E _2)}
\ar[r] ^-{\mathrm{adj}}
&
{\R \underline{\Gamma} ^\dag _{X'}  f  ^!\E _2}
\ar[r] ^-{\mathrm{adj}}
&
{\R \underline{\Gamma} ^\dag _{X'}  f  ^! f _{ +}\R \underline{\Gamma} ^\dag _{X'} f  ^! (\E _2).}
}
\end{equation}
sont l'identité. Or, le carré de droite de \ref{isosurcohpfdiag2} est commutatif
par fonctorialité, tandis que celle du carré de gauche se vérifie en remarquant
qu'il correspond à l'image par
$\R \underline{\Gamma} ^\dag _{X'}  f  ^! \theta$
du carré de droite du diagramme de gauche de \ref{isosurcohpfdiag1}.
Ainsi, \ref{isosurcohpfdiag2} est commutatif.
D'où :
$\R \underline{\Gamma} ^\dag _{X'}  f _{T} ^! f _{T +} \phi =
\R \underline{\Gamma} ^\dag _{X'}  f _{T} ^!f _{T +}\R \underline{\Gamma} ^\dag _{X'}  f _{T} ^! \theta$.

Or,
$f _{T +}$ :
$F \text{-}\mathrm{Isoc} ^{\dag \dag}( \PP',T',X' /K)
\rightarrow F \text{-}\mathrm{Isoc} ^{\dag \dag}( \PP,T,X /K)$ est fidèle.
En effet, soit $\alpha$ : $\E ' _1 \rightarrow \E '_2$ tel que
$f _{T +} (\alpha) =0$.
D'après \cite[4.1.2]{tsumono} et \cite{kedlaya_full_faithfull}, $\alpha =0$
si et seulement si sa restriction est nulle sur un ouvert
$\PP ' _1$ de $\PP'$ tel que $P ' _1 \cap Y'$ soit dense dans $Y'$.
On se ramène ainsi à prouver que lorsque $b $ est fini et étale et $\U$ est affine,
le foncteur $g _+$ : $F\text{-}\mathrm{Isoc} ^{\dag \dag }( \U ',Y' /K)
\rightarrow
F\text{-}\mathrm{Isoc} ^{\dag \dag }( \U ,Y /K)$ est fidèle.
Cela se voit via l'isomorphisme canonique
$\sp _{Y \hookrightarrow \U,+} (b _* (E')) \riso
g _+ \sp _{Y \hookrightarrow \U,+} (E')$ valable pour
$E ' \in F \text{-}\mathrm{Isoc} (Y'/K) $.

De plus, via un diagramme de la forme de celui de gauche de \ref{isosurcohpfdiag1}, on vérifie
que $\R \underline{\Gamma} ^\dag _{X'}  f _{T} ^!$ :
$F \text{-}\mathrm{Isoc} ^{\dag \dag}( \PP,T,X /K)
\rightarrow F \text{-}\mathrm{Isoc} ^{\dag \dag}( \PP',T',X' /K)$ est fidèle.
Il en découle $\phi = \R \underline{\Gamma} ^\dag _{X'}  f _{T} ^! \theta$.

On a ainsi prouvé que le foncteur de \ref{isosurcohpfpre} est pleinement fidèle.
\hfill \hfill \qed \end{proof}

\begin{prop}\label{pffidcasgen}
 Le foncteur restriction \newline $|_\U$ :
  $F \text{-} \mathrm{Isoc} ^{\dag \dag} ( \PP, T  ,X /K) \rightarrow
  F \text{-} \mathrm{Isoc} ^{\dag \dag} ( \U, Y /K) $ est pleinement fidèle.
\end{prop}
\begin{proof}
Avec l'aide de \ref{isosurcohpf0}, on peut supposer que $ X$ est irréductible.
Soient $\E _1$ et $\E _2$ deux objets de $F\text{-}\mathrm{Isoc} ^{\dag \dag}( \PP,T , X/K)$
et $\psi $ : $\E _1 |_{\U}\rightarrow \E _2 |_{\U}$ un morphisme de
$F \text{-} \mathrm{Isoc} ^{\dag \dag} ( \U, Y /K) $.
Avec les notations de \ref{defidagdagpxt?} et \ref{plfidked},
il découle de \ref{theo-surcohlisse} l'égalité suivante
$F \text{-}\mathrm{Isoc} ^{\dag \dag}( \PP', T',X' /K)= \linebreak
F \text{-}\mathrm{Isoc} ^{\dag}( \PP', T',X' /K)$.
Grâce au théorème de pleine fidélité de Kedlaya (voir \ref{plfidked}),
il existe un morphisme
$\phi$ :
$\R \underline{\Gamma} ^\dag _{X'} f ^! (\E _1) \rightarrow \R \underline{\Gamma} ^\dag _{X'} f ^! (\E _2)$
tel que $\phi |_{\U'}=\R \underline{\Gamma} ^\dag _{Y'} g ^! (\psi)$.
On conclut grâce à \ref{isosurcohpfpre}.
\hfill \hfill \qed \end{proof}

Intéressons-nous maintenant à l'indépendance en $\PP$ et $T$ de la catégorie
$F\text{-}\mathrm{Isoc} ^{\dag \dag}( \PP , T , X/K)$. 
\begin{lemm}\label{f+sp=sp}
  Soient $f$ : $\PP ' \rightarrow \PP$ un morphisme de $\V$-schémas formels lisses,
  $X$ un $k$-schéma lisse, $u $ : $X \hookrightarrow \PP$ et $u'$ : $X \hookrightarrow \PP '$ des immersions
  fermées telles que $f \circ u '=u$, $T  $ un diviseur de $P $ tel que $f ^{-1} (T )$ et $T _X := T \cap X$ soient
  respectivement des diviseurs de $P'$ et $X$.
  Pour tout isocristal sur $X \setminus T _X$ surconvergent  le long de $T _X$, on dispose d'un isomorphisme
  canonique : $f _{T ,+} \sp _{X \hookrightarrow \PP',f ^{-1}( T),+} (E)\riso \sp _{X \hookrightarrow \PP,T,+} (E)$.
\end{lemm}
\begin{proof}
Si $\PP$ est affine, il existe des relèvements
$u $ : $\X \hookrightarrow \PP$ et $u' $ : $\X \hookrightarrow \PP'$.
On a alors, en notant $\sp $ : $\X _K \rightarrow \X$ le morphisme de spécialisation
de $\X$,  $\sp _{X \hookrightarrow \PP,T,+} (E)\riso u _{T +} \sp _* (E)$
et $\sp _{X \hookrightarrow \PP',T',+} (E) \riso u ' _{f ^{-1} (T) +} \sp _* (E)$, ceux-ci ne dépendant pas,
à isomorphisme canonique près, du choix de $u$ et $u'$.
Il en résulte
$f _{T ,+} \sp _{X \hookrightarrow \PP',T',+} (E) \riso \sp _{X \hookrightarrow \PP,T,+} (E)$,
par composition des images directes.

Comme les isomorphismes de la forme $f _{T ,+} \circ u ' _{f ^{-1} (T) +} \riso u _{T,+} $
sont compatibles aux isomorphismes de recollement (induits par \cite[2.1.5]{Be2}),
on conclut la preuve par recollement.

\hfill \hfill \qed \end{proof}

\begin{theo} \label{theoindpdt1}
  Soient $f$ : $\PP ' \rightarrow \PP $ un morphisme propre et lisse de $\V$-schémas formels séparés et lisses et,
  $T $ (resp. $T'$) un diviseur et $X$ (resp. $X'$) un sous-schéma fermé de $P $ (resp. $P'$).
  On suppose de plus $X  \setminus T $ lisse et $f$ induisant l'isomorphisme
  $X  \setminus T  \riso X '\setminus T '$.
  On se donne $\E  \in F\text{-}\mathrm{Isoc} ^{\dag \dag} (\PP  , T  , X /K)$ et
$\E ' \in F\text{-}\mathrm{Isoc} ^{\dag \dag} (\PP ', T ', X '/K)$.
  Alors,
  \begin{enumerate}
    \item \label{1theoindpdt1}Pour tout $k \neq 0$, $\mathcal{H} ^k (f _+) (\E ') =0$ et
   $\mathcal{H} ^k  (\R \underline{\Gamma} ^\dag _{X '} \circ  f ^!(\E ))=0$ ;
    \item \label{2theoindpdt1} Les morphismes canoniques $f _+\circ \R \underline{\Gamma} ^\dag _{X '} \circ  f ^!(\E ) \rightarrow \E $
    et
    $\E' \rightarrow \R \underline{\Gamma} ^\dag _{X '} \circ  f ^!\circ f _+ (\E ')$ sont des isomorphismes.
  \end{enumerate}
  Les foncteurs $f_+$ et
  $\R \underline{\Gamma} ^\dag _{X '} \circ f ^!$ induisent
  des équivalences quasi-inverses entre les catégories
  $F\text{-}\mathrm{Isoc} ^{\dag \dag} (\PP ', T ', X '/K)$ et
  $F\text{-}\mathrm{Isoc} ^{\dag \dag} (\PP , T , X /K)$.
\end{theo}

    \begin{proof}
    Notons $\U  = \PP  \setminus T $, $\U '= \PP '\setminus T '$,
    $g$ : $\U ' \rightarrow \U $ le morphisme induit par $f$ et
    $Y =X  \setminus T $. Grâce à \ref{propindX}, on peut supposer $X $ (resp. $X'$)
    égal à l'adhérence schématique de $Y$ dans $P $ (resp. $P '$).

    Il résulte de \cite[4.3.12]{Be1} et de \cite[3.2.3]{caro_surcoherent}
    les égalités de \ref{1theoindpdt1}).

    \'Etablissons maintenant que l'équivalence de catégorie mentionnée dans l'énoncé
    a un sens.

    Comme $f$ est propre, l'immersion $Y \hookrightarrow P ' \setminus f ^{-1} (T  )$ est fermée.
    Comme $X '$ est l'adhérence de $Y$ dans $P '$, il en dérive que $Y =X' \setminus f ^{-1} (T  ) $.
    Grâce à \ref{propindX}, on obtient
    $F\text{-}\mathrm{Isoc} ^{\dag \dag} (\PP ', T ', X '/K) =
    F\text{-}\mathrm{Isoc} ^{\dag \dag} (\PP ', f ^{-1}( T ), X '/K)$.

    Il découle alors de \ref{isoctilstaiinvlis} que
    $\R \underline{\Gamma} ^\dag _{X '} \circ  f ^!(\E ) \in F\text{-}
    \mathrm{Isoc} ^{\dag \dag} (\PP ', T ', X '/K)$.
    De plus, par \ref{PPPTstable+propre}, $f _+ (\E ')$ vérifie $P _{\PP , T }$.
    Comme il existe un isocristal convergent $E '$ sur $Y$ tel que
    $\E ' |_{\U '} \riso \sp _{Y \hookrightarrow \U ',+}(E ')$, il résulte de \ref{f+sp=sp}
l'isomorphisme,
    $ g _+ \sp _{Y \hookrightarrow \U ',+}(E ') \riso \sp _{Y \hookrightarrow \U ,+}(E ')$.
    On a donc prouvé que $f _+ (\E ')  $ appartient à $F\text{-}\mathrm{Isoc} ^{\dag \dag} (\PP , T , X /K)$.

    Prouvons à présent \ref{2theoindpdt1}).
    Par adjonction, on a
    $\E '\tilde {\leftarrow}\R \underline{\Gamma} ^\dag _{X '}\E' \rightarrow
    \R \underline{\Gamma} ^\dag _{X '}  \circ f ^! \circ f _+ (\E')$.
    Comme ce composé est
    un isomorphisme au dessus de $\U '$ (\cite[3.2.3]{caro_surcoherent}),
    celui-ci est un isomorphisme (\cite[4.3.12]{Be1}).
    Avec les mêmes arguments, si $\E $ est un objet de
    $F\text{-}\mathrm{Isoc} ^{\dag \dag} (\PP , T , X /K)$,
    on établit que le morphisme canonique :
    $f _+ \circ \R \underline{\Gamma} ^\dag _{X '}  \circ f ^! \E  \rightarrow \E $ est un isomorphisme.
    \hfill \hfill \qed \end{proof}

\begin{vide}\label{prolong}
Soit $b $ : $Y ' \rightarrow Y$ un morphisme de $k$-schémas lisses.
On suppose qu'il existe $\PP$, $\PP'$ deux $\V$-schémas formels propres et lisses,
 $T$ (resp. $T'$) un diviseur de $P$ (resp. $P'$), $\U := \PP \setminus T$ et $\U':= \PP' \setminus T'$
 et des immersions fermées $Y \hookrightarrow \U$ et $Y' \hookrightarrow \U'$.

 On note $X$ (resp. $X'$) l'adhérence schématique de $Y$ (resp. $Y'$) dans $P$ (resp. $P'$),
$\PP'' := \PP '\times \PP $, $\U '' := \U '\times \U$, $X ''$ l'adhérence schématique de $Y $ (immergé via
le graphe de $b$) dans $X '\times X$, $f_1$ : $\PP''\rightarrow \PP'$ et
$f _2$ : $\PP ''\rightarrow \PP$ les projections, et $T '':= f _1^{-1} (T') \cup f _2^{-1} (T)$.
Grâce à \ref{theoindpdt1},
les foncteurs $f _{1+}$ et $\R \underline{\Gamma} ^\dag _{X''} f _{1} ^!$
induisent des équivalences quasi-inverses
entre les catégories
 $F\text{-}\mathrm{Isoc} ^{\dag \dag}( \PP'', T'', X''/K)$ et
 $F\text{-}\mathrm{Isoc} ^{\dag \dag}( \PP ', T', X'/K)$.

 Lorsque $b = id$, $\R \underline{\Gamma} ^\dag _{X''} f _{2} ^!$ et
 $f _{2 +}$ induisent des équivalences quasi-inverses
 entre $F\text{-}\mathrm{Isoc} ^{\dag \dag}( \PP , T, X/K)$
 et $F\text{-}\mathrm{Isoc} ^{\dag \dag}( \PP'', T'', X''/K)$.
 On obtient l'équivalence canonique de catégorie
 $$f _{1+} \circ \R \underline{\Gamma} ^\dag _{X''} f _{2} ^! \
 :\ F\text{-}\mathrm{Isoc} ^{\dag \dag}( \PP , T, X/K)\cong
 F\text{-}\mathrm{Isoc} ^{\dag \dag}( \PP ', T', X'/K)$$
 dont $f _{2+} \circ \R \underline{\Gamma} ^\dag _{X''} f _{1} ^! $ est un foncteur quasi-inverse.
 On remarque que lorsque $b$ se relève en un morphisme
 $f$ : $\PP ' \rightarrow  \PP$, le foncteur
 $f _{1+} \circ \R \underline{\Gamma} ^\dag _{X''} f _{2} ^!$
 (resp. $f _{2+} \circ \R \underline{\Gamma} ^\dag _{X''} f _{1} ^! $)
  est canoniquement isomorphisme à $\R \underline{\Gamma} ^\dag _{X'} f  ^!$ (resp. $ f_+$).
 De plus, ces isomorphismes sont transitifs en les choix $\PP$, $T$ et $Y \hookrightarrow \U$.
On écrira ainsi $F\text{-}\mathrm{Isoc} ^{\dag \dag}(Y/K)$ à la place de
$F\text{-}\mathrm{Isoc} ^{\dag \dag}( \PP, T, X/K)$.
De plus, en notant $p _{\PP}$ : $\PP \rightarrow \S$ le morphisme structural,
on vérifie que $p _{\PP +}[-d _Y]$ : $F\text{-}\mathrm{Isoc} ^{\dag \dag}( \PP, T, X/K) \rightarrow
D ^\mathrm{b} (K)$ est indépendant du choix de $(\PP, T, X)$. Il sera alors noté
$$H ^\bullet _{\mathrm{DR}}(-/K) \ : \ F\text{-}\mathrm{Isoc} ^{\dag \dag}(Y/K)
\rightarrow D ^\mathrm{b} (K). $$

Retournons au cas où $b$ n'est plus forcément l'identité.
Modulo l'équivalence de catégories précédente,
on peut toujours supposer que $b$ se prolonge en un morphisme
$f$ : $\PP ' \rightarrow \PP$ propre et lisse tel que $T ' \supset f ^{-1} (T)$.
Un tel prolongement sera appelé {\it bon prolongement}.
On dispose du foncteur
$\R \underline{\Gamma} ^\dag _{X'} f  ^! [-d _{Y'/Y}]$ :
$F \text{-}\mathrm{Isoc} ^{\dag \dag}( Y/K) \rightarrow
F \text{-}\mathrm{Isoc} ^{\dag \dag}( Y'/K)$.
Celui-ci ne dépend pas du choix du prolongement $f$ de $b$.
En effet, si $f _1$ : $\PP ' _1 \rightarrow \PP _1$ et
$f _2$ : $\PP '_2 \rightarrow \PP _2$ sont deux bons prolongements de $b$,
en notant
$f _3 = f _1  \times f _2$ : $ \PP ' _1 \times \PP ' _2 \rightarrow \PP _1\times \PP _2$,
$X '_1$, $X ' _2$ et $X' _3$ les adhérences schématiques respectives de
$Y'$ dans $P '_1$, $P' _2$ et $P ' _1 \times P ' _2$,
on vérifie que les foncteurs
$\R \underline{\Gamma} ^\dag _{X' _1} f _1 ^! [-d _{Y'/Y}]$
et $\R \underline{\Gamma} ^\dag _{X'_3} f _3 ^! [-d _{Y'/Y}]$ se correspondent
modulo l'équivalence canonique de catégorie entre
$F\text{-}\mathrm{Isoc} ^{\dag \dag}( \PP' _1, T '_1, X'_1/K)$ et
$F\text{-}\mathrm{Isoc} ^{\dag \dag}( \PP' _1 \times \PP ' _2 , T'_3, X'_3 /K)$,
où $T '_3$ est le diviseur réduit dont le support est le complémentaire de $\U' _1 \times \U ' _2$
dans $\PP' _1 \times \PP ' _2$. Et de même en remplaçant, "$1$" par "$2$".
Le foncteur $\R \underline{\Gamma} ^\dag _{X'} f  ^! [-d _{Y'/Y}]$ sera donc simplement noté
de la façon suivante
$$b ^* \ : F \text{-}\mathrm{Isoc} ^{\dag \dag}( Y/K) \rightarrow
F \text{-}\mathrm{Isoc} ^{\dag \dag}( Y'/K).$$

Enfin, les bifoncteurs de la forme $-\smash{\overset{\L}{\otimes}} ^\dag _{\O _{\PP} (\hdag T) _\Q}-[d _{X/P}]$
$$
F \text{-}\mathrm{Isoc} ^{\dag \dag}( \PP, T, X/K) \times
F \text{-}\mathrm{Isoc} ^{\dag \dag}( \PP, T, X/K)
\rightarrow
F \text{-}\mathrm{Isoc} ^{\dag \dag}( \PP, T, X/K)$$
commutent aux isomorphismes canoniques
$F \text{-}\mathrm{Isoc} ^{\dag \dag}( \PP, T, X/K) \cong \linebreak
F \text{-}\mathrm{Isoc} ^{\dag \dag}( \PP', T', X'/K)$.
En effet, on peut toujours supposer qu'il existe un bon prolongement
$f $ : $\PP' \rightarrow \PP$ de l'identité de $Y$.
L'isomorphisme canonique
$F \text{-}\mathrm{Isoc} ^{\dag \dag}( \PP, T, X/K) \cong
F \text{-}\mathrm{Isoc} ^{\dag \dag}( \PP', T', X'/K)$
est fourni par $\R \underline{\Gamma} ^\dag _{X'} f  ^!$.
Pour tous $\E _1,\E _2 \in F \text{-}\mathrm{Isoc} ^{\dag \dag}( \PP, T, X/K)$, on a alors
\begin{gather}
  \notag
  \R \underline{\Gamma} ^\dag _{X'} f  ^!
(\E _1 \smash{\overset{\L}{\otimes}} ^\dag _{\O _{\PP} (\hdag T) _\Q} \E _2 [d _{X/P}]  )
\riso
\R \underline{\Gamma} ^\dag _{X'}
(f ^! \E _1 \smash{\overset{\L}{\otimes}} ^\dag _{\O _{\PP'} (\hdag T') _\Q} f ^!  \E _2 )[d _{X/P'}]  )
\\ \riso
(\R \underline{\Gamma} ^\dag _{X'} f ^! \E _1 )
\smash{\overset{\L}{\otimes}} ^\dag _{\O _{\PP'} (\hdag T') _\Q}
(\R \underline{\Gamma} ^\dag _{X'} f ^! \E _2 )  [d _{X/P'}]  ),\notag
\end{gather}
le dernier isomorphisme résultant de \cite[2.1.8]{caro_surcoherent}.

On obtient ainsi le foncteur noté
$$-\smash{\overset{\L}{\otimes}} ^\dag _{\O _Y } -\ :\
F \text{-}\mathrm{Isoc} ^{\dag \dag}( Y/K) \times
F \text{-}\mathrm{Isoc} ^{\dag \dag}(Y/K)
\rightarrow
F \text{-}\mathrm{Isoc} ^{\dag \dag}(Y/K).$$

\end{vide}

\begin{vide}
  Soient $Y$ un $k$-schéma lisse, $\U$ un $\V$-schéma formel lisse et
  $Y \hookrightarrow \U$ une immersion fermée.

On dispose de l'équivalence de catégorie
$\sp _{Y \hookrightarrow \U, +}$ :
$F \text{-} \mathrm{Isoc} ^\dag ( Y,Y/K) \cong
F \text{-} \mathrm{Isoc} ^{\dag \dag} ( \U, Y/K)$.
On en déduit que
la catégorie $F \text{-} \mathrm{Isoc} ^{\dag \dag} (\U,Y/K)$ est indépendante de $\U$ à isomorphisme
canonique près.
  En effet, si $Y \hookrightarrow \U'$ est un second choix, en notant
  $g _1$ : $\U \times \U ' \rightarrow \U$ et
  $g _2$ : $\U \times \U ' \rightarrow \U'$, il découle du lemme \ref{f+sp=sp}
  et de \cite[2.2.27]{caro_surholonome}
  que les foncteurs
  $\R \underline{\Gamma} ^\dag _{Y} g _1 ^!$ et $g _{1+}$
  (resp. $\R \underline{\Gamma} ^\dag _{Y} g _2 ^!$ et $g _{2+}$)
  sont des équivalences quasi-inverses
  entre $F \text{-} \mathrm{Isoc} ^{\dag \dag} (\U,Y/K)$
  (resp. $F \text{-} \mathrm{Isoc} ^{\dag \dag} (\U',Y/K)$)
  et $F \text{-} \mathrm{Isoc} ^{\dag \dag} (\U \times \U',Y/K)$.
  Le foncteur $g _{2 +} \R \underline{\Gamma} ^\dag_{Y} g _1 ^!$ :
$F \text{-} \mathrm{Isoc} ^{\dag \dag} (\U,Y/K) \rightarrow
F \text{-} \mathrm{Isoc} ^{\dag \dag} (\U',Y/K)$ fournit une équivalence
de catégorie canonique. On écrira donc sans ambiguïté
$F \text{-} \mathrm{Isoc} ^{\dag \dag} ( Y, Y/K)$ à la place de
$F \text{-} \mathrm{Isoc} ^{\dag \dag} ( \U, Y/K)$.

  De plus,
le foncteur
  $\sp _{Y \hookrightarrow \U, +}$ :
$F \text{-} \mathrm{Isoc} ^\dag ( Y,Y/K) \cong
F \text{-} \mathrm{Isoc} ^{\dag \dag} ( \U, Y/K)$,
grâce à \ref{f+sp=sp} et à \cite[2.2.27]{caro_surholonome},
commute à ces isomorphismes. On écrira alors
$\sp _{Y ,Y, +}$ :
$F \text{-} \mathrm{Isoc} ^\dag ( Y,Y/K) \cong
F \text{-} \mathrm{Isoc} ^{\dag \dag} ( Y, Y/K)$
ce foncteur.

Soit $b$ : $ Y' \rightarrow Y$ un morphisme de $k$-schémas lisses
tel qu'il existe $\U$ et $\U'$ des $\V$-schémas formels lisses et des immersions fermées
$Y \hookrightarrow \U$ et $Y ' \hookrightarrow \U'$.
On dispose d'un foncteur canonique
$b ^*$ : $F \text{-} \mathrm{Isoc} ^{\dag \dag} ( Y, Y/K) \rightarrow
F \text{-} \mathrm{Isoc} ^{\dag \dag} ( Y', Y'/K)$ bien défini.
En effet, pour tout
$\E \in F \text{-} \mathrm{Isoc} ^{\dag \dag} ( \U, Y/K)$, on pose
$b ^* (\E):= g _{2+} \R \underline{\Gamma} _{Y'} ^\dag g _1 ^! (\E)$,
où $g _1$ : $\U  \times \U ' \rightarrow \U $ et
$g _2$ : $\U  \times \U ' \rightarrow \U '$.
Ceux-ci sont indépendants des choix de
$\U$ et de $\U'$.
De plus, on dispose d'un isomorphisme canonique
$ b ^* \circ \sp _{Y ,Y, +} \riso \sp _{Y ',Y', +} \circ b ^*$.

En reprenant les notations de la section \ref{prolong},
on remarque que le foncteur $| _\U$ :
$F \text{-} \mathrm{Isoc} ^{\dag \dag} ( \PP, T,X/K)
\rightarrow F \text{-} \mathrm{Isoc} ^{\dag \dag} ( \U, Y/K)$
est indépendant des choix faits à isomorphisme canonique près.
On le notera $\mathrm{cv} _{_Y}$ : $F \text{-} \mathrm{Isoc} ^{\dag \dag} ( Y/K)
\rightarrow F \text{-} \mathrm{Isoc} ^{\dag \dag} ( Y, Y/K)$.
Pour tout $\E \in F \text{-} \mathrm{Isoc} ^{\dag \dag} ( Y/K)$,
on a l'isomorphisme canonique
$ \mathrm{cv} _{_{Y'}} \circ b ^* (\E) \riso b ^* \circ \mathrm{cv} _{_Y} (\E)$
fonctoriel en $\E$.
\end{vide}

 \begin{vide}
 Soient $X$ un $k$-schéma lisse,  $\PP$ un $\V$-schéma formel propre et lisse,
 $u$ : $X \hookrightarrow \PP$ une immersion fermée et
 un diviseur $T$ de $P$ tel que $T \cap X$ soit un diviseur de $X$.

   Le foncteur $\sp _{X \hookrightarrow \PP, T,+}$ :
   $F \text{-} \mathrm{Isoc} ^\dag ( Y/K) \rightarrow
   F \text{-} \mathrm{Isoc} ^{\dag \dag} ( \PP, T,X/K)$, qui est une équivalence
   de catégorie d'après \ref{dag2=daglisse}, est indépendant des choix
   de $X \hookrightarrow \PP$ et de $T$ vérifiant les conditions ci-dessus.
  En effet, si $X ' \hookrightarrow \PP '$ et $T'$ sont un second choix, on note
  $\U := \PP \setminus T$, $\U ':= \PP ' \setminus T'$,
  $f _1$ : $\PP \times \PP ' \rightarrow \PP$, $f _2$ : $\PP \times \PP ' \rightarrow \PP '$,
  $g _1 $ : $\U \times \U ' \rightarrow \U$, $g _2 $ : $\U \times \U ' \rightarrow \U'$
  et $X ''$ l'adhérence schématique de $Y$ dans $\PP \times \PP'$.
  Il s'agit de vérifier l'isomorphisme :
  $\R \underline{\Gamma} ^\dag _{X''} f _2 ^! \sp _{X' \hookrightarrow \PP', T',+} (E)
  \riso \R \underline{\Gamma} ^\dag _{X''} f _1 ^!  \sp _{X \hookrightarrow \PP, T,+}(E)$,
  pour tout $E \in F \text{-} \mathrm{Isoc} ^\dag ( Y/K)$.
  D'après \ref{pffidcasgen} et \ref{theoindpdt1}, comme ceux-ci sont
  dans $F \text{-}\mathrm{Isoc} ^\dag ( \PP \times \PP', f _1 ^{-1} (T ) \cup f _2^{-1} (T '), X''/K)$
  il s'agit de le voir au dessus de $\U \times \U'$, ce qui résulte,
  en notant $\widehat{E}$ l'isocristal convergent sur $Y$ associé à $E$,
  des isomorphismes $\R \underline{\Gamma} ^\dag _{Y} g _1 ^!  \sp _{Y \hookrightarrow \U, +}(\widehat{E})
  \riso \sp _{Y \hookrightarrow \U \times \U ', +}(\widehat{E})$
  et \linebreak
  $\R \underline{\Gamma} ^\dag _{Y} g _2 ^!  \sp _{Y \hookrightarrow \U ', +}(\widehat{E})
  \riso
  \sp _{Y \hookrightarrow \U \times \U ', +}(\widehat{E})$.

On le notera alors
$\sp _{Y +}$ :
$F \text{-} \mathrm{Isoc} ^\dag ( Y/K) \cong F \text{-} \mathrm{Isoc} ^{\dag \dag} ( Y/K)$.
On bénéficie de l'isomorphisme 
$\mathrm{cv} _{_Y} \circ \sp _{Y+} \riso \sp _{Y,Y+} \circ \mathrm{cv} _Y $, où
$\mathrm{cv} _{_Y}$ est le foncteur canonique
$F \text{-} \mathrm{Isoc} ^\dag ( Y/K) \rightarrow F \text{-} \mathrm{Isoc} ^\dag (Y, Y/K)$
(et de même en remplaçant {\og $^\dag $\fg} par {\og $^{\dag \dag}$\fg}).

 \end{vide}

\begin{theo}\label{isosurcohpf}
On suppose $\PP$ propre.
  Il existe un foncteur (non canonique) pleinement fidèle
  $$\rho _Y \ : \ F\text{-}\mathrm{Isoc} ^{\dag \dag}(Y/K)\rightarrow F\text{-}\mathrm{Isoc} ^\dag( Y /K).$$

  Les objets de $F\text{-}\mathrm{Isoc} ^{\dag \dag}(Y/K)$ serons alors appelés
  {\og $F$-isocristaux surcohérents sur $Y $ \fg}. 
\end{theo}
\begin{proof}
  Par \ref{isosurcohpf0} et avec ses notations,
  comme $F\text{-}\mathrm{Isoc} ^\dag( Y /K) \cong \oplus _r F\text{-}\mathrm{Isoc} ^\dag( Y_r /K)$,
  on se ramène au cas où $Y$ est irréductible.
  D'après \ref{isosurcohpfpre} et avec ses notations, le foncteur
$(b ^* , \mathrm{cv} _Y)$ :
  \begin{equation}\notag 
 F \text{-}\mathrm{Isoc} ^{\dag \dag }( Y /K)\rightarrow
F \text{-}\mathrm{Isoc} ^{\dag \dag}(Y' /K) \times _{F \text{-}\mathrm{Isoc} ^{\dag \dag} (Y ', Y' /K) }
F \text{-}\mathrm{Isoc} ^{\dag \dag} (Y, Y /K)
\end{equation}
  est alors pleinement fidèle.
Or, il résulte de \cite[Théorème 3]{Etesse-descente-etale} et de \cite{kedlaya_full_faithfull},
que le foncteur canonique
\begin{equation}\notag 
  F \text{-}\mathrm{Isoc} ^\dag( Y /K) \rightarrow
F \text{-}\mathrm{Isoc} ^{\dag }(Y' /K) \times _{F \text{-}\mathrm{Isoc} ^{\dag } (Y ', Y' /K) }
F \text{-}\mathrm{Isoc} ^{\dag } (Y, Y /K)
\end{equation}
est une équivalence de catégorie.

Comme $X ' $ est lisse,
$\sp _{Y '+}$ :
$F \text{-} \mathrm{Isoc} ^\dag ( Y'/K) \cong  F \text{-} \mathrm{Isoc} ^{\dag \dag} ( Y'/K)$.
De plus,
$\sp _{Y ,Y, +}$ :
$F \text{-} \mathrm{Isoc} ^\dag ( Y,Y/K) \cong
F \text{-} \mathrm{Isoc} ^{\dag \dag} ( Y, Y/K)$, et de même en remplaçant $Y $ par $Y'$.

\hfill \hfill \qed \end{proof}

\begin{rema}
\label{rhocv=cv}
Soit $b $ : $Y ' \rightarrow Y$ un morphisme de $k$-schémas lisses.
On suppose qu'il existe $\PP$, $\PP'$ deux $\V$-schémas formels propres et lisses,
 $T$ (resp. $T'$) un diviseur de $P$ (resp. $P'$),
$Y \hookrightarrow \U$ et $Y' \hookrightarrow \U'$ des immersions fermées,
avec $\U := \PP \setminus T$ et $\U':= \PP' \setminus T'$.

  Pour tout $\E \in F \text{-}\mathrm{Isoc} ^{\dag \dag}(Y/K)$,
$\rho _{Y} (\E)$ est l'unique (à isomorphisme près) $F$-isocristal $E$ surconvergent sur $Y$
tel que $\mathrm{cv} _Y  ( \E ) \riso \sp _{Y,Y+} ( \widehat{E})$, où
$\widehat{E}$ est le $F$-isocristal convergent sur $Y$ associé à $E$.
Il en résulte un isomorphisme
 $\rho _{Y'} \circ b ^* (\E) \riso b ^* \circ \rho _Y (\E)$.
De même, on obtient, pour tous $\E _1,  \E_2 \in F \text{-}\mathrm{Isoc} ^{\dag \dag}(Y/K)$,
$ \rho _Y (\E _1 \smash{\overset{\L}{\otimes}} ^\dag _{\O _Y } \E _2)
\riso \rho _Y (\E _1) \otimes \rho _Y (\E _2)$.

\end{rema}

\subsection{Surcohérence générique des $F$-isocristaux surconvergents}

\begin{defi}\label{def-notation-surcoh-gen}
  Soient $Y _0$ un $k$-schéma affine et lisse,
  $Y _0 \hookrightarrow X _0$ une immersion ouverte et $Y ^\dag$ un $\V$-schéma formel faible affine et lisse
  relevant $Y _0$.
 L'immersion ouverte $Y _0 \hookrightarrow X_0$ {\og se désingularise idéalement \fg}
  s'il existe un morphisme surjectif $a _0$ : $X'_0 \rightarrow X _0$, qui se décompose
  en une immersion fermée $X '_0 \hookrightarrow \P ^r _{X _0}$ suivie de la projection
  canonique $\P ^r _{X _0} \rightarrow X _0$, tel que
  \begin{enumerate}
    \item \label{def-notation-surcoh-genitem1} $X ' _0$ est lisse ;
    \item \label{def-notation-surcoh-genitem2} Le morphisme $Y ' _0 := a _0 ^{-1} (Y _0) \rightarrow Y _0$ induit par $a_0$ est
    fini et étale ;
    \item \label{def-notation-surcoh-genitem3} Le morphisme $Y '_0 \hookrightarrow \P ^r _{Y _0}$ induit par
    $X '_0 \hookrightarrow \P ^r _{X _0}$
   se relève en un morphisme de $\V$-schémas formels faibles lisses de la forme
   $ Y ^{'\dag} \rightarrow \P ^{r \dag} _{Y ^\dag}$.
  \end{enumerate}

  L'immersion ouverte $Y _0 \hookrightarrow X_0$ se désingularise {\it localement idéalement} s'il
  existe un recouvrement ouvert $( X _{0,i}) _{i\in I}$ de $X _0$ tel que
  $Y _0\cap X _{0,i} \hookrightarrow  X _{0,i}$ se désingularise idéalement.
\end{defi}

\begin{rema}\label{remalocide}
Par exemple et avec les notations de \ref{def-notation-surcoh-gen},
lorsque $X _0$ est lisse, $Y _0 \hookrightarrow X_0$ se désingularise idéalement.
  De plus, la condition \ref{def-notation-surcoh-genitem3}
  est indépendante du choix
  du relèvement $ Y^\dag$ car ceux-ci sont isomorphes.
  Enfin, si
$Y _0 \hookrightarrow X_0$ se désingularise idéalement alors il en est de même,
pour tout ouvert affine $\widetilde{Y} _0$ de $Y _0$,
de l'immersion ouverte induite $\widetilde{Y} _0 \hookrightarrow X_0$.
En effet, le morphisme $\smash{\widetilde{Y}} ' _0 :=
a _0 ^{-1} (\widetilde{Y} _0) \rightarrow \widetilde{Y} _0$ induit par $a_0$ est fini et étale.
De plus, en notant $\widetilde{Y} ^\dag $ l'ouvert de $Y ^\dag$ d'espace sous-jacent $\widetilde{Y} _0$,
la projection $ Y ^{'\dag} \times _{\P ^{r \dag} _{Y ^\dag}} \P ^{r \dag} _{\widetilde{Y} ^\dag}
\rightarrow \P ^{r \dag} _{\widetilde{Y} ^\dag}$ est un relèvement de
$\smash{\widetilde{Y}} ' _0 \rightarrow \P ^r _{\widetilde{Y} _0}$.
\end{rema}

\begin{lemm}\label{theo-ideal->oKlemm}
  Soit le diagramme de $\V$-schémas formels faibles lisses
\begin{equation}\label{theo-ideal->oK-diagm1}
  \xymatrix @R=0,4cm {
{Y ^{'\dag}}
\ar[r] ^{v'} \ar[d] ^b
&
{U ^{'\dag} }
\ar[r] ^{j'} \ar[d] ^g
&
{P ^{'\dag} }
\ar[d] ^f
\\
{Y ^{\dag}}
\ar[r] ^{v}
&
{U ^{\dag} }
\ar[r] ^{j}
&
{P ^{\dag} ,}
}
\end{equation}
  où $f$ est propre, $j$ est une immersion ouverte, le carré de droite est cartésien,
  $v$ et $v'$ sont des immersions fermées, $b$ est fini et étale et $Y ^{\dag}$ est affine.
  On suppose de plus que $T _0 := P _0\setminus U _0$ et $ T ' _0:= P '_0\setminus U '_0$ sont
  les supports de diviseurs.
  On désigne par $b _*$ le foncteur canonique de la catégorie des
isocristaux surconvergents sur $Y '_0$ dans
celle des isocristaux surconvergents sur $Y _0$ (cela a un sens car $b$ est fini, étale).

  Pour tout isocristal surconvergent $E'$ sur $Y '_0$, on dispose alors d'un isomorphisme canonique
  \begin{equation}\label{theo-ideal->oK-lemmmorph}
\sp _{Y ^{\dag}  \hookrightarrow U ^{\dag},T _0 , +} (b _* E ')
\riso
f _{T _0, + ^\dag} (\sp _{Y ^{\prime \dag}  \hookrightarrow U ^{\prime \dag}, T '_0 , +} (E')).
\end{equation}
En outre, ceux-ci sont compatibles aux compositions de diagrammes
de la forme \ref{theo-ideal->oK-diagm1}.
\end{lemm}
\begin{proof}
  On note $\E '$ le $\D _{Y ^{\prime \dag} ,\Q}$-module associé
à $E'$, $\E ^{\prime (0)}$, un $\D ^{(0)} _{Y ^{\prime \dag}}$-module globalement de présentation finie et
$\O _{Y ^{\prime \dag} }$-cohérent tel que $\E ^{\prime (0)} _\Q \riso \E'$ et
on pose $\FF ^{\prime (0)}:= v ' _+ (\E  ^{\prime (0)})$.
Un modèle $\D ^{(0)} _{Y ^{\dag}}$-globalement de présentation finie et $\O _{Y ^{\dag}}$-cohérent de $b _*(E')$ est
donné par $b _+ (\E  ^{\prime (0)})$.

Comme $g _+ (\FF ^{\prime (0)}) =g _+ v ' _+ (\E  ^{\prime (0)}) \riso v _+ b _+ (\E  ^{\prime (0)})$,
il découle de \ref{globpresenfiniu+} que
$\FF ^{\prime (0)} $ (resp. $g _+ (\FF ^{\prime (0)})$)
est localement en $P ^{\prime \dag}$ (resp. $P ^\dag$) de présentation finie.
D'après \ref{g-+commdag*f+}, on dispose alors du morphisme
\begin{gather}\notag
  \D ^\dag _{\PP } (\hdag * T _0) \otimes _{j _* \D ^{(0)} _{U ^\dag}} j _* g _+ (\FF ^{\prime (0)})
  \\
\rightarrow
  \R f _* (
  \D ^\dag _{\PP \leftarrow \PP '} (\hdag *T _0 )
  \otimes _{\D ^\dag _{\PP '} (\hdag * T '_0 )} ^\L
  \D ^\dag _{\PP '} (\hdag{* T '_0 })
  \otimes ^\L  _{j '_* \D ^{(0)} _{U ^{\prime \dag}}} j '_* \FF ^{\prime (0)}).\notag 
\end{gather}
En lui appliquant $\Q \otimes _\Z -$, 
on obtient le morphisme
\ref{theo-ideal->oK-lemmmorph}.
La cohérence différentielle étant préservée par image directe par un morphisme propre,
\ref{theo-ideal->oK-lemmmorph} est un morphisme de $\D ^\dag _{\PP } (\hdag T _0 )$-modules cohérents.
Comme \ref{theo-ideal->oK-lemmmorph} est un isomorphisme au dessus de $\U$,
par \cite[4.3.12]{Be1}, celui-ci est bien un isomorphisme.
\hfill \hfill \qed \end{proof}

\begin{lemm}
  Soient $b$ : $ \Y' \rightarrow \Y$ un morphisme fini et étale de $\V$-schémas formels lisses,
  $\widehat{E}$ un isocristal convergent sur $Y _0$ et $\widehat{\E}:= \sp _* (\widehat{E})$.
On dispose du diagramme commutatif suivant
\begin{equation}\label{rho=rho}
  \xymatrix @R=0,3cm {
{b _{+^\dag} b ^{!^\dag} \widehat{\E}}
\ar[r] ^-\rho
\ar[d] ^-\sim
&
{\widehat{\E}}
\ar@{=}[d]
\\
{ \sp _* ( b _* b ^*\widehat{E})}
\ar[r] ^-\rho
&
{\sp _* (\widehat{E}),}
}
\end{equation}
où le morphisme du haut est le morphisme d'adjonction, celui du bas
résulte du morphisme trace $b _* b ^* (\widehat{E}) \rightarrow \widehat{E}$
et celui de gauche découle de \ref{prop-comm-^inv2} et \ref{theo-ideal->oK-lemmmorph}.
\end{lemm}
\begin{proof}
Par construction du morphisme Trace de Virrion (\cite{Vir04}),
le diagramme
$\xymatrix {
{b_* (\omega _{\Y' ,\Q})} \ar@/^1pc/[rr] ^{\mathrm{Tr} _{b}} \ar[r] _-\sim  &
{b _+ (\omega _{\Y' ,\Q}) } \ar[r]  _-{\mathrm{Tr} _{+,b}}
& {\omega _{\Y, \Q}}}$ est commutatif.
On termine la preuve par construction des morphismes d'adjonction et trace respectifs.
\end{proof}

\begin{vide}\label{notation-surcoh-gen}
  Soient $P ^\dag$ un $\V$-schéma formel faible lisse et séparé, $T _0$ un diviseur de $P _0$, $U ^\dag$ l'ouvert de
$P ^\dag $ complémentaire de $T _0$, $j$ : $U ^\dag \hookrightarrow P ^\dag$ l'immersion
  ouverte et $v$ : $Y ^\dag \hookrightarrow U ^\dag$ une immersion fermée
de $\V$-schémas formels faibles.
On suppose en outre $Y ^\dag$ affine et lisse et on désigne par $X _0$ l'adhérence schématique de $Y _0$ dans $P _0$.

On remarque que lorsque $X _0$ est lisse, alors $T _0 \cap X _0$ est un diviseur de $X _0$. En effet,
comme $X _0$ est la somme directe de ses composantes irréductibles, on se ramène à traiter le cas où $X _0$
est irréductible. Celui-ci résulte du fait que l'on peut avoir $T _0\supset X _0$.
  \end{vide}

\begin{theo}\label{theo-ideal->oK}
Avec les notations \ref{notation-surcoh-gen}, on suppose que
$Y _0 \hookrightarrow X_0$ se désingularise localement idéalement.

Alors, pour tout isocristal $ E$ sur $Y _0$ surconvergent, 
le $\D ^\dag _{\PP } (\hdag T _0) _\Q$-cohérent $\sp _{Y ^{\dag}  \hookrightarrow U ^{\dag},T _0 , + }(E)$
 vérifie $P _{\PP, T_0}$ (voir \ref{PPPT}).

 D'autre part, l'isomorphisme canonique de commutation à Frobenius
  $\phi$ : $\sp _{Y ^{\dag}  \hookrightarrow U ^{\dag},T _0 , + } (F ^* E)
  \riso F ^* \sp _{Y ^{\dag}  \hookrightarrow U ^{\dag},T _0 , + } ( E)$
  existe (voir \ref{isocanfrob}), i.e., le diagramme ci-dessous
  \begin{equation}
    \label{theo-ideal->oKdiagp}
\xymatrix @R=0,4cm {
{\sp _{Y ^{\dag}  \hookrightarrow U ^{\dag},T _0 , + } (F ^* E) |_{\U}}
\ar[r] ^{\phi |_{\U }} _-\sim \ar[d] ^-\sim
&
{ F ^* (\sp _{Y ^{\dag}  \hookrightarrow U ^{\dag},T _0 , + } (E) |_{\U})}
\ar[d] ^-\sim
\\
{v _{+ ^{\dag}} F ^* \widehat{\E}}
\ar[r] ^-\sim
&
{F ^* v _{+ ^{\dag}} \widehat{\E},}
}
  \end{equation}
  où les isomorphismes verticaux sont \ref{spyu+rest} et celui du bas est l'isomorphisme de commutation
  à Frobenius de l'image directe, est commutatif.
  Ainsi, on dispose d'un foncteur \newline
  $\sp _{Y ^{\dag}  \hookrightarrow U ^{\dag},T _0 , + }$ :
  $F \text{-}\mathrm{Isoc} ^\dag (Y _0/K) \rightarrow F \text{-}\mathrm{Isoc} ^{\dag \dag}  (\PP, T _0, X _0/K)$.

\end{theo}
\begin{proof}
  Comme la propriété $P _{\PP, T_0}$ est locale en $\PP$ (voir \ref{remaPPPT}) ainsi
  que l'existence de l'isomorphisme canonique de commutation à Frobenius, on peut supposer
  que $Y _0 \hookrightarrow X_0$ se désingularise idéalement.
  Par hypothèse, il existe alors un morphisme surjectif $a _0$ : $X'_0 \rightarrow X _0$, qui se décompose
  en une immersion fermée $X '_0 \hookrightarrow \P ^r _{X _0}$ suivie de la projection
  canonique $\P ^r _{X _0} \rightarrow X _0$,
  tel que $X ' _0$ soit lisse, le morphisme
  $ Y' _0 = a _0 ^{-1}(Y _0) \rightarrow Y _0$ soit fini et étale
  et tel qu'il existe un relèvement
  $ Y ^{'\dag} \rightarrow \P ^{r \dag} _{Y ^\dag}$
  de $Y '_0 \hookrightarrow \P ^r _{Y _0}$.
  En composant ce relèvement avec la projection canonique
  $\P ^{r \dag} _{Y ^\dag} \rightarrow Y ^\dag$ (resp. avec l'immersion fermée
  $\P ^{r \dag} _{Y ^\dag} \hookrightarrow \P ^{r \dag} _{U ^\dag}$), on obtient un morphisme fini, étale et surjectif
  $b$ : $ Y ^{'\dag} \rightarrow Y ^\dag $ (resp. une immersion fermée
  $v'$ : $ Y ^{'\dag} \hookrightarrow \P ^{r \dag} _{U ^\dag}$).
En notant $U ^{'\dag} := \P ^{r \dag} _{U ^\dag}$, $P ^{'\dag} := \P ^{r \dag} _{P ^\dag}$,
$j '$ : $U ^{'\dag} \hookrightarrow P ^{'\dag}$ l'immersion ouverte,
$f$ : $P^{'\dag} \rightarrow P ^\dag$ et
$g$ : $U^{'\dag} \rightarrow U ^\dag$ les projections, on
obtient un diagramme commutatif de la forme \ref{theo-ideal->oK-diagm1}, avec en plus $b$ surjectif.

Soit $E$ un isocristal surconvergent sur $Y _0$.
On note $\E$ le $\D _{Y ^\dag ,\Q}$-module associé
à $E$ et $\E ^{(0)}$, un $\D ^{(0)} _{Y ^\dag}$-module globalement de présentation finie et $\O _{Y ^\dag }$-cohérent
tel que $\E ^{(0)} _\Q \riso \E$.
On pose $\FF ^{(0)} := v _+ ( \E ^{(0)})$,
$\E  ^{\prime (0)} := b ^* (\E ^{(0)})$ et $\FF ^{\prime (0)}:= v ' _+ (\E  ^{\prime (0)})$.
On remarque que $\E  ^{\prime (0)}$
est un modèle $\D ^{(0)} _{Y ^{\prime \dag}}$-globalement de présentation finie et $\O _{Y ^{\prime \dag}}$-cohérent
de $b ^* (E)$
et $b _+ (\E  ^{\prime (0)})$
est un modèle $\D ^{(0)} _{Y ^{\dag}}$-globalement de présentation finie et $\O _{Y ^{\dag}}$-cohérent
de $b _* b ^* (E)$.

D'après \ref{theo-ideal->oKlemm}, on dispose de l'isomorphisme canonique
  \begin{equation}\label{theo-ideal->oK-morph2}
\sp _{+} (b _* b ^* E)
\riso
f _{T _0, + ^\dag} ( \sp _+ ^\prime (b ^* E)),
\end{equation}
où $\sp _+ =\sp _{Y ^{\dag}  \hookrightarrow U ^{\dag},T _0 , +}$
et $\sp _+ ^\prime = \sp _{Y ^{\prime \dag}  \hookrightarrow U ^{\prime \dag}, f ^{-1} T _0 , +}$.

Comme $X ' _0$ est lisse, avec la remarque de \ref{notation-surcoh-gen} et via la description
de l'image essentielle de $\widetilde{\sp} ' _+ := \sp _{X '_0 \hookrightarrow \PP ', f ^{-1} T _0 , +}$,
il existe un isocristal $E'$ sur $Y '_0$ surconvergent le long de $f ^{-1} T _0$
et vérifiant
$\smash{\widetilde{\sp}} ' _+  (E') \riso \sp _+ ^\prime  (b ^* E)$.
Or, d'après \ref{theo-surcohlisse}, $\smash{\widetilde{\sp}} ' _+  (E') $
vérifie $P _{\PP ',f ^{-1} T _0}$.
Avec \ref{PPPTstable+propre} et \ref{theo-ideal->oK-morph2},
il en dérive que
$\sp _+  (b _* b ^* E)$ vérifie
$P _{\PP,T _0}$.

Or, puisque $b$ est fini, étale et surjectif, $E$ est un facteur direct de
$b _* b ^* E$. Le module
$\sp _+  (E)$
est donc facteur direct de
$\sp _+  (b _* b ^* E)$.
D'où le premier résultat.

Prouvons à présent l'isomorphisme de commutation à Frobenius.
Dans la suite de la preuve, nous ne considérons plus d'opérations cohomologiques
définies sur les schémas formels faibles et
nous ne noterons plus le symbole {\og $^\dag$ \fg} pour
désigner les opérations cohomologiques (pour les $\D ^\dag$-complexes)
car aucune confusion ne sera à craindre.

Notons $\theta$ le morphisme défini par le diagramme commutatif ci-dessous
\begin{equation}
  \label{theo-ideal->oKdiagp1}
\xymatrix @R =0,3cm {
{\sp _{ + } b _* b ^* (F ^* E)}
\ar@{.>}[ddd] ^-\sim _{\theta} \ar[r] ^-\sim
&
{f _{T _0 +} (\sp _+ ^\prime  (b ^* F ^*  E))}
\ar[d] ^-\sim
\\
&
{f _{T _0 +} (\sp _+ ^\prime  (F ^* b ^*   E))}
\ar[d] ^-\sim
\\
&
{f _{T _0 +}F ^* (\sp _+ ^\prime  ( b ^*   E))}
\ar[d] ^-\sim
\\
{F ^*\sp _{ + } b _* b ^* (E)}
\ar[r] ^-\sim
&
{F ^* f _{T _0 +} (\sp _+ ^\prime  ( b ^*   E)),}
}
\end{equation}
où l'isomorphisme à droite du milieu découle de \ref{sp+froblisse} et ceux horizontaux
résultent de \ref{theo-ideal->oKlemm}.

Soit $\phi$ l'unique
morphisme rendant commutatif le diagramme ci-après
\begin{equation}
  \label{theo-ideal->oKdiagp2}
  \xymatrix @R=0,3cm {
  {\sp _+  (F ^* E)}
  \ar[r]  \ar@{.>}[d] ^{\phi}
  &
  {\sp _+  b _* b ^* (F ^* E)}
  \ar[r] ^-\rho \ar[d] ^{\theta}
  &
  {\sp _+  (F ^* E)}
  \ar@{.>}[d] ^{\phi}
  \\
  {F ^*\sp _+  ( E)}
  \ar[r]
  &
  {F ^*\sp _+  b _* b ^* (E)}
  \ar[r] ^-\rho
  &
  {F ^* \sp _+  (E),}
}
\end{equation}
où les morphismes horizontaux de gauche (resp. de droite)
découlent des morphismes canoniques $id \rightarrow b _* b ^*$
(resp. des morphismes traces $b _* b ^* \rightarrow id$).
En effet, cela est toujours possible car le composé des morphismes horizontaux sont des isomorphismes.
Par \cite[4.3.12]{Beintro2}, il suffit de prouver que $\phi$ induit
le diagramme commutatif \ref{theo-ideal->oKdiagp}.

Considérons le diagramme ci-après :
\begin{equation}
\label{theo-ideal->oKdiagp3}
  \xymatrix @R =0,3cm @C=0,25cm {
 {(\sp _{ _+ } b _* b ^* F ^* E) _{|\U}}
\ar[ddd] ^-\sim _{\theta |_{\U}}
\ar[r] ^-\sim
&
{g _{_+ } v ' _{_+ } b ^{! } F ^* \widehat{\E}}
\ar[r]\ar[d]
&
{v _{_+ } b _{_+ } b ^{! } F ^* \widehat{\E}}
\ar[d] \ar@<1ex>[rdd] ^-{\rho}
\\
&
{g _{_+ } v ' _{_+ } F ^*  b ^{! } \widehat{\E}}
\ar[r]\ar[d]
&
{v _{_+ } b _{_+ } F ^* b ^{! }  \widehat{\E}}
\ar[d]
\\
&
{g _{_+ } F ^* v ' _{_+ }   b ^{! } \widehat{\E}}
\ar[d]
&
{v _{_+ } F ^* b _{_+ }  b ^{! }  \widehat{\E}}
\ar[d]\ar[r]
&
{v _{_+ } F ^* \widehat{\E}}
\ar[d] ^\psi
\ar[r]
&
{(\sp _{ _+ } F ^* E) _{|\U}}
\ar[d] ^{\phi |_\U}
\\
{(F ^*\sp _{ _+ } b _* b ^* E) _{|\U}}
\ar[r] ^-\sim
&
{F ^* g _{_+ }  v ' _{_+ }   b ^{! } \widehat{\E}}
\ar[r]
&
{F ^* v _{_+ } b _{_+ }  b ^{! }  \widehat{\E}}
\ar[r] ^-\rho
&
{F ^* v _{_+ }  \widehat{\E}}
\ar[r]
&
{(F ^*\sp _{ _+ }  E) _{|\U},}
}
\end{equation}
où $\widehat{\E} := \sp _* (\widehat{E})$ et
$\psi$ : $v _{+ } F ^* \widehat{\E} \riso F ^* v _{+ } \widehat{\E}$
est l'isomorphisme canonique.
Via \ref{theo-ideal->oKdiagp1} et avec les isomorphismes canoniques compatibles à Frobenius
$\sp _* b ^* (\widehat{E}) \riso b ^{! } \sp _* \widehat{E}$ (déduit de \ref{prop-comm-^inv2}) 
et $\sp _+ ^\prime  (\widehat{E}') \riso  v '   _{+ }\sp _*  (\widehat{E}')$ (voir \ref{spyu+rest})
fonctoriels respectivement en $E \in F\text{-}\mathrm{Isoc} ^\dag (Y_0, X _0/K)$ et
$E '\in F\text{-}\mathrm{Isoc} ^\dag (Y'_0, X ' _0/K)$,
on vérifie que le grand rectangle de gauche est commutatif.
Comme le morphisme d'adjonction
$b _+ b ^! \rightarrow id $ est compatible à Frobenius,
le triangle (de la troisième colonne) est commutatif.
De plus, comme l'isomorphisme de commutation à la composition
des images directes est compatible à Frobenius (voir \cite{caro_courbe}),
le rectangle de la deuxième colonne est commutatif.
Enfin, les deux carrés autres que celui de droite le sont par fonctorialité.
Or, avec l'aide de \ref{rho=rho}, on vérifie que le grand contour du diagramme
\ref{theo-ideal->oKdiagp3} est l'image par $|_\U$ du carré commutatif de droite de
\ref{theo-ideal->oKdiagp2}.
Comme les flèches du carré de droite de
\ref{theo-ideal->oKdiagp3} sont des isomorphismes, il en résulte que
celui-ci est commutatif.

\hfill \hfill \qed \end{proof}

\begin{rema}
  \label{remspcv=cv}
  Avec les notations et hypothèses de \ref{theo-ideal->oK},
  on vérifie grâce à \ref{spyu+rest} et \ref{theo-ideal->oK} que le diagramme
  $$\xymatrix @ C=2cm @R=0,4cm {
{F \text{-}\mathrm{Isoc} ^{\dag}(Y _0/K)}
\ar[r] ^-{\sp _{Y ^{\dag}  \hookrightarrow U ^{\dag},T _0 , + }}
\ar[d]
&
{F \text{-}\mathrm{Isoc} ^{\dag \dag}  (\PP, T _0, X _0/K)}
\ar[d] ^{| _{\U}}
\\
{F \text{-}\mathrm{Isoc} (Y _0/K)}
\ar[r] ^{\sp _{Y _0 \hookrightarrow \U,+}} _\cong
&
{F \text{-}\mathrm{Isoc} ^{\dag \dag} (\U, Y _0/K),}
}$$
est essentiellement commutatif.
Comme le foncteur de gauche est pleinement fidèle (\cite{kedlaya_full_faithfull}),
il en dérive celle 
de $\sp _{Y ^{\dag}  \hookrightarrow U ^{\dag},T _0 , + }$.
Par \ref{pffidcasgen}, tous les foncteurs de ce diagramme sont ainsi pleinement fidèles.
\end{rema}

\begin{vide}\label{frobb^*}
On reprend les notations et hypothèses de \ref{sp+plssstabiminv}.
On suppose de plus que $P _0$ est séparé et que
$Y  _0 \subset X  _0$ se désingularise localement idéalement.
Par \ref{sp+plssstabiminv} et \ref{theo-ideal->oK},
on dispose alors d'un isomorphisme canonique
\begin{equation}
\label{sp+plssstabiminvfrob}
  \sp _{Y ^{\prime \dag}  \hookrightarrow U ^{\prime \dag},T '_0 , + } (b ^*E)[d_{X '_0/X _0}]
\riso
\R \underline{\Gamma} ^\dag _{X ' _0} f _{T _0} ^!
(\sp _{Y ^{\dag}  \hookrightarrow U ^{\dag},T _0 , + }(E) )
\end{equation}
fonctoriel en $E \in \mathrm{Isoc} ^\dag (Y_0/K)$.
On définit alors l'isomorphisme $\phi '$ via le diagramme ci-après,
\begin{equation}
  \label{sp+plssstabiminvfrob2}
  \xymatrix @R=0,3cm {
  {\sp _{Y ^{\prime \dag}  \hookrightarrow U ^{\prime \dag},T '_0 , + } (F ^* b ^*  E)[d_{X '_0/X _0}]}
  \ar[d] ^-\sim \ar@{.>}[r] _{\phi '}
  &
  {F ^* \sp _{Y ^{\prime \dag}  \hookrightarrow U ^{\prime \dag},T '_0 , + } (b ^*  E)[d_{X '_0/X _0}]}
  \ar[d] ^-\sim
  \\
  {\sp _{Y ^{\prime \dag}  \hookrightarrow U ^{\prime \dag},T '_0 , + } (b ^* F ^* E)[d_{X '_0/X _0}]}
  \ar[d] ^-\sim
  &
     {F ^* \R \underline{\Gamma} ^\dag _{X ' _0} f _{T _0} ^!
(\sp _{Y ^{\dag}  \hookrightarrow U ^{\dag},T _0 , + }( E) )}
\ar[d] ^-\sim
   \\
  {\R \underline{\Gamma} ^\dag _{X ' _0} f _{T _0} ^!
(\sp _{Y ^{\dag}  \hookrightarrow U ^{\dag},T _0 , + }(F ^*  E) )}
 \ar[r] ^-\sim _\phi
 &
  {\R \underline{\Gamma} ^\dag _{X ' _0} f _{T _0} ^! F ^*
(\sp _{Y ^{\dag}  \hookrightarrow U ^{\dag},T _0 , + }( E) ),}
}
\end{equation}
où l'isomorphisme du bas dérive de \ref{theo-ideal->oK}
et
où les isomorphismes à droite en haut et à gauche en bas découlent de \ref{sp+plssstabiminvfrob} appliqué
respectivement à $E$ et $F ^* E$.
On vérifie ensuite que $\phi '[-d_{X '_0/X _0}] $
est l'isomorphisme canonique de commutation à Frobenius
de $\sp _{Y ^{\prime \dag}  \hookrightarrow U ^{\prime \dag},T '_0 , + } (b ^*  E)$ (voir \ref{isocanfrob}).
Enfin, par construction, l'isomorphisme
\ref{sp+plssstabiminvfrob} commute à Frobenius.
\end{vide}

\begin{vide}\label{isocanfrobf+}
Reprenons à présent les notations et hypothèses de \ref{theo-ideal->oKlemm}.
Si l'isomorphisme canonique de commutation à Frobenius de
$\sp _{Y ^{\prime \dag}  \hookrightarrow U ^{\prime \dag}, f ^{-1} T _0 , +} $ en $E '$ existe
et est noté $\phi '$ alors, la flèche $\phi$ construite via le diagramme
  \begin{equation} \notag 
  \xymatrix @R=0,4cm {
{\sp _{Y ^{\dag}  \hookrightarrow U ^{\dag},T _0 , +} (F ^* b _*   E ')}
\ar[d] ^-\sim \ar@{.>}[r] _\phi
&
{F ^* \sp _{Y ^{\dag}  \hookrightarrow U ^{\dag},T _0 , +} (b _* E ')}
\ar[d] ^-\sim
\\
{\sp _{Y ^{\dag}  \hookrightarrow U ^{\dag},T _0 , +} ( b _*  F ^* E ')}
\ar[d] ^-\sim
&
{F ^* f _{T _0, +} (\sp _{Y ^{\prime \dag}  \hookrightarrow U ^{\prime \dag}, f ^{-1} T _0 , +} (E'))}
\ar[d] ^-\sim
\\
{f _{T _0, +} (\sp _{Y ^{\prime \dag}  \hookrightarrow U ^{\prime \dag}, f ^{-1} T _0 , +} (F ^* E'))}
\ar[r] ^-\sim _{\phi '}
&
{f _{T _0, +} F ^*  (\sp _{Y ^{\prime \dag}  \hookrightarrow U ^{\prime \dag}, f ^{-1} T _0 , +} (E')),}
}
\end{equation}
où la flèche de droite du haut et celle de gauche du bas résultent de \ref{theo-ideal->oK-lemmmorph}
tandis que les deux autres verticales sont les isomorphismes de commutation à Frobenius des images directes,
est
l'isomorphisme canonique de commutation à Frobenius de
$\sp _{Y ^{\dag}  \hookrightarrow U ^{\dag},T _0 , +} (b _* E ')$.
Par construction, l'isomorphisme  
$\sp _{Y ^{\dag}  \hookrightarrow U ^{\dag},T _0 , +} (b _* E ')
\riso
f _{T _0  {} _+ } (\sp _{Y ^{\prime \dag}  \hookrightarrow U ^{\prime \dag}, T '_0 , +} (E'))$
canonique \linebreak commute à Frobenius.
\end{vide}

\begin{vide}\label{frobb^*2}
Soit le diagramme de $\V$-schémas formels faibles lisses
\begin{equation}\label{frobb^*2-diagm1}
  \xymatrix @R=0,3cm {
{Y ^{'\dag}}
\ar[r] ^{v'} \ar[d] ^b
&
{U ^{'\dag} }
\ar[r] ^{j'} \ar[d] ^g
&
{P ^{'\dag} }
\ar[d] ^f
\\
{Y ^{\dag}}
\ar[r] ^{v}
&
{U ^{\dag} }
\ar[r] ^{j}
&
{P ^{\dag} ,}
}
\hfill
\xymatrix @R=0,3cm {
{Y ^{''\dag}}
\ar[r] ^{v''} \ar[d] ^{b'}
&
{U ^{''\dag} }
\ar[r] ^{j''} \ar[d] ^{g '}
&
{P ^{''\dag} }
\ar[d] ^{f'}
\\
{Y ^{\prime \dag}}
\ar[r] ^{v'}
&
{U ^{\prime \dag} }
\ar[r] ^{j''}
&
{P ^{\prime \dag} ,}
}
\end{equation}
  où $f$ et $f '$ sont propres, $j$ est une immersion ouverte, les carrés de droite respectifs sont cartésiens,
  $v$, $v'$ et $v''$ sont des immersions fermées, $b$ et $b'$ sont finis et étales et $Y ^{\dag}$ est affine.
  On suppose de plus que $T _0 := P _0\setminus U _0$, $ T ' _0:= P '_0\setminus U '_0$
  et $ T '' _0:= P '' _0\setminus U '' _0$ sont les supports de diviseurs.
  Enfin, en notant $X _0$ (resp. $X' _0$ et $X '' _0$) l'adhérence de $Y _0$ (resp. $Y '_0$ et $Y ''_0$),
  on suppose que $Y _0\subset X _0$ se désingularise localement idéalement.

En notant $\sp _+ =\sp _{Y ^{\dag}  \hookrightarrow U ^{\dag},T _0 , +}$,
pour tout $F$-isocristal surconvergent $E$ sur $Y _0$,
le carré de droite du diagramme
\begin{equation}
  \label{theo-ideal->oK-lemmmorphfrobtr}
\xymatrix @R=0,3cm @C=0,5cm {
{\sp _+  (b ^* E )}
\ar@{=}[d]
\ar[r]
&
{\sp _+  (b '_* b ^{\prime *} b ^* E )}
\ar[r] \ar[d] ^-\sim
&
{\sp _+  (b ^* E )}
\ar@{=}[d]
\\
{\sp _+  (b ^* E )}
\ar@{.>}[r]
&
{f ' _{T '_0, +}  \R \underline{\Gamma} ^\dag _{X '' _0} f ^{\prime !} _{T' _0} \sp _+  (b ^* E ) }
\ar[r]
&
{\sp _+  (b ^* E ).}
}
\end{equation}
où l'isomorphisme vertical du milieu est le composé de
\ref{isocanfrobf+} et de \ref{sp+plssstabiminvfrob}
et où les morphismes
horizontaux de droite se déduisent des morphismes traces canoniques.
En effet, grâce à \ref{rho=rho}, celui-ci l'est au-dessus $\U'$.
\end{vide}

\begin{lemm}\label{lem-DRfrob}
  Avec les notations et hypothèses de \ref{notation-surcoh-gen}, on suppose
$P _0$ propre et $X _0$ lisse.
  Pour tout $F$-isocristal surconvergent $E$ sur $Y  _0$,
  on bénéficie d'un isomorphisme canonique commutant à Frobenius
$$H ^\bullet _{DR} ( \sp _{Y ^{\dag}  \hookrightarrow U ^{\dag},T _0  +} (E)/K)
\riso H ^\bullet _{\mathrm{rig}} ( E/K).$$
\end{lemm}
\begin{proof}
Notons $\sp _+ := \sp _{Y ^{\dag}  \hookrightarrow U ^{\dag},T _0  +}$
et $\smash{\widetilde{\sp}} _+ := \sp _{X _0 \hookrightarrow \PP ,T _0 , +}$.
Comme $X _0$ est lisse,
d'après \ref{conjvalF}, on dispose de l'isomorphisme ca-\linebreak nonique
$\sp _{+} ( E)
\riso
\smash{\widetilde{\sp}} _+ (   E)$, celui-ci commutant à Frobenius.
Par \cite[4.3.6.3]{Beintro2} et
par descente cohomologique par un recouvrement ouvert,
on vérifie
$  H ^\bullet _{DR} ( \smash{\widetilde{\sp}} _+ ( E)/K)
\riso
H ^\bullet _{\mathrm{rig}} (  E/K).$
\hfill \hfill \qed \end{proof}

\begin{prop}
  \label{prop-DRfrob}
  Avec les notations \ref{notation-surcoh-gen}, on suppose que
$Y _0 \hookrightarrow X_0$ se désingularise localement idéalement et
$P _0$ propre.

  Pour tout $F$-isocristal surconvergent $E$ sur $Y  _0$,
  on bénéficie d'un isomorphisme canonique
$$H ^\bullet _{DR} ( \sp _{Y ^{\dag}  \hookrightarrow U ^{\dag}, T _0 , +} (E)/K)
\riso H ^\bullet _{\mathrm{rig}} ( E/K)$$
  commutant à Frobenius.
\end{prop}
\begin{proof}
On reprend la construction et les notations
du premier paragraphe de la preuve de \ref{theo-ideal->oK} concernant $f$, $b$ etc.
De plus, on note, pour tous $r \geq 1$ et $i=0,\dots, r$, $P ^{\prime \dag  r}:=
P ^{\prime \dag  } \times _{P ^\dag}  P ^{\prime \dag } \times \cdots \times _{P ^{\dag}}P ^{\prime \dag  }$ ($r$-fois),
$f ^{i,r}$ : $P ^{\prime \dag  r+1} \rightarrow P ^{\prime \dag  r}$ la $i$-ième projection,
$f ^{r}$ : $P ^{\prime \dag  r} \rightarrow P ^{\dag  }$ la projection
et $f ^{\prime r}$ : $P ^{\prime \dag  r} \rightarrow P ^{\prime \dag  }$ la $0$-ième projection.
De même, en remplaçant $P$ par $Y$ et $f$ par $b$. On note
$X ^{\prime r} _{0}$, l'adhérence de $Y ^{\prime r} _{0}$ dans $P ^{\prime r} _{0}$,
$b ^{i,r!}  := \R \underline{\Gamma} ^\dag _{X ^{'r} _0} f  ^{i,r!}$ (on omet d'indiquer le diviseur)
et
$b ^{i,r} _+ := f  ^{i,r} _+$, et aussi en remplaçant {\og $i,r$ \fg} par
{\og $r$ \fg} ou {\og $\prime r$ \fg} ou rien.
Notons $\sp _+ := \sp _{Y ^{\dag}  \hookrightarrow U ^{\dag},T _0  +}$
et
$\sp ^\prime _+ := \sp _{Y ^{\prime \dag}  \hookrightarrow U ^{\prime \dag}, f ^{-1}( T _0 ) +}$.
Pour tout $F$-isocristal $E '$ surconvergent sur $Y' _0$, on définit
${H ^\bullet _{\mathrm{rig}} (b _* E'/K)} \riso
{H ^\bullet _{DR} ( b _+\sp _{+} ^\prime (E')/K)}$
via le diagramme
\begin{equation}
  \label{prop-DRfrob2}
  \xymatrix @R=0,3cm {
  {H ^\bullet _{\mathrm{rig}} (b _* E'/K)}
  \ar@{.>}[d]
  \ar[r] ^-\sim
  &
  {H ^\bullet _{\mathrm{rig}} ( E'/K)}
  \ar@{.>}[d]
  \ar[d] ^-\sim
\\
{H ^\bullet _{DR} ( b _+\sp _{+} ^\prime (E')/K)}
  \ar[r] ^-\sim
&
  {H ^\bullet _{DR} ( \sp _{+} ^\prime (E')/K)}
  }
\end{equation}
Il en découle le carré commutatif du haut de
\begin{equation}
  \label{prop-DRfrob3}
\xymatrix @R=0,3cm  {
{H ^\bullet _{\mathrm{rig}} (b _* ^r b ^{r*} E )}
\ar[d]
\ar[r]
&
{H ^\bullet _{\mathrm{rig}} (b ^r _*  b^{i,r} _* b ^{i,r*} b ^{r*} E /K)}
\ar[d] ^-\sim
\\
{H ^\bullet _{DR} (b _+\sp _{+} ^\prime b _* ^{\prime r} b ^{r*} E / K)}
\ar[r]
\ar[d]
&
{H ^\bullet _{DR} (b _+\sp _{+} ^\prime b ^{\prime r} _*  b^{i,r} _* b ^{i,r*} b ^{r*} E /K)}
\ar[d]
\\
{H ^\bullet _{DR} (b _+ b _+ ^{\prime r} b ^{r!}\sp _{+}  E / K)}
\ar[r]
&
{H ^\bullet _{DR} (b _+ b ^{\prime r} _+  b^{i,r} _+ b ^{i,r!} b ^{r!}\sp _{+} E /K)}
}
\end{equation}
dont les morphismes horizontaux sont les morphismes d'adjonction.
La commutativité du carré du bas se prouve par fonctorialité et via
le carré de gauche \ref{theo-ideal->oK-lemmmorphfrobtr}.
On conclut par descente cohomologique finie et étale appliquée
  au recouvrement $ Y ^{\prime \dag} \rightarrow Y ^\dag$
  (on utilise la résolution \cite[2.19]{milne}).
\hfill \hfill \qed \end{proof}

\begin{vide}\label{boxplus}
    Soient $P ^\dag$, $P ^{\prime \dag}$ deux $\V$-schémas formels faibles lisses et séparés,
    $T _0$ (resp. $T '_0$) un diviseur de $P _0$ (resp. $P '_0$), $U ^\dag$ (resp. $U ^{\prime \dag }$)
    l'ouvert de $P ^\dag $ (resp. $P ^{\prime \dag}$)
    complémentaire de $T _0$ (resp. $T '_0$), $j$ : $U ^\dag \hookrightarrow P ^\dag$
    (resp. $j'$ : $U ^{\prime \dag} \hookrightarrow P ^{\prime \dag}$)
    l'immersion ouverte et $v$ : $Y ^\dag \hookrightarrow U ^\dag$
    (resp. $v '$ : $Y ^{\prime \dag} \hookrightarrow U ^{\prime \dag}$)
    une immersion fermée de $\V$-schémas formels faibles.
On suppose en outre $Y ^\dag$ et $Y ^{\prime \dag}$ affines et lisses et on désigne par $X _0$ (resp. $X '_0$)
l'adhérence schématique de $Y _0$ (resp. $Y ' _0$) dans $P _0$ (resp. $P ' _0$).

On note $P ^{\prime \prime \dag} := P ^\dag \times P ^{\prime \dag}$,
$U ^{\prime \prime \dag} := U ^\dag \times U ^{\prime \dag}$,
$T ''_0$ le diviseur réduit de $P '' _0$ d'espace topologique $P'' _0 \setminus U ''_0$,
$Y ^{\prime \prime \dag} := Y ^\dag \times Y ^{\prime \dag}$,
$f _1$ : $P ^{\prime \prime \dag} \rightarrow  P ^{\dag}$,
$f _2$ : $P ^{\prime \prime \dag} \rightarrow  P ^{\prime \dag}$,
$g _1$ : $U ^{\prime \prime \dag} \rightarrow  U ^{\dag}$,
$g _2$ : $U ^{\prime \prime \dag} \rightarrow  U ^{\prime \dag}$,
$b _1$ : $Y ^{\prime \prime \dag} \rightarrow  Y ^{\dag}$,
$b _2$ : $Y ^{\prime \prime \dag} \rightarrow  Y ^{\prime \dag}$.

Si $E \in \mathrm{Isoc} ^\dag ( Y _0 /K)$ et
$E' \in \mathrm{Isoc} ^\dag ( Y '_0 /K)$, on pose
$E \boxtimes E ' := b _1 ^* (E ) \otimes  b _2 ^* (E ')$,
le produit tensoriel se calculant dans $\mathrm{Isoc} ^\dag ( Y ''_0 /K)$.

\end{vide}

\begin{prop}\label{boxplusprop}
  Avec les notations \ref{boxplus}, on suppose que
  $Y _0 \subset X _0$ et $Y ' _0\subset X '_0$ se désingularisent localement idéalement.

  On dispose alors d'un isomorphisme canonique
$$ \sp _{Y ^{\prime \prime \dag} \hookrightarrow U ^{\prime \prime \dag} , T'' _0 +}
(E \boxtimes E ') \riso
\sp _{Y ^{\dag} \hookrightarrow U ^{\dag} , T _0 +} (E)
\underset{\O _{\PP''} (\hdag T_0 '') _\Q}{\smash{\overset{\L}{\boxtimes}}^\dag}
\sp _{Y ^{\prime \dag} \hookrightarrow U ^{\prime \dag} , T '_0 +} (E').$$
\end{prop}
\begin{proof}
Considérons les diagrammes commutatifs
\begin{equation} \notag
  \xymatrix @C=0,5cm @R=0,3cm {
{\smash{\widetilde{Y}} ^{\dag}}
\ar[r] ^{\tilde{v}} \ar[d] ^{b}
&
{\smash{\widetilde{U}} ^{\dag} }
\ar[r] ^{\tilde{j} } \ar[d] ^{g}
&
{\smash{\widetilde{P}} ^{\dag} }
\ar[d] ^{f }
\\
{Y  ^{\dag}}
\ar[r] ^{v}
&
{U ^{\dag} }
\ar[r] ^{j  }
&
{P ^{\dag} ,}
}
\
\xymatrix @C=0,5cm @R=0,3cm {
{\smash{\widetilde{Y}} ^{\prime \dag}}
\ar[r] ^{\smash{\tilde{v}} '} \ar[d] ^{b' }
&
{\smash{\widetilde{U}} ^{\prime \dag} }
\ar[r] ^{\smash{\tilde{j}} ' } \ar[d] ^{g '}
&
{\smash{\widetilde{P}} ^{\prime \dag} }
\ar[d] ^{f '}
\\
{Y  ^{\prime \dag}}
\ar[r] ^{v '}
&
{U ^{\prime \dag} }
\ar[r] ^{j ' }
&
{P ^{\prime \dag} ,}
}
\
\xymatrix @C=0,5cm @R=0,3cm {
{\smash{\widetilde{Y}} ^{\prime \prime \dag}}
\ar[r] ^{\smash{\tilde{v}}  ^{\prime \prime} } \ar[d] ^{b^{\prime \prime }}
&
{\smash{\widetilde{U}} ^{\prime \prime \dag} }
\ar[r] ^{\smash{\tilde{j}} ^{\prime \prime} } \ar[d] ^{g ^{\prime \prime }}
&
{\smash{\widetilde{P}} ^{\prime \prime \dag} }
\ar[d] ^{f ^{\prime \prime }}
\\
{Y  ^{\prime \prime \dag}}
\ar[r] ^{v ^{\prime \prime }}
&
{U ^{\prime \prime \dag} }
\ar[r] ^{j ^{\prime \prime }}
&
{P ^{\prime \prime \dag} ,}
}
\end{equation}
où les deux de gauche sont construits de manière analogue à celui du début de la preuve
de \ref{theo-ideal->oK},
celui de droite s'en déduisant en posant
$\smash{\widetilde{Y}} ^{\prime \prime \dag} := \smash{\widetilde{Y}} ^{\dag} \times \smash{\widetilde{Y}} ^{\prime \dag}$,
$\smash{\widetilde{U}} ^{\prime \prime \dag} := \smash{\widetilde{U}} ^{\dag} \times \smash{\widetilde{U}} ^{\prime \dag}$,
$\smash{\widetilde{P}} ^{\prime \prime \dag} := \smash{\widetilde{P}} ^{\dag} \times \smash{\widetilde{P}} ^{\prime \dag}$,
$b^{\prime \prime } = b \times b^{\prime }$,
$g^{\prime \prime } = g \times g^{\prime }$,
$f^{\prime \prime } = f \times f^{\prime }$,
$j  ^{\prime \prime } = j  \times  j  ^{\prime }$,
$v  ^{\prime \prime } = v  \times  v  ^{\prime }$,
$\smash{\tilde{j}}  ^{\prime \prime } = \smash{\tilde{j}}  \times  \smash{\tilde{j}}  ^{\prime }$,
$\smash{\tilde{v}}  ^{\prime \prime } = \smash{\tilde{v}}  \times  \smash{\tilde{v}}  ^{\prime }$.
On note $\smash{\widetilde{X}} _0$ (resp. $\smash{\widetilde{X}} '_0$)
l'adhérence de $\smash{\widetilde{Y}} _0$ (resp. $\smash{\widetilde{Y}} '_0$) dans
celle de $\smash{\widetilde{P}} _0$ (resp. $\smash{\widetilde{P}} '_0$),
$\smash{\widetilde{X}}'' _0 := \smash{\widetilde{X}} _0 \times \smash{\widetilde{X}}' _0$,
$\smash{\tilde{f}} _1 $ : $\smash{\widetilde{P}} ^{\prime \prime \dag} \rightarrow  \smash{\widetilde{P}} ^{\dag}$ et
$\smash{\tilde{f}} _2$ : $\smash{\widetilde{P}} ^{\prime \prime \dag} \rightarrow  \smash{\widetilde{P}} ^{\prime \dag}$,
$\smash{\tilde{b}} _1$ : $\smash{\widetilde{Y}} ^{\prime \prime \dag} \rightarrow  \smash{\widetilde{Y}} ^{\dag}$,
$\smash{\tilde{b}} _2$ : $\smash{\widetilde{Y}} ^{\prime \prime \dag} \rightarrow  \smash{\widetilde{Y}} ^{\prime \dag}$.
De plus, on pose $\sp _{_+} := \sp _{Y ^{\dag} \hookrightarrow U ^{\dag} , T _0 { }_+}$,
$\sp ' _{_+}:= \sp _{Y ^{\prime \dag} \hookrightarrow U ^{\prime \dag} , T '_0 { }_+}$,
$\sp '' _{_+}:= \sp _{Y ^{\prime \prime \dag} \hookrightarrow U ^{\prime \prime \dag} , T ''_0 { }_+}$,
$\smash{\widetilde{\sp }} _{_+} :=
\sp _{\smash{\widetilde{Y}} ^{\dag} \hookrightarrow \smash{\widetilde{U}} ^{\dag},f ^{-1} (T _0) { }_+}$,
$\smash{\widetilde{\sp }} ^\prime _{_+} :=
\sp _{\smash{\widetilde{Y}} ^{\prime \dag} \hookrightarrow \smash{\widetilde{U}} ^{\prime \dag} ,
f ^{\prime -1} (T '_0) { }_+}$,
$\smash{\widetilde{\sp }} ^{\prime \prime} _{_+} :=
\sp _{\smash{\widetilde{Y}} ^{\prime \prime \dag} \hookrightarrow \smash{\widetilde{U}} ^{\prime \prime \dag} ,
f ^{\prime \prime -1} (T '' _0) { }_+}$.

Par adjonction puis via \ref{theo-ideal->oK-lemmmorph},
on dispose des morphismes
\begin{equation}
  \label{boxpluspropeq1}
\sp '' _+  ( E \boxtimes E')
\rightarrow
\sp '' _+ b '' _* b ^{\prime \prime *} ( E \boxtimes E')
\riso
f ^{\prime \prime} _+ (
\smash{\widetilde{\sp }} '' _+ b ^{\prime \prime *}
( E \boxtimes E')) .
\end{equation}
Considérons l'isomorphisme composé suivant
\begin{gather}\notag
\smash{\widetilde{\sp }} ^{\prime \prime} _+ b ^{\prime \prime *}
( E \boxtimes E')
=
\smash{\widetilde{\sp }} ^{\prime \prime} _+
(b ^{\prime \prime *}  b _1 ^* E \otimes b ^{\prime \prime *} b _2 ^{\prime *} E')
\riso
\smash{\widetilde{\sp }} ^{\prime \prime} _+
(\smash{\tilde{b}}_1 ^* b ^* E \otimes \smash{\tilde{b}}_2 ^* b ^{\prime *} E')
 \\ \label{boxpluspropgath0}
\riso
(\smash{\widetilde{\sp }} ^{\prime \prime} _+ (\smash{\tilde{b}}_1 ^* b ^* E) )
\smash{\overset{\L}{\otimes}} ^\dag
_{\O _{\smash{\widetilde{\PP}}''} (\hdag \smash{\widetilde{T}}_0'') _\Q}
 \smash{\widetilde{\sp }} ^{\prime \prime} _+ (\smash{\tilde{b}}_2 ^* b ^{\prime *} E')
[d _{\smash{\widetilde{X}} '' _0/ \smash{\widetilde{P}} ''_0}]
 \\
 \label{boxpluspropgath}
 \riso
\smash{\tilde{f}}  _1  ^!
(\smash{\widetilde{\sp}} _+ ( b ^* E) )
\smash{\overset{\L}{\otimes}} ^\dag
_{\O _{\smash{\widetilde{\PP}}''} (\hdag \smash{\widetilde{T}}_0'') _\Q}
\smash{\tilde{f}} _2 ^!( \smash{\widetilde{\sp}} ' _+ ( b ^{\prime *} E') )
[-d _{\smash{\widetilde{P}} ''_0}]
 \\
 {\label{boxpluspropgath1}
 \riso
\smash{\tilde{f}}  _1  ^!
\R \underline{\Gamma} ^\dag _{\smash{\widetilde{X}} _0 } f ^!  \sp _+ (E)
\smash{\overset{\L}{\otimes}} ^\dag
_{\O _{\smash{\widetilde{\PP}}''} (\hdag \smash{\widetilde{T}}_0'') _\Q}
\smash{\tilde{f}} _2 ^!
\R \underline{\Gamma} ^\dag _{\smash{\widetilde{X}} '_0}
f ^{\prime !} \sp ' _+ (E'))
[-d _{\smash{\widetilde{P}} ^{\prime \prime } _0}]}
 \\  \notag
\riso
\R \underline{\Gamma} ^\dag _{\smash{\widetilde{X}} _0 \times \smash{\widetilde{P}} '_0}
\smash{\tilde{f}}  _1  ^!  f ^!  \sp _+ (E)
\smash{\overset{\L}{\otimes}} ^\dag
_{\O _{\smash{\widetilde{\PP}}''} (\hdag \smash{\widetilde{T}}_0'') _\Q}
\R \underline{\Gamma} ^\dag _{\smash{\widetilde{P}} _0\times \smash{\widetilde{X}} '_0}
\smash{\tilde{f}} _2 ^!  f ^{\prime !}
\sp ' _+ (E'))
[-d _{\smash{\widetilde{P}} ^{\prime \prime } _0}]
 \\
 \notag
\riso \R \underline{\Gamma} ^\dag _{\smash{\widetilde{X}} _0 ^{\prime \prime}} f ^{\prime \prime !}
(  f ^! _1 \sp _+ (E)
\smash{\overset{\L}{\otimes}} ^\dag _{\O _{\PP''} (\hdag T _0'') _\Q}
 f ^! _2 \sp ' _+ (E'))
 [-d _{P '' _0}]
 \\
\notag =
\R \underline{\Gamma} ^\dag _{\smash{\widetilde{X}} _0 ^{\prime \prime}} f ^{\prime \prime !}
  (\sp _+ (E)
\smash{\overset{\L}{\boxtimes}} ^\dag _{\O _{\PP''} (\hdag T_0'') _\Q}
\sp ' _+  (E'))
\end{gather}
où \ref{boxpluspropgath0} se déduit de \ref{lissestableotimes},
\ref{boxpluspropgath} résulte de \cite{caro_unite} et
\ref{boxpluspropgath1} découle de
\ref{sp+plssstabiminvfrob}.
En lui appliquant $f ^{\prime \prime} _+$ cela donne :
\begin{equation}\label{boxpluspropeq2}
  f ^{\prime \prime} _+ (
\smash{\widetilde{\sp }} '' _+ b ^{\prime \prime *}
( E \boxtimes E'))
\riso
f ^{\prime \prime} _+
\R \underline{\Gamma} ^\dag _{\smash{\widetilde{X}} _0 ^{\prime \prime}} f ^{\prime \prime !}(
  \sp _+ (E)
\smash{\overset{\L}{\boxtimes}} ^\dag _{\O _{\PP''} (\hdag T_0 '') _\Q}
\sp ' _+  (E') ).
\end{equation}
Via le morphisme d'adjonction
$f ^{\prime \prime} _+
\R \underline{\Gamma} ^\dag _{\smash{\widetilde{X}} _0 ^{\prime \prime}} f ^{\prime \prime !} \rightarrow id$,
en composant \ref{boxpluspropeq1} et \ref{boxpluspropeq2}, on obtient
\begin{equation}\label{boxpluspropeq3}
 \sp '' _+  ( E \boxtimes E')
\rightarrow
  \sp _+ (E) \smash{\overset{\L}{\boxtimes}} ^\dag _{\O _{\PP''} (\hdag T_0 '') _\Q}
\sp ' _+  (E') .
\end{equation}
Comme \ref{boxpluspropeq3} est un morphisme de
$\D ^\dag _{\PP''} (\hdag T_0 '') _\Q$-modules cohérents, il suffit de
le vérifier au dessus de $\U \times \U'$, ce qui découle de \ref{lissestableotimes}
et de \ref{spyu+rest}.

\hfill \hfill \qed \end{proof}

\begin{coro}\label{stbintn}
  Avec les notations et hypothèses de \ref{boxplusprop}, pour tous
  $E _1, E _2 \in \mathrm{Isoc} ^\dag ( Y _0 /K)$,
en notant $\sp _+ =\sp _{Y ^{\dag} \hookrightarrow U ^{\dag} , T _0 +}$,
  on dispose d'un isomorphisme canonique :
  $$\sp _{+} (E _1 \otimes E _2)
  \riso
  \sp _{ +} (E _1)
\smash{\overset{\L}{\otimes}} ^\dag _{\O _{\PP} (\hdag T) _\Q}
\sp _{ +} (E_2)[d _{X /P}].$$
\end{coro}
\begin{proof}
D'après \ref{boxplusprop} appliqué au cas particulier
où les objets avec des primes sont égaux à ceux sans primes,
on dispose de l'isomorphisme :
$$ \sp _{Y ^{\prime \prime \dag} \hookrightarrow U ^{\prime \prime \dag} , T'' _0 +}
(E _1 \boxtimes E _2) \riso
\sp _{ +} (E _1)
\smash{\overset{\L}{\boxtimes}} ^\dag _{\O _{\PP ''} (\hdag T'') _\Q}
\sp _{ +} (E _2).$$
Or, en notant $\delta$ l'immersion diagonale $\PP \hookrightarrow \PP \times \PP = \PP ''$,
on obtient :
$
\delta ^!   (\sp _+ (E _1)
\smash{\overset{\L}{\boxtimes}} ^\dag _{\O _{\PP ''} (\hdag T'') _\Q}
\sp _+ (E _2))
\riso
\sp _+ (E _1)
\smash{\overset{\L}{\otimes}} ^\dag _{\O _{\PP } (\hdag T) _\Q}
\sp _+ (E _2)[- d _{P _0}].$
De plus, via \ref{sp+plssstabiminvfrob},
$$\sp _+ ( E _1 \otimes E _2) [- d _{X _0}]
\riso
 \delta ^! (\sp _{Y ^{\prime \prime \dag} \hookrightarrow U ^{\prime \prime \dag} , T'' _0 +}
(E _1 \boxtimes E _2) ).$$
D'où le résultat.
\hfill \hfill \qed \end{proof}

\begin{vide}\label{sp+y0}
  On garde les notations et hypothèses de \ref{notation-surcoh-gen}.
On suppose que $P ^\dag$ est propre et que $Y _0\subset X _0$ se désingularise localement idéalement.

Le foncteur $\sp _{Y ^{\dag} \hookrightarrow U ^{\dag} , T _0 +}$ ne dépend pas des choix
de $P ^\dag$ et $v $ : $Y ^{\dag} \hookrightarrow U ^{\dag}$ choisis tels,
que $Y _0 \subset X _0$ se désingularise localement idéalement, i.e., celui-ci
est compatible avec les équivalences canoniques entre les catégories de la forme
$F\text{-}\mathrm{Isoc} ^{\dag \dag}(\PP, T _0, X _0/K)$.
En effet, 
cela découle 
de \ref{sp+plssstabiminvfrob}.
On obtient alors un foncteur,
$F\text{-}\mathrm{Isoc} ^\dag( Y _0/K)\rightarrow F\text{-}\mathrm{Isoc} ^{\dag \dag}(Y _0/K) $,
bien défini à isomorphisme canonique près et que l'on notera $\sp _{Y _0 +}$.
\end{vide}

\begin{theo}\label{coro-ideal->oK}
On garde les notations et hypothèses de \ref{sp+y0}.

Les foncteurs
 $\sp _{Y _0 +}$ :
 $F\text{-}\mathrm{Isoc} ^\dag( Y _0/K)\rightarrow F\text{-}\mathrm{Isoc} ^{\dag \dag}(Y _0/K) $
 et $\rho _{Y_0}$ (voir \ref{isosurcohpf}) sont des équivalences quasi-inverses.
\end{theo}
\begin{proof}
Cela résulte de \ref{rhocv=cv}, de \ref{remspcv=cv}.

\hfill \hfill \qed \end{proof}

\begin{prop}\label{dense-ideal}
  Avec les notations \ref{notation-surcoh-gen}, on suppose $X _0$ intègre et $P ^\dag$ propre.

  Il existe alors un diviseur $\smash{\widetilde{T}} _0$ de
  $P _0$  contenant $T _0$ tel que $\widetilde{Y} _0:= (P _0 \setminus \smash{\widetilde{T}} _0  ) \cap X _0$
  soit affine, dense dans $Y _0$ et
  l'immersion ouverte $\widetilde{Y} _0 \hookrightarrow X _0$ se désingularise idéalement.
\end{prop}
\begin{proof}Grâce au théorème de désingularisation de de Jong (\cite{dejong}) et quitte à remplacer $Y _0$ par
un ouvert affine et dense,
il existe un morphisme projectif et surjectif $a _0$ : $X'_0 \rightarrow X _0$, qui se décompose
  en une immersion fermée $X '_0 \hookrightarrow \P ^r _{X _0}$ suivie de la projection
  canonique $\P ^r _{X _0} \rightarrow X _0$, tel que
  \begin{enumerate}
    \item $X ' _0$ est intègre et lisse ;
    \item le morphisme $b _0$ : $Y ' _0 := a_0 ^{-1} (X _0) \rightarrow Y _0$ induit par $a_0$ est
    fini et étale.
    \end{enumerate}
Il reste à vérifier que quitte à rétrécir à nouveau $Y _0$, la propriété $3$ de
\ref{def-notation-surcoh-gen} est validée.
Notons $v ' _0$ l'immersion fermée $Y '_0 \hookrightarrow \P ^r _{Y _0}$.
  Si $x _0 ,\dots, x _r$ sont les coordonnées projectives de
  $\P ^r _k$, on note, pour tout entier $\alpha \in \{ 0,\dots,r\}$,
  $D _\alpha$ le diviseur de $\P ^r _k$ défini par l'équation $x _\alpha = 0$ et
  $D _{\alpha, Y _0}:=D _{\alpha} \times _{\P ^r _k }\P ^r _{Y _0} $.
  Comme l'intersection des diviseurs $D _{\alpha, Y _0}$ est vide,
  il existe un entier $\alpha _0$ tel que $v ' _0 (Y ' _0)$ ne soit pas inclus dans
  $D _{\alpha _0, Y _0}$. Comme $Y_0 '$ est intègre, on obtient
  $\dim Y ' _0 \cap D _{\alpha _0, Y _0} < \dim Y ' _0$.
  Via \cite[5.4.2]{EGAIV2}, la finitude de $b_0$ implique l'égalité
  $\dim b _0 ( Y ' _0 \cap D _{\alpha _0, Y _0} ) = \dim Y ' _0 \cap D _{\alpha _0, Y _0} $.
  Comme $b_0$ est en outre surjectif, on a aussi $\dim Y ' _0 = \dim Y _0$ et donc
  $\dim b_0 ( Y ' _0 \cap D _{\alpha _0, Y _0} ) < \dim Y _0$.
  Il résulte de cette dernière inégalité qu'il existe un diviseur ${\widetilde{T}}_0$ de $P _0$ tel
  que l'ouvert $(P _0 \setminus {\widetilde{T}} _0) \cap X _0$ de $X _0$ soit affine et inclus dans
  $Y _0\setminus b_0 ( Y ' _0 \cap D _{\alpha _0, Y _0} ) $ (le fait que l'on peut le choisir affine
  résulte de \cite[6.3.1]{caro}).

Posons ${\widetilde{T}} _{X _0} = {\widetilde{T}} _0\cap X_0 $ et $ {\widetilde{T}} _{X_0 '} := a _0^{-1} ({\widetilde{T}} _{X _0})$.
L'inclusion $b_0 ( Y ' _0 \cap D _{\alpha _0, Y _0} ) \subset {\widetilde{T}} _{X _0}$
(resp. $X _0 \setminus {\widetilde{T}} _{X _0} \subset  Y _0$) implique
alors $Y ' _0 \cap D _{\alpha _0, Y _0} \subset {\widetilde{T}} _{X_0 '}$
(resp. $X '_0 \setminus {\widetilde{T}} '_{X _0} \subset  Y '_0$).
Il en résulte la factorisation $ X_0 ' \setminus {\widetilde{T}} _{X_0 '}
\hookrightarrow (\P _k ^r \setminus D _{\alpha _0} )\times  (X _0 \setminus {\widetilde{T}} _{X _0})$.
Celle-ci se relève en un morphisme de $\V$-schémas formels faibles affines et lisses.
Il en découle que le morphisme canonique
$X_0 ' \setminus {\widetilde{T}} _{X_0 '} \rightarrow \P ^r _{X_0  \setminus {\widetilde{T}} _{X _0} }$
se relève en un morphisme de $\V$-schémas formels faibles lisses.
L'immersion ouverte $X _0 \setminus {\widetilde{T}} _{X _0} \hookrightarrow X_0$ se désingularise donc idéalement.
\hfill \hfill \qed \end{proof}

\begin{defi}\label{defmodide}
  Soit $Y _0$ un $k$-schéma affine et lisse. On dit que $Y _0$ possède un {\it modèle idéal}
  s'il existe un $\V$-schéma formel faible propre et lisse $P ^\dag$, un diviseur $T _0$ de $P _0$,
 un sous-schéma fermé $X _0$ de $P _0$ tels que $Y _0 \riso  X _0 \setminus T _0$ et
 l'immersion ouverte $Y _0 \hookrightarrow X _0$ se désingularise localement idéalement.

 D'après \ref{coro-ideal->oK}, si $Y _0$ est un $k$-schéma affine, lisse et possédant un modèle idéal,
 le foncteur
$\sp _{Y _0 +}$ :
$F\text{-}\mathrm{Isoc} ^{ \dag}(Y _0/K) \rightarrow F\text{-}\mathrm{Isoc} ^{\dag \dag}( Y _0/K)$ 
est une équivalence de catégorie.
\end{defi}

\begin{theo}\label{densemodide}
  Soit $Y _0$ un $k$-schéma lisse.
  Il existe $\widetilde{Y} _0$ un ouvert affine, dense dans $Y _0$, possédant un modèle idéal et,
  en particulier, 
  le foncteur canonique
$\sp _{Y _0 +}$ :
$F\text{-}\mathrm{Isoc} ^{ \dag}(Y _0/K) \rightarrow F\text{-}\mathrm{Isoc} ^{\dag \dag}( Y _0/K)$
 est une équivalence de catégorie.
\end{theo}
\begin{proof}
Comme $Y _0$ est la somme directe de ses composantes irréductibles,
il ne coûte rien de supposer $Y _0$ affine, lisse et intègre. 
 Grâce à Elkik (\cite{elkik}), il existe un
$\V$-schéma affine et lisse $Y$ dont la réduction modulo $\pi$ est isomorphe à $Y _0$.
Il existe alors une immersion fermée $Y  \hookrightarrow \A _\V ^r$.
On note $P := \P^r _\V$,
$U :=\A _\V ^r$, $X$ l'adhérence schématique de $Y$ dans $\P ^r _\V$
et $T:=P \setminus U$.
On obtient
une immersion fermée $Y ^\dag\hookrightarrow U ^\dag$ de $\V$-schémas formels faibles lisses.

Par \ref{dense-ideal}, il existe un diviseur $\widetilde{T} _0$ de $P _0$ contenant $T _0$ tel que,
en notant $\widetilde{Y}$ l'ouvert de $Y$ complémentaire de $ \widetilde{T} _0$,
$\widetilde{Y} _0$ soit affine et dense dans $Y _0$ et
$\widetilde{Y} _0 \subset X _0$ se désingularise idéalement.
Notons $\widetilde{U} := U \setminus \widetilde{T} _0$.
L'immersion fermée $Y \hookrightarrow U$ induit la suivante $\widetilde{Y} ^\dag \hookrightarrow \widetilde{U} ^\dag$.
Le triplet $( P ^\dag, \widetilde{T} _0, \widetilde{Y} ^\dag \hookrightarrow \widetilde{U} ^\dag)$
fournit donc un modèle idéal à $\widetilde{Y} _0$.
\hfill \hfill \qed \end{proof}

\section{$F$-complexes de $\mathcal{D}$-modules arithmétiques dévissables}

\subsection{Définitions et lien avec la surholonomie}
Dans cette section, $\PP$ est un $\V$-schéma formel propre et lisse et $T$ est un diviseur de $P$.

\begin{defi}\label{defi-dev}
Soient $\E$
un objet de $F\text{-}D ^\mathrm{b} _\mathrm{coh} (\D ^\dag _{\PP} (\hdag T) _\Q)$, $X$ son support
et $Y := X \setminus T$.
On dit que $\E$ {\og se dévisse en $F$-isocristaux surconvergents \fg}
s'il existe des diviseurs $T _1,\dots, T _{r+1}$ de $P$ contenant $T$ tels que
\begin{enumerate}
  \item $Y _1 = X \setminus T _1$ est affine, lisse et possède un modèle idéal (\ref{defmodide}) ;
  \item Pour $1 \leq i\leq r$, $Y _{i+1} := (X \cap T _1 \cap \dots \cap T _i) \setminus T _{i+1}$,
  est affine, lisse et possède un modèle idéal ;
  \item Soit $Y _{r+1} :=X \cap T _1 \cap \dots \cap T _{r+1}$ est lisse soit
  $Y _{r+1} \setminus T$ est vide ;
  \item[a)] Les espaces de cohomologie de $(\hdag T_1) (\E)$ sont associés
  à des $F$-isocristaux surconvergents sur $Y _1 $,
  i.e., sont des objets de
  $F\text{-}\mathrm{Isoc} ^{\dag \dag}( \PP, T _{1}, X /K)$ (voir \ref{coro-ideal->oK});
  \item[b)] Les espaces de cohomologie de
  $\R \underline{\Gamma} ^\dag _{T _1 \cap \dots \cap T _i} (\hdag T_{i+1}) (\E)$
  sont associés à des $F$-isocristaux surconvergents sur $Y _i$, i.e., appartiennent à la catégorie
  $F\text{-}\mathrm{Isoc} ^{\dag \dag}( \PP, T _{i+1}, X\cap T _1 \cap \dots \cap T _i/K)$ ;
  \item[c)] Les espaces de cohomologie de
  $\R \underline{\Gamma} ^\dag _{T _1 \cap \dots \cap T _{r+1}} (\E)$
  sont associés via \cite{caro_surholonome} à des $F$-isocristaux surconvergents
  sur $X \cap T _1 \cap \dots \cap T _{r+1}\setminus T$.
\end{enumerate}
En gros, on dispose d'une stratification $Y = \cup _{i=1,\dots r+1} Y _i$, telle que,
les espaces de cohomologie de la restriction de $\E$ au-dessus de chaque strate $Y _i$
soient associés à des $F$-isocristaux surconvergents sur $Y _i$.
On dira aussi que $\E$ est un $F\text{-}\D ^\dag _{\PP} (\hdag T) _\Q$-complexe {\it dévissable}
ou {\it se dévisse au-dessus de la stratification $Y = \cup _{i=1,\dots r+1} Y _i$}.
On notera $F\text{-}D ^\mathrm{b} _\mathrm{\text{dév}} (\D ^\dag _{\PP} (\hdag T) _\Q)$ la sous-catégorie pleine de
$F\text{-}D ^\mathrm{b} _\mathrm{coh} (\D ^\dag _{\PP} (\hdag T) _\Q)$ des
$F\text{-}\D ^\dag _{\PP} (\hdag T) _\Q$-complexes dévissables.
\end{defi}

\begin{prop}\label{finif+dev}
  Soit $( \E,\ \Phi )\in F\text{-}D ^\mathrm{b} _\mathrm{\text{dév}} (\D ^\dag _{\PP} (\hdag T) _\Q)$.
Les espaces de cohomologie de $f _{T+} (\E )$ sont des $K$-espaces vectoriels de dimension finie.
\end{prop}
\begin{proof}
  Cela découle par dévissage de \ref{prop-DRfrob} et de la finitude la cohomologie rigide
  \cite{kedlaya-finiteness}.
\hfill \hfill \qed \end{proof}

\begin{rema}\label{rema-caradev}
  Soit $\E \in D ^\mathrm{b} _\mathrm{coh} (\D ^\dag _{\PP} (\hdag T) _\Q)$.
  Les deux propriétés ci-après sont équivalentes :
  \begin{enumerate}
  \item[(a)]
  Pour tout diviseur $T'$ de $P$, $\DD_{\PP,T} (\hdag T') (\E)\in
D ^\mathrm{b} _\mathrm{surcoh} (\D ^\dag _{\PP} (\hdag T) _\Q)$ ;
  \item [(b)] 
  Pour tous sous-schémas fermés $Z$ et $Z'$ de $P$,
   $\DD_{\PP,T}\R \underline{\Gamma} ^\dag _Z (\hdag Z') (\E)$ est un objet de
   $D ^\mathrm{b} _\mathrm{surcoh} (\D ^\dag _{\PP} (\hdag T) _\Q)$ ;
  \end{enumerate}
  En effet, supposons que $\E$ satisfasse (a).
  En utilisant le dualisé du triangle de localisation en $T _1$ de
  $\R \underline{\Gamma} ^\dag _{Z _1}  (\hdag T _2) (\E)$ suivant
$$\R \underline{\Gamma} ^\dag _{T _1\cap Z _1}  (\hdag T _2) (\E)
\rightarrow
\R \underline{\Gamma} ^\dag _{Z _1}  (\hdag T _2) (\E)
\rightarrow
\R \underline{\Gamma} ^\dag _{Z _1}  (\hdag T _1 \cup T _2 ) (\E)
\rightarrow
+1 $$
valable pour tout sous-schéma fermé $Z _1$, et tous diviseurs $T _1 $ et $T _2$,
  on vérifie d'abord par récurrence sur le nombre minimal
  de diviseurs d'intersection $Z$ (on convient que celui-ci est nul lorsque $Z=X$),
  que, pour tout sous-schéma fermé $Z$ et tout diviseur $T'$,
  le complexe $\DD_{\PP,T}\R \underline{\Gamma} ^\dag _Z (\hdag T') (\E)$ appartient à
 $  D ^\mathrm{b} _\mathrm{surcoh} (\D ^\dag _{\PP} (\hdag T) _\Q)$.
  Ensuite, en appliquant le foncteur $\DD _{\PP, T}$ au triangle de localisation
  en $Z'$ de $\R \underline{\Gamma} ^\dag _Z (\E)$, on prouve que la propriété $(b)$ est validée pour $\E$.
  La réciproque est une tautologie.

  Il résulte des équivalences entre $(a)$ et $(b)$ que la propriété $(a)$ est préservée,
  pour tout sous-schéma fermé $Z$ de $P$, par les foncteurs
  $\R \underline{\Gamma} ^\dag _Z$ est $(\hdag Z)$.
\end{rema}

\begin{theo}\label{caradev}
Soit 
$\E \in F\text{-}D ^\mathrm{b} _\mathrm{coh} (\D ^\dag _{\PP} (\hdag T) _\Q)$.
On suppose que pour tout diviseur $T'$ de $P$, $\E$ et $\DD_{\PP,T} (\hdag T') (\E)$ sont dans
$F\text{-}D ^\mathrm{b} _\mathrm{surcoh} (\D ^\dag _{\PP} (\hdag T) _\Q)$.
Alors $\E$ se dévisse en $F$-isocristaux surconvergents.

\end{theo}
\begin{proof}
En reprenant les arguments du deuxième cas de la preuve \cite[3.2.4]{caro_surcoherent}, on vérifie
que le support de $\E$, noté $X$, est une compactification de $Y$.
On procède par récurrence sur la dimension de $X$.
Via \ref{densemodide},
il existe un diviseur
$T_1$ contenant $T$ tel que $Y _1:= X \setminus T _1$ soit affine, lisse et dense dans $X$ et possède un modèle idéal.
Notons $u$ : $\Y _1 \hookrightarrow \PP$ un relèvement de l'immersion $Y _1\hookrightarrow P$.
Comme le faisceau $u ^! (\E)$ étant $\D ^\dag_{\Y _1,\Q}$-surcohérent,
quitte à remplacer $T_1$ par un diviseur plus grand,
on peut supposer que les espaces de cohomologie de
$u ^! (\E)$ sont $\O _{\Y _1,\Q}$-cohérents (\cite[2.2]{caro_courbe}).
Il en résulte que les espaces de cohomologie de $(\hdag T _1) (\E) $
soient dans $ F\text{-}\mathrm{Isoc} ^{\dag \dag}( \PP, T _{1}, X /K)$.
Grâce à la remarque \ref{rema-caradev} et comme la dimension de $X \cap T _1$ est strictement inférieure
à celle de $X$,
par hypothèse de récurrence,
$\R \underline{\Gamma} ^\dag _{T _1} (\E)$ se dévisse en $F$-isocristaux surconvergents.
D'où le résultat.
\hfill \hfill \qed \end{proof}

Le théorème \ref{caradev} implique aussitôt l'inclusion ci-après.
\begin{theo}\label{surholincdev}
  On a l'inclusion $F\text{-}D ^\mathrm{b} _\mathrm{surhol} (\D ^\dag _{\PP, \Q})
  \subset F\text{-}D ^\mathrm{b} _\mathrm{\text{dév}} (\D ^\dag _{\PP, \Q})$.
\end{theo}

\begin{conj}\label{conjsp+surhol}
  Il est raisonnable de penser que l'inclusion \ref{surholincdev} est une égalité.
  Cela résulte aussitôt de la conjecture suivante :

  {\og Soit $Y $ un $k$-schéma affine, lisse et possédant un bon modèle.
  Pour tout $F$-isocristal surconvergent sur $Y$,
  $\sp _{Y+} (E)$ est un $F\text{-}\D _Y$-module arithmétique surholonome (voir \cite[1.2.7]{caro_surholonome}). \fg}

\end{conj}
\begin{rema}
Si la conjecture de Berthelot \cite[5.3.6.D]{Beintro2} sur la stabilité de l'holonomie
est vraie alors la conjecture \ref{conjsp+surhol} l'est aussi.
De plus, le théorème qui suit nous permet de penser que \ref{conjsp+surhol} est réaliste.
\end{rema}
\begin{theo}
  Soit $Y$ un $k$-schéma affine, lisse et possédant un bon modèle.
  Pour tout $F$-isocristal unité surconvergent sur $Y$,
  $\sp _{Y+} (E)$ est un $F\text{-}\D _Y$-module arithmétique surholonome.
\end{theo}
\begin{proof}
Via la remarque \ref{isosurcohpf0},
  on se ramène au cas où $Y$ est irréductible.
  Grâce à \ref{defidagdagpxt?}, on peut supposer que $Y$ possède une compactification
  lisse, ce qui a été traité dans \cite[2.3.2]{caro_unite}.
\hfill \hfill \qed \end{proof}

\begin{theo}
Soient $X$ un sous-schéma fermé de $P$ et $Y := X \setminus T$.
Si la conjecture \ref{conjsp+surhol} est validée,
alors la catégorie des $F$-complexes de $\D _Y$-modules arithmétiques surholonomes,
$F\text{-}D ^\mathrm{b} _\mathrm{surhol} (\D _{Y })$, est stable par produit tensoriel interne.
\end{theo}
\begin{proof}
On procède par récurrence sur la dimension de $X$. 
Soient $\E ,\E'$ deux $F$-complexes de
 $F\text{-}D ^\mathrm{b} _\mathrm{coh} (\D ^\dag _{\PP} (\hdag T) _\Q) \cap
 F\text{-}D ^\mathrm{b} _\mathrm{surhol} (\D ^\dag _{\PP, \Q})$ à supports dans $X$.
Via les mêmes arguments qu'au début de la preuve de \ref{caradev}, il existe
un diviseur
$T_1$ contenant $T$ tel que $Y _1:= X \setminus T _1$ soit affine, lisse et dense dans $X$ et possède un modèle idéal
et tel que
les espaces de cohomologie de $(\hdag T _1) (\E) $ et de $(\hdag T _1) (\E ') $
sont dans $ F\text{-}\mathrm{Isoc} ^{\dag \dag}( \PP, T _{1}, X /K)$.

Par \ref{stbintn}, on conclut par hypothèse de récurrence en utilisant les triangles de localisation
en $T _1$ de $\E$ et $\E'$.
\hfill \hfill \qed
\end{proof}

\begin{prop}
Soient $\E, \E' \in F\text{-}D ^\mathrm{b} _\mathrm{coh} (\D ^\dag _{\PP} (\hdag T) _\Q)$
et $\U := \PP \setminus T$.
On suppose que,
pour tout diviseur $T'$ de $P$, les faisceaux $\E$,
et $\DD_{\PP,T} (\hdag T') (\E)$ sont des objets de
$F\text{-}D ^\mathrm{b} _\mathrm{surcoh} (\D ^\dag _{\PP} (\hdag T) _\Q)$
et de même pour $\E'$.

Soit $\phi $ : $\E  |_{\U} \rightarrow \E ' |_{\U} $ un morphisme $\D ^\dag _{\U,\Q} $-linéaire. Il
existe alors un morphisme $\psi$ : $\E \rightarrow \E'$ tel que $\psi |_\U =\phi$.
\end{prop}
\begin{proof}
  On procède par récurrence sur la dimension de $X$, la réunion des supports de $\E $ et $\E'$.
  De manière analogue à la preuve de \ref{caradev},
  il existe un diviseur $T_1$ contenant $T$ tel que
  $Y _1:= X \setminus T _1$ soit affine, lisse et dense dans $X$ et possède un modèle idéal,
  et tel que les espaces de cohomologie de
  $( \hdag T _1) (\E ), ( \hdag T _1) (\E ')$
  soient des objets de $ F\text{-}\mathrm{Isoc} ^{\dag \dag}( \PP, T _{1}, X /K )$.
  Comme le foncteur canonique
  $F\text{-}\mathrm{Isoc} ^{\dag}( Y _1/K) \rightarrow
  F\text{-}\mathrm{Isoc} ^{\dag } ( Y _1, Y/K) $ est pleinement fidèle
  (\cite{kedlaya_full_faithfull}), avec les notations de \ref{plfidked},
  le foncteur canonique de restriction $|_{\U}$ :
  $F\text{-}\mathrm{Isoc} ^{\dag \dag}( \PP, T _{1} , X /K ) \rightarrow
  F\text{-}\mathrm{Isoc} ^{\dag}( \U, T _{1}\cap U , Y /K )$
  l'est aussi.
  Il en résulte (par dévissage en utilisant \cite[I.7.2]{HaRD})
  qu'il existe un morphisme  $\psi$ : $( \hdag T _1) (\E ) \rightarrow ( \hdag T _1) (\E ')$
  tel que $\psi |_\U = ( \hdag T _1) (\phi)$.
  Avec la remarque \ref{rema-caradev} et puisque $\dim X \cap T _1 < \dim X$,
  on conclut
par hypothèse de récurrence et
via les triangles de localisation en $T _1$ de $\E$ et $\E'$. 
\hfill \hfill \qed \end{proof}

\subsection{Application aux fonctions $L$}
\begin{vide}
  \label{notafctL}
  On suppose ici $k =\F _q$ où $q=p ^s$.
On se donne $\PP$ un $\V$-schéma formel propre et lisse, $T $ un diviseur de $ P$,
$\U$ l'ouvert de $\PP$ complémentaire de $T$ et $j $ : $\U \hookrightarrow \PP$ l'immersion ouverte correspondante.

Si $x$ est un point fermé de $P$, on désigne par $k(x)$ le corps résiduel de $x$, $\deg x$ son degré,
$i_x$ : $\S (x)=\Spf \V(x)  \hookrightarrow \PP $ un relèvement $\V$-linéaire de l'immersion fermée canonique
$\Spec k(x) \hookrightarrow P$, $K(x)$ le corps des fractions de $\V(x)$
et $f _x$ : $\S(x) \rightarrow \S$ le morphisme structural.
\end{vide}

La définition des fonctions $L$ (voir \cite{caro_courbe} ou \cite{caro_surcoherent})
s'étend aux duaux des $F$-complexes dévissables :
\begin{defi}
On note $F\text{-}D ^{\mathrm{b}*} _\mathrm{\text{dév}} (\D ^\dag _{\PP} (\hdag T) _\Q)$,
la sous-catégorie pleine de
$F\text{-}D ^{\mathrm{b}} _\mathrm{coh} (\D ^\dag _{\PP} (\hdag T) _\Q)$
des complexes $\E$ tels que $\DD _{T} (\E) \in
F\text{-}D ^{\mathrm{b}} _\mathrm{\text{dév}} (\D ^\dag _{\PP} (\hdag T) _\Q)$.

Soient $S$ un sous-ensemble de $Y$, $S^{\mathrm{0}}$ l'ensemble des points fermés de $S$ et
$\E\in F\text{-}D ^{\mathrm{b}*} _\mathrm{\text{dév}} (\D ^\dag _{\PP} (\hdag T) _\Q)$.
La fonction $L$ associée à $\E$ au dessus de $S$ est définie en posant :
$$L( S, \E, t )=\prod _{ x \in S^{\mathrm{0}} } \prod_{r \in \Z }
 \det _K \left( 1 - t^{\deg x } F ^{\deg x } _{| H^r ( f _{x+} i_{x,T}^+ (\E)) }\right)^{(-1)^{r+1+d _{X}} / \deg x }.$$

 De plus, sa fonction cohomologique $P$ s'écrit :
$$P(\U , \E, t ) := \prod _{r\in \Z} \det_K \left( 1- t F_{|H^r   (  f _{T,!} \E ) }\right)^{(-1)^{r+1+d _X}}.$$
On remarque que grâce à \ref{finif+dev}, cette fonction a bien un sens et est de surcroît une fraction
rationnelle.
\end{defi}

\begin{vide}
\label{notaFE}
Soit $X$ un sous-schéma fermé de $P$ tel que $Y := X \setminus T$ soit affine et lisse.
On suppose en outre que $Y \subset X$ se désingularise localement idéalement et on se donne
$E $, un $F$-isocristal surconvergent sur $Y$. 

Pour tout point fermé $x$ de $Y$,
on note $E _x := H ^{\mathrm{0}} _{\mathrm{rig}} ( \Spec k(x) , i ^* _{x, K} E )$ la fibre de $E$ en $x$ et
$F _{| E _x}$ son automorphisme de Frobenius.

Dans \cite[2.3]{E-LS}, la fonction $L$ associée à $E $ est donnée par
$$L(Y, E, t ):=\prod _{x \in Y ^0} \det _K ( 1 - t ^{\deg x }  F ^{\deg x}_{| E _x} ) ^{ -1 / \deg x}.$$
\end{vide}

\begin{lemm}
\label{L=L}
Avec les notations \ref{notaFE}, on a l'égalité :
$$L ( Y , E ^\vee, t ) =L (\Y , \DD _{\PP ,T} (\sp _{Y+} E ) , t).$$
\end{lemm}
\begin{proof}
Analogue à \cite[3.3.1]{caro_courbe}.
\hfill \hfill \qed \end{proof}

\begin{prop}\label{L=Pisoc}
Avec les notations \ref{notaFE}, on pose $\E :=\DD _{\PP ,T} (\sp _{Y+} E$).
L'égalité $L(\U , \E, t )= P(\U , \E, t )$ est satisfaite.
\end{prop}
\begin{proof}
  De manière analogue à \cite[3.3.3]{caro_courbe},
  cela découle de \ref{prop-DRfrob}, \ref{L=L} et de la formule cohomologique d'\'Etesse et Le Stum :
\begin{equation}
\notag 
L ( Y , E ^{\vee},t) =\prod _{r=0} ^{2 d _{\X}}
\det _K ( 1 -t q ^{d _{\X}}
F^{-1} _{| H ^r _{\mathrm{rig}} ( Y , E  )}) ^{(-1) ^{r+1}}. \hfill \hfill \qed
\end{equation}
\end{proof}

\begin{theo}
  Pour tout $\E \in F\text{-}D ^{\mathrm{b}*} _\mathrm{\text{dév}} (\D ^\dag _{\PP} (\hdag T) _\Q)$,
  la formule $L(\U , \E, t )= P(\U , \E, t )$
  est validée.
\end{theo}
\begin{proof}
  Cela découle par dévissage de \ref{L=Pisoc}.
\hfill \hfill \qed \end{proof}

\subsection{Analogue $p$-adique Weil II}
On reprend les notations et hypothèses de \ref{notafctL}.
De plus, $\iota$ : $K ^\mathrm{alg} \hookrightarrow \C$ désigne un plongement d'une clôture algébrique
de $K$ dans le corps des complexes.
La proposition qui suit est un cas particulier du théorème \cite[6.6.2.(b)]{kedlaya-weilII} de Kedlaya
(en effet {\og $\iota$-pure \fg} implique {\og $\iota$-realizable \fg}).
\begin{prop}[Kedlaya]\label{kedweilII}
Soient $Y$ un $k$-schéma lisse et
$E$ un $F$-iso-\linebreak cristal surconvergent sur $Y$. Si $E$ est de poids $\iota $-pur $w$, alors,
pour tout $r \in \Z$, $H ^r _{\mathrm{rig}} (Y, E/K)$ est de poids $\iota$-mixte $\geq w + r$.
\end{prop}

La notion de poids pour les $F$-isocristaux surconvergents s'étend aux complexes
à fibres extraordinaires finies (voir la définition dans \cite{caro_hol-unite}) :
\begin{defi}
  Soient $\E \in F\text{-}D ^\mathrm{b} _\mathrm{coh} (\D ^\dag _{\PP} (\hdag T) _\Q)$ à fibres extraordinaires finies
  et $x$ un point fermé de $P$.
  \begin{enumerate}
  \item $\E$ est de poids $\iota$-pur $w$ en $x$ si, pour tout $r \in \Z$,
$\mathcal{H} ^r( \DD  i ^! _x \DD ^* _T (\E))$ est de poids $\iota$-pur $w +r$
(au sens des $F$-isocristaux surconvergents) ;
  \item $\E$ est de poids $\iota$-mixte $\geq w$ en $x$ si, pour tout $r \in \Z$,
$\mathcal{H} ^r( \DD  i ^! _x \DD ^* _T (\E))$ est de poids $\iota$-mixte $\geq w +r$.
\end{enumerate}
De plus, $\E$ est de poids $\iota$-pur $w$ (resp. de poids $\iota$-mixte $\geq w$) s'il l'est en tout
point fermé de $P$.
\end{defi}

Dans le langage des $\D$-modules arithmétiques, la proposition \ref{kedweilII} se traduit par la suivante.
\begin{prop}\label{kedweilIID}
  Avec les notations \ref{prop-DRfrob},
  pour tout $F$-isocristal surconvergent $E$ sur $Y  _0$,
  de poids $\iota$-pur $w$, $\E := \sp _{Y ^{\dag}  \hookrightarrow U ^{\dag}, T _0 , +} (E)$
  est de poids $\iota$-pur $w+d _Y$. De plus,
  $f _+ (\E)$ est de poids $\iota$-mixte $w+d _Y$.
\end{prop}
\begin{proof}
  Soient $x$ un point fermé de $Y$ et $u _x $ : $\Spec k (x) \hookrightarrow Y _0$ l'immersion fermée
  canonique induite. Grâce à l'analogue $p$-adique de Ber-\linebreak thelot du théorème de Kashiwara, avec les théorèmes
  de dualité relative et de bidualité, et via l'isomorphisme
  de comparaison des foncteurs duaux d'un $F$-isocristal surconvergent(voir \cite{caro_comparaison}, on vérifie
  l'isomorphisme canonique $\DD  i ^! _x \DD ^* _T (\E) \riso u ^* _x (E) [d _Y]$
  compatible à Frobenius. Ainsi, $\E$ est de poids $\iota$-pur $w+d _Y$.
  La dernière assertion découle ensuite de \ref{prop-DRfrob} (notations \ref{prolong}).
  \hfill \hfill  \qed
\end{proof}

Le théorème qui suit est un analogue $p$-adique de Weil II ponctuel
(on ne s'intéresse pas ici au cas relatif).

\begin{theo}
 Soit $\E$ un complexe de $ F\text{-}D ^{\mathrm{b}} _\mathrm{\text{dév}} (\D ^\dag _{\PP} (\hdag T) _\Q)$.
  Si $\E$ est de poids $\iota$-pur $w$
  alors, $ f _{T +} (\E)$ est de poids $\iota$-mixte $\geq w$.
  \end{theo}
\begin{proof}
  Cela résulte par dévissage de \ref{kedweilII} et de \ref{kedweilIID}.
\hfill \hfill  \qed
\end{proof}

\begin{acknowledgement}
  L'auteur a bénéficié du soutien de l'Université de Sydney ainsi que de celui
du réseau européen TMR \textit{Arithmetic Algebraic Geometry}
(contrat numéro UE MRTN-CT-2003-504917).
\end{acknowledgement}

\bibliographystyle{smfalpha}
\bibliography{bib1}

\end{document}